\documentclass[11pt]{amsart}
\pdfoutput=1
\usepackage{latexsym, graphicx, epsfig, amsmath, amsfonts,amssymb}
\usepackage{multirow}
 \usepackage{diagbox}
\usepackage{color}
 \usepackage[percent]{overpic}
\usepackage{algorithm}
\usepackage{algorithmicx}
\usepackage{algpseudocode}
\floatname{algorithm}{Algorithm}
\usepackage{threeparttable}
\usepackage{booktabs}

\usepackage{epstopdf,mathtools,bbm}
\usepackage[left=1.25in,right=1.25in,top=1in]{geometry}

\def\mb{\mathbf}
\def\R{\mathbb{R}}

\newcommand{\bs}[1]{\boldsymbol{#1}}

\DeclareMathOperator*{\argmax}{arg\,max}

\newcounter{Rownumber}

\title{Optimal design for  kernel  interpolation:  applications to uncertainty quantification}

\author{Akil Narayan}
\thanks{Department of Mathematics and Scientific Computing and Imaging Institute,
	University of Utah, Salt Lake City, UT 84112.  Email: akil@sci.utah.edu.  A. Narayan is partially supported by NSF DMS-1848508 and AFOSR FA9550-20-1-0338. }

\author{Liang Yan}
\thanks{School of Mathematics, Southeast University, Nanjing, China. Email: yanliang@seu.edu.cn. L. Yan is supported by NSF of China (No.11771081), the science challenge project (No. TZ2018001) and  Zhishan Young Scholar Program of SEU}

\author{Tao Zhou}
\thanks{LSEC, Institute of Computational Mathematics, Academy of Mathematics and Systems
Science, Chinese Academy of Sciences, Beijing 100190, China. Email: tzhou@lsec.cc.ac.cn. T. Zhou is partially supported  by the National Key R$\&$D Program of China(No. 2020YFA0712000), the NSF of China (under grant numbers 11822111, 11688101and 11731006), the science challenge project (No. TZ2018001),  the Strategic Priority Research Program of Chinese Academy of Sciences (No.  XDA25000404) and youth innovation promotion association (CAS)}

\begin{document}

\maketitle

\begin{abstract}
The paper is concerned with classic kernel interpolation methods, in addition to approximation methods that are augmented by gradient measurements. To apply kernel interpolation using radial basis functions (RBFs) in a stable way, we propose a type of quasi-optimal interpolation points, searching from a large set of \textit{candidate} points, using a procedure similar to designing Fekete points or power function maximizing points that use pivot from a Cholesky decomposition. The proposed quasi-optimal points results in smaller condition number, and thus mitigates the instability of the interpolation procedure when the number of points becomes large. Applications to parametric uncertainty quantification are presented, and it is shown that the proposed interpolation method can outperform sparse grid methods in many interesting cases. We also demonstrate the new procedure can be applied to constructing gradient-enhanced Gaussian process emulators.
\end{abstract}


\pagestyle{myheadings}
\thispagestyle{plain}
\markboth{AKIL NARAYAN, LIANG YAN, AND TAO ZHOU}
{OPTIMAL DESIGN FOR KERNEL INTERPOLATION: APPLICATIONS
TO UQ}

\maketitle

\section{Introduction}
The field of uncertainty quantification (UQ) has received intensive attention in recent  years as theoreticians and practitioners have tackled problems in the diverse areas of stochastic analysis, high-dimensional approximation, and Bayesian learning.   A fundamental problems in UQ is to approximation   of a multivariate function $u(Z)$, where the parameters $Z=(Z^{(1)},\cdots, Z^{(d)})$ is a $d$-dimensional random vector. The function $u$ might be a solution resulting from a stochastic PDE problem or a derived quantities of interest (QoI) from such a system.  The most popular approach is to expand the function $u$ in a generalized polynomial chaos (gPC) basis \cite{Ghanem_book_1991,Xiu_2002wiener} and approximate its coefficients with a Galerkin projection. In this case, a serious drawback is that existing deterministic solvers must be rewritten.

In recent years, more attention has been devoted to the construction of a surrogate of $u$ based on data $D_N \coloneqq \{(z_j, u(z_j))\}_{j=1}^N.$ This in general refers to stochastic collocation methods, as the sample points $\{z_j\}$ are chosen from the parametric domain. The procedure is also known as non-intrusive response construction.``Non-intrusive" effectively means that existing black-box tools can be used in their current form. Most existing work concerns the collocation based polynomial approximations, this include the interpolation technic based on sparse grid \cite{Zabaras_2007sparsegrid,narayan_adaptive_2014,Fabio,Xiu_2005highorder}, the least-squares projection onto polynomial spaces \cite{Cohen_interpolation,Guo2018weighted,Narayan_2014Christoffel,GNZ_2020}, the compressive sampling method with $\ell_1$ minimization \cite{Doostan_2011nonada,JNZ_2016,CMR_2020,Yan_GX_IJUQ_2012}.  Another approach to construct the surrogate using $D_N$ is the so called Gaussian process (GP)  regression \cite{Bilionis+Zabaras2016,Bilionis+Zabaras2012,Sacks1989,Tripathy+Bilionis2016}. The GP provides an analytically tractable Bayesian framework where prior information about the $f$ can be encoded in the covariance function, and the uncertainty about the prediction is easily quantified given measurements.  
One drawback of GP is the  number of samples required for an accurate surrogate increased exponentially as the number of input parameters grows. One method of mitigating this limitation is the incorporation of gradient information into the training of the surrogate \cite{deBaar2014,Lockwood2012,Morris1993,Ulaganathan2016}.  By incorporating derivative values, the cost associated with training an accurate surrogate can be greatly reduced.  Unfortunately, gradient-enhanced GP emulators raise more computational challenges compared to ordinary GP emulators \cite{He+Chien2018}.  Especially,  the condition number of the gradient-enhanced GP emulator will increase much faster than the original GP emulator.

In this work, we approximate $u(Z)$ using a kernel interpolation  through a limited number of support points. Over the last four decades, kernel methods using radial basis functions (RBFs) have been successfully applied to scattered data interpolation/approximation in hight dimensions  \cite{Fasshauer2007, Wendland2005}. It is known that the kernel  interpolation process is equivalent  to finding the unbiased estimator for a Gaussian process with covariance from measurements \cite{Scheuerer+Schaback2013,Schaback2006kernel}. The theoretical connections between kernel methods and Gaussian process regressions for the classical interpolations has recently been highlighted \cite{Fasshauer2015kernel,Scheuerer+Schaback2013}.  Similar to GP progression, the numerical performance of kernel interpolation  is depended on the number of sample points $N$ and the  so-called shape parameter $\epsilon$. If $\epsilon$ is too large, the interpolant or solution will usually be inaccurate. On the other hand, if $\epsilon$ is too  small, the condition number of the resultant matrix system will become so bad that linear solvers may perform poorly.  Searching for the best shape parameter is still an open problem. On the other hand, for some fixed $\epsilon$, the problem of ill-conditioning can still arise either if the total number of points (commonly called {\it centers}) is too large or if some points are closely clustered.  How to select the optimal  centers in the kernel methods is still a challenging work. In this paper, we consider a more practical and solvable problem:
\vskip 5pt
\begin{center}
  \textit{Given a large set of candidate points, how can we choose an ``optimal" subset of samples for interpolative kernel approximation?}
\end{center}
\vskip 5pt
To this end, we borrow the notion of \textit{Fekete points} in the polynomial interpolation community.  The selection strategy we propose only involves linear algebra, and essentially amounts to Cholesky-type decompositions for the design matrix. More precisely, the strategy is data-independent, so it can proceed offline before data is collected. We shown that the selected quasi-optimal points  results in a smaller condition number, and thus postpone dramatically the instability of the interpolation procedure when the number of points goes to large.   Comparisons with other sampling strategies are presented, and it is noticed that our method demonstrates much improved stability property.   Applications to parametric uncertainty quantification are also presented, and it is shown that the proposed kernel interpolation method over perform the sparse grid methods in many interesting cases.

The remained of this paper is organized as follows. Section \ref{sec:setup} summaries the problem formulation, where we start with an abstract problem and provide detains on the kernel interpolation. Our approach for choosing the optimal sub-set of samples is introduced in Section \ref{sec:method}. Several numerical examples are considered in Section \ref{sec:tests}  to study the empirical performance of the new approach. The conclusions are drawn in the final section.

\section {Problem setup and kernel interpolation} \label{sec:setup}
Here we consider a general setting for PDEs with random input parameters.  Let $Z=(Z^{(1)},\cdots,Z^{(d)}) \in I_Z \subseteq \R^d, d\geq 1$ be a parameter domain representing the uncertain inputs of the system. Consider
\begin{eqnarray}
\begin{cases}
   u_t(x,t,Z) = \mathcal{L} (u),   &  \mathcal{D}\times (0, T] \times I_Z \\
   \mathcal{B} (u)=0,  &  \partial{\mathcal{D}}\times  (0, T] \times I_Z  \\
  u=u_0 &    \mathcal{D}\times \{t=0\} \times I_Z
\end{cases}
\end{eqnarray}
where $\mathcal{D} \in \R^l, l=1,2, 3$, is the physical domain, and $T>0$ is the terminal time. Here  $\mathcal{L}$ is a differential operator and $\mathcal{B}$ a boundary operator that may involve differential operators with respect to the spacial variable $x.$  In a probabilistic framework, $Z$ is modeled as a random variable on a complete probability space $(\Omega, \mathcal{F}, \mathcal{P})$ equipped with the probability  distribution function $F_Z(z)=P(Z\leq z)$, with $z\in \R^d$.  Note that $Z$ reflects the complete parameterization of uncertain inputs and can include both physical parameters of the system and hyperparameters that characterize certain input processes.

Throughout the paper, we assume that conditioned on the $i$th independent sample of $Z$, denoted by $z_i$, a numerical solution to this problem may be identified by a fixed deterministic solver, such as finite element solves, finite difference solvers, ect.  Let $\Xi =\{z_1,\cdots, z_N\}\subset I_Z, N\geq 1$, be a set of sample points in $I_Z$ .  For any fixed $x \in \mathcal{D}$ and $t>0$,  we shall denote $u_j =u(x,t, z_j)$ for simplicity. Hereafter we will suppress the notions of $x$ and $t$ whenever possible, with the understanding that our statements are made for all fixed $x$ and $t$. Once the pairings $(z_j, u_j), j=1,\cdots, N$,  are obtained, the objective is to construct a function $u_N(Z)$ such that $u_N(Z)\approx u(Z)$ in a proper sense.

\subsection{Kernel interpolation in parameter space}
In the current work, we shall consider the classic kernel interpolations using RBF to construct $u_N$ based on unstructured meshes. The pairing information will be enforced exactly by requiring $u_N(z_j)=u_j$ for all $j=1,\cdots,N$. In the kernel interpolation framework,  we first pick a translation-invariant radial kernel $K=\Phi(\epsilon\|\cdot-\cdot\|): \R^d\times \R^d \rightarrow \R$.  Here, $\Phi: \R \rightarrow \R$ known as the radial basis functions (RBF), $\epsilon$ is a shape parameter and $\|\cdot\|$ is, for example, Euclidean distance.  Using the set  $\Xi =\{z_j\}^N_{j=1}\subset I_Z$, which is usually referred to as the {\it {trial centers}} in kernel interpolation, we can define a finite-dimensional trial space in the form of
\begin{eqnarray}\label{trspace}
\mathcal{U}_{\Xi}=\mathcal{U}_{\Xi,  K}:= \textmd{span}\{\Phi(\epsilon\|\cdot-z_j\|)\vert z_j \in \Xi\}.
\end{eqnarray}
In applications, it is trivial that having a suitable trial space is essential for the performance of kernel methods. A good trial space should contain some good approximation to the solution.  In the present work we shall confine ourselves to the case of positive definite kernels/RBFs.  The important role of positive definiteness for kernel interpolation was pointed out in \cite{Micchelli1986} and, for example, guarantees existence and uniqueness of the approximation. Widely used positive definite kernels/RBFs are the Gaussians, with $\Phi(r)=\exp(-r^2)$, the inverse multiquadrics (IMQ) with $\Phi(r)=1/\sqrt{1+r^2}$, and the family of compactly supported (CS) RBFs \cite{Wendland1995,Wendland2005}.

Now, let us consider kernel interpolation problems. Note that any  trial function $u_N$ is a linear combination of the basis used in defining \eqref{trspace} and is in the form of
 \begin{eqnarray} \label{interplant}
u_N(Z)=\sum^N_{j=1} c_{j}\, \Phi (\epsilon\|Z-z_j\|),
 \end{eqnarray}
for some coefficients $\mb{c}=[c_1,\cdots, c_N]^T \in \R^N$.
Using the interpolation conditions $u_N(z_i)=u_i, i=1,\cdots, N$, we obtain the following linear system
\begin{eqnarray}\label{mateq}
\begin{bmatrix}
\Phi(\epsilon\|z_1-z_1\|) & \cdots & \Phi(\epsilon\|z_1-z_N\|) \\
\vdots &  \ddots & \vdots\\
\Phi(\epsilon\|z_N-z_1\|) & \cdots & \Phi(\epsilon\|z_N-z_N\|)
\end{bmatrix}
\begin{bmatrix} c_1\\ \vdots \\ c_N \end{bmatrix}
=\begin{bmatrix} u_1\\ \vdots \\ u_N \end{bmatrix}  \quad \mbox{or} \quad
\mb{Ac =u},
\end{eqnarray}
where the matrix $\mb{A}=K(\Xi,\Xi) $ is symmetric with entries $\mb{A}_{ij}=\Phi(\epsilon\|z_i-z_j\|)$ for $z_i,z_j \in \Xi$. The coefficients $\mb{c}$ are unique, if the interpolation matrix $\mb{A}$ is invertible.

\textbf{Remark 1:} It is known that the kernel interpolation process in trial space \eqref{trspace} is equivalent to finding the unbiased estimator for a Gaussian process with covariance $K$ from realization $\mb{u}$ at $\Xi$, see  \ref{app: GP}.

\textbf{Remark 2:}  If two trial centers $z_i,z_j \in \Xi$ are too closed to each other, identically shaped trial basis functions centered at these centers have nearly equal values at $\Xi$. This leads to two columns of nearly identical values and the problem of ill-conditioning. Yet, we need to decide what the experimental configuration $\Xi$ is in order to effectively use the kernel interpolation.

\subsection{Leave-One-Out Cross Validation (LOOCV)}
Noted that the accuracy of the solution of Eq. (\ref{mateq}) and the well-conditioning of the matrix $\mb{A}$ depends on the shape parameter $\epsilon$.  Much work addressese this issue, see, e.g.  \cite{Fasshauer2007,Fasshauer+Zhang2007,Rippa1999,Scheuerer2011,Yang2018}, and the references therein. It is out of our scope of this paper to thoroughly analyze choosing optimal shape parameters.   In this work, we use the  leave-one out cross validation (LOOCV) algorithm \cite{Fasshauer+Zhang2007,Rippa1999} to  select the value of the shape parameter.
In the RBF interpolation setting the LOOCV algorithm targets solution of the underlying linear system $\mb{Ac}=\mb{u}$, where the entries of the matrix $\mb{A}$ depend on the shape parameter $\epsilon$ which we seek to optimize. The cost function defined as the norm of the error vector $e(\epsilon)$ has entries \cite{Fasshauer+Zhang2007,Rippa1999},
\begin{equation*}
e_i(\epsilon)=\frac{c_i}{\mb{A}^{-1}_{ii}},
\end{equation*}
where $c_i$ is the $i$th component of the vector $\mb{c}$ and $\mb{A}^{-1}_{ii} \coloneqq \left(\mb{A}^{-1}\right)_{ii}$ is the $i$th diagonal element of the inverse of the coefficient matrix.  Thus, the optimal value of the shape parameter is considered as the one which minimizes the cost function $e(\epsilon)$. We define the optimal shape parameter as the value $\epsilon^*$ such that
\begin{equation}\label{errorvector}
\|e(\epsilon^*)\|=\min_{\epsilon} \|e(\epsilon)\|.
\end{equation}

\subsection{A Gradient enhanced approach}
We consider inclusion of gradient measurements in the kernel interpolation. Consider the availability of the following data:
$\{z_i, u(z_i)\}^N_{i=1}$ and $\{z_i, u'_m(z_i)\}_{i=1}^N$ for all $m=1,2,\ldots, d)$. Here $u'_m(z)=\frac{\partial f(z)}{\partial z^{(m)}}$ stands for the partial derivative with respect to the m-th variable $z^{(m)}$.  We try to find an  interpolant of the form
\begin{equation*}
u_N(Z) =\sum^N_{j=1} c_j \Phi (\epsilon \|Z-z_j\|) -\sum^d_{m=1} \sum^N_{j=1} \beta_{m,j} \Phi'_m (\epsilon \|Z-z_j\|)
\end{equation*}
with appropriate radial basis functions $\Phi$ so that $u_N$ satisfies the generalized interpolation conditions
$$ \lambda_i u= \lambda_i  u_N, \, i=1,\ldots,N(d+1).$$ Here, $\lambda_i$ could denote point evaluation at the point $z_i$, or it could denote evaluation of some derivative at the point $z_i$.

After enforcing the interpolation conditions the system matrix is given by
\begin{eqnarray*}
\mb{B}=
\begin{bmatrix}
\mb{A}_{0,0}& \mb{A}_{0,1}  & \cdots & \mb{A}_{0,d} \\
\mb{A}_{1,0} & \mb{A}_{1,1} & \cdots & \mb{A}_{1,d} \\
\vdots & \vdots & \ddots & \vdots\\
\mb{A}_{d,0}  & \mb{A}_{d,1} & \cdots & \mb{A}_{d,d} \\
\end{bmatrix}
\end{eqnarray*}
where $\mb{A}_{m,n}$ is the $N\times N$ matrix with the $(i,j)$th element
\begin{eqnarray}
  \left(\bs{A}_{m,n}\right)_{i,j} =
  \begin{cases}
  \Phi(\epsilon \|z_i-z_j\|),   &  m=n=0, \\
  -\Phi'_n(\epsilon \|z_i-z_j\|),  & m=0, n\neq 0,  \\
  \Phi'_m(\epsilon \|z_i-z_j\|),  & m\neq0, n=0,  \\
   -\Phi''_{m,n}(\epsilon \|z_i-z_j\|),  & m\neq0, n\neq 0,
\end{cases}
\end{eqnarray}
and $\mb{B}$ is a symmetric $(d+1)N \times (d+1)N$ matrix.

\begin{figure}[htbp]
\begin{center}
    \includegraphics[width=4.2cm]{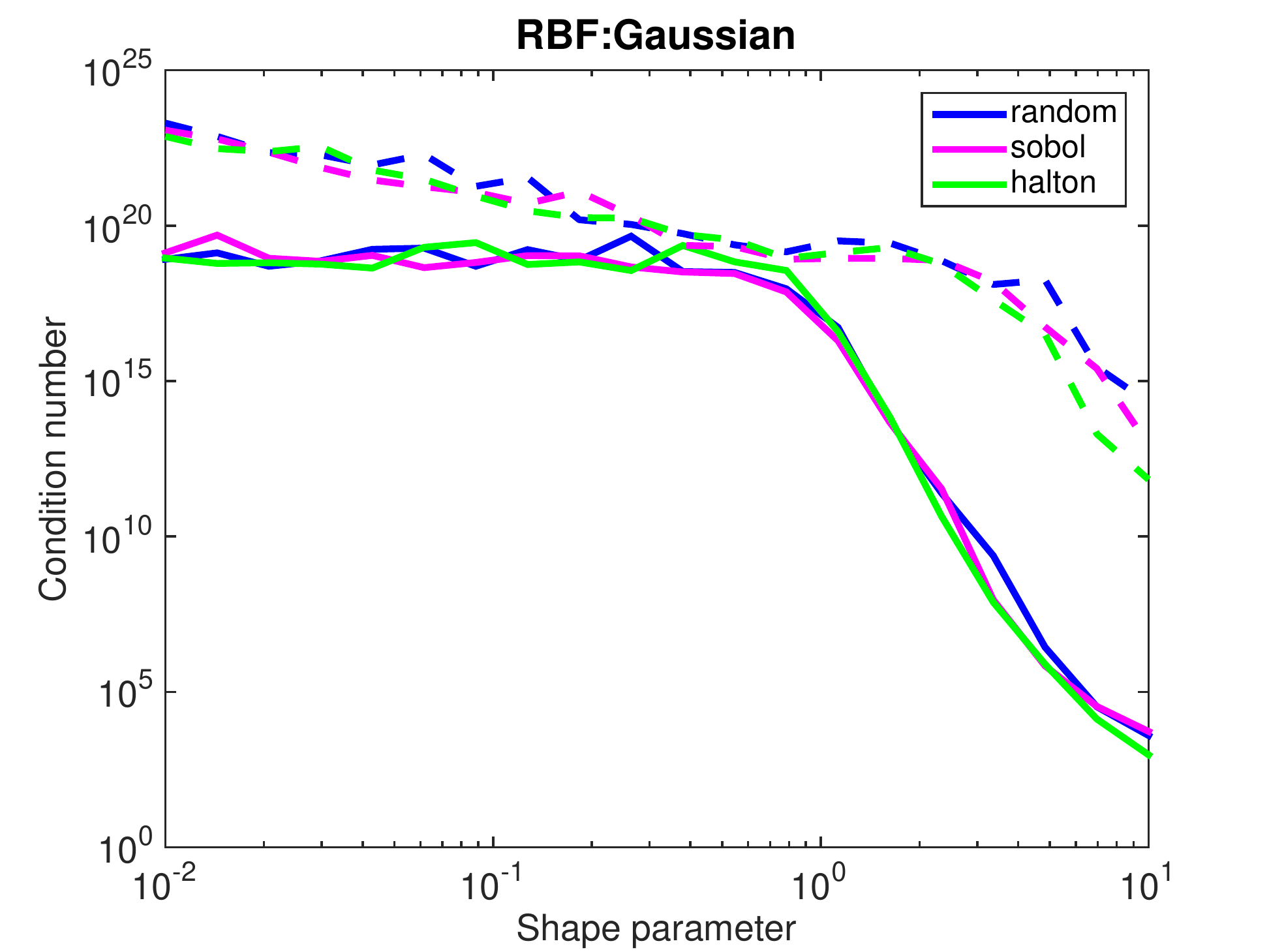}
    \includegraphics[width=4.2cm]{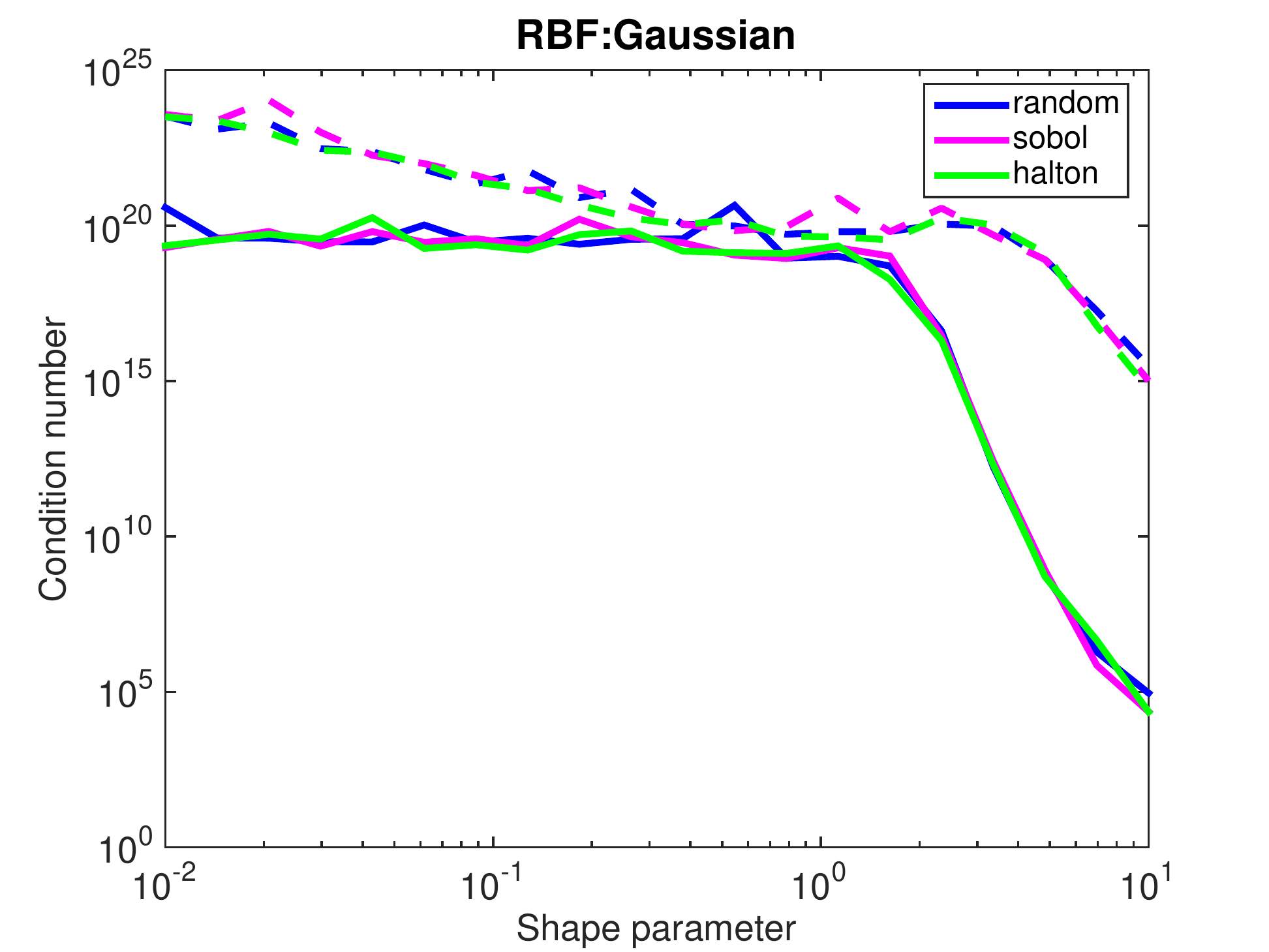}
    \includegraphics[width=4.2cm]{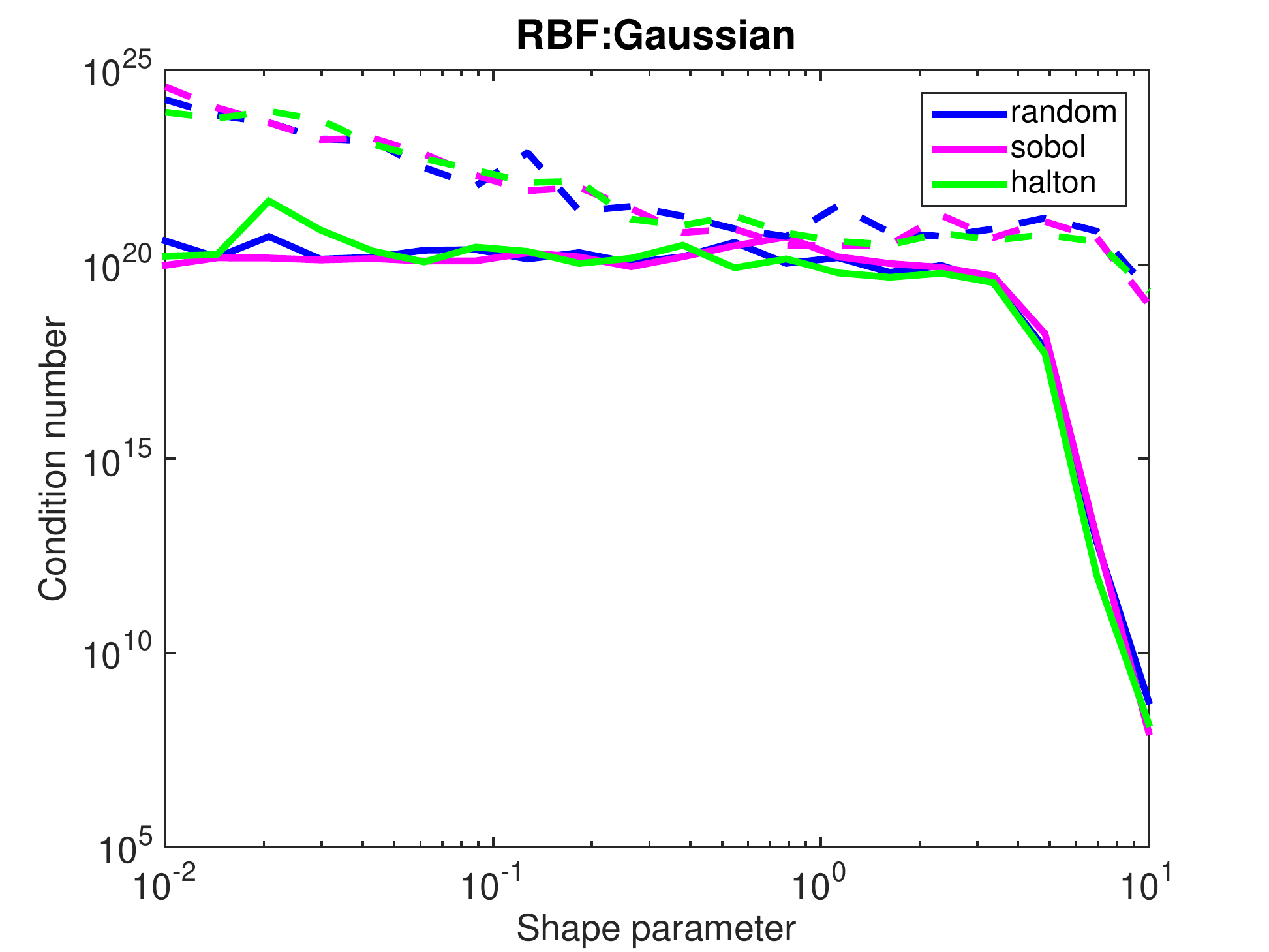}
     \includegraphics[width=4.2cm]{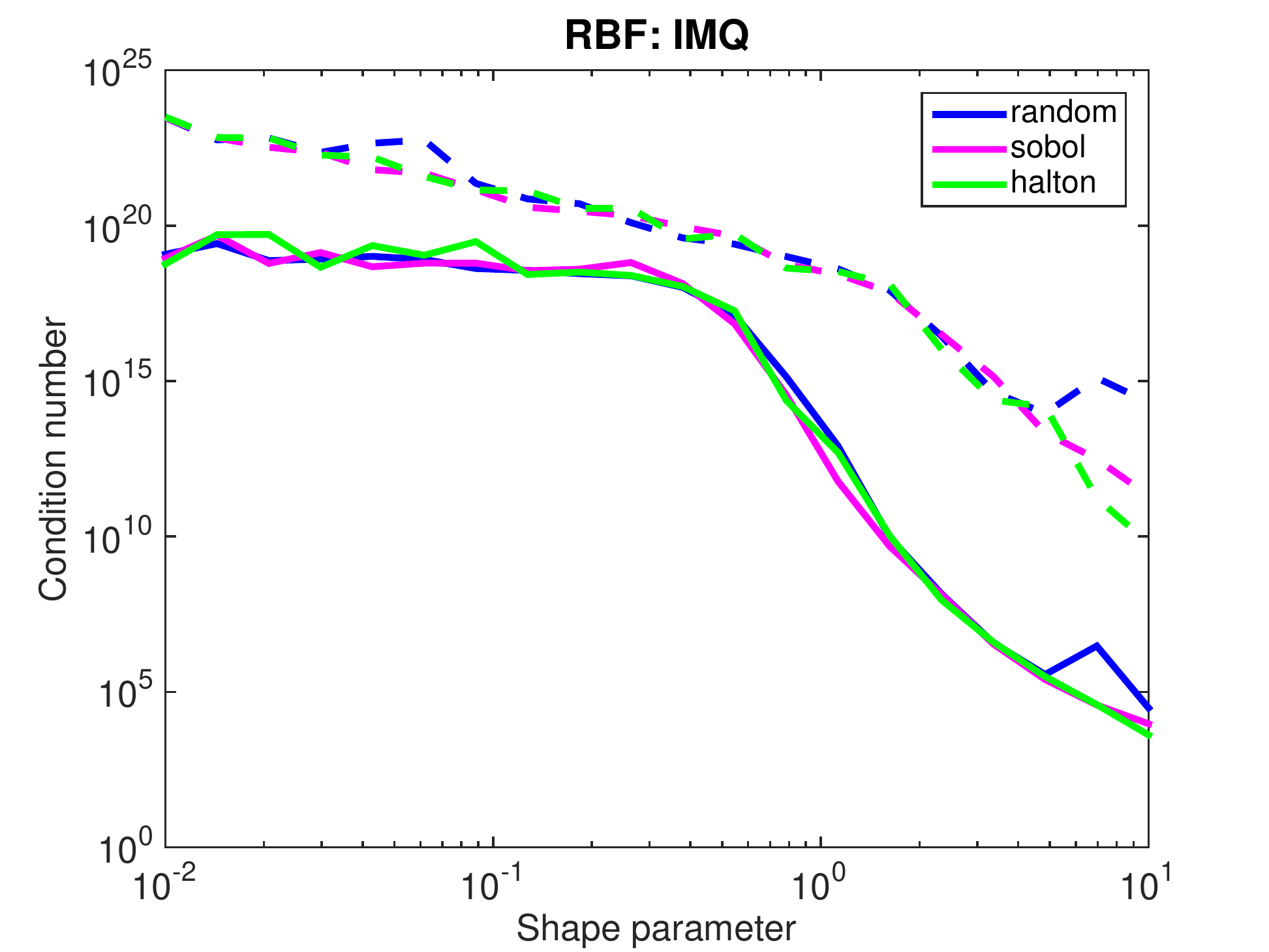}
   \includegraphics[width=4.2cm]{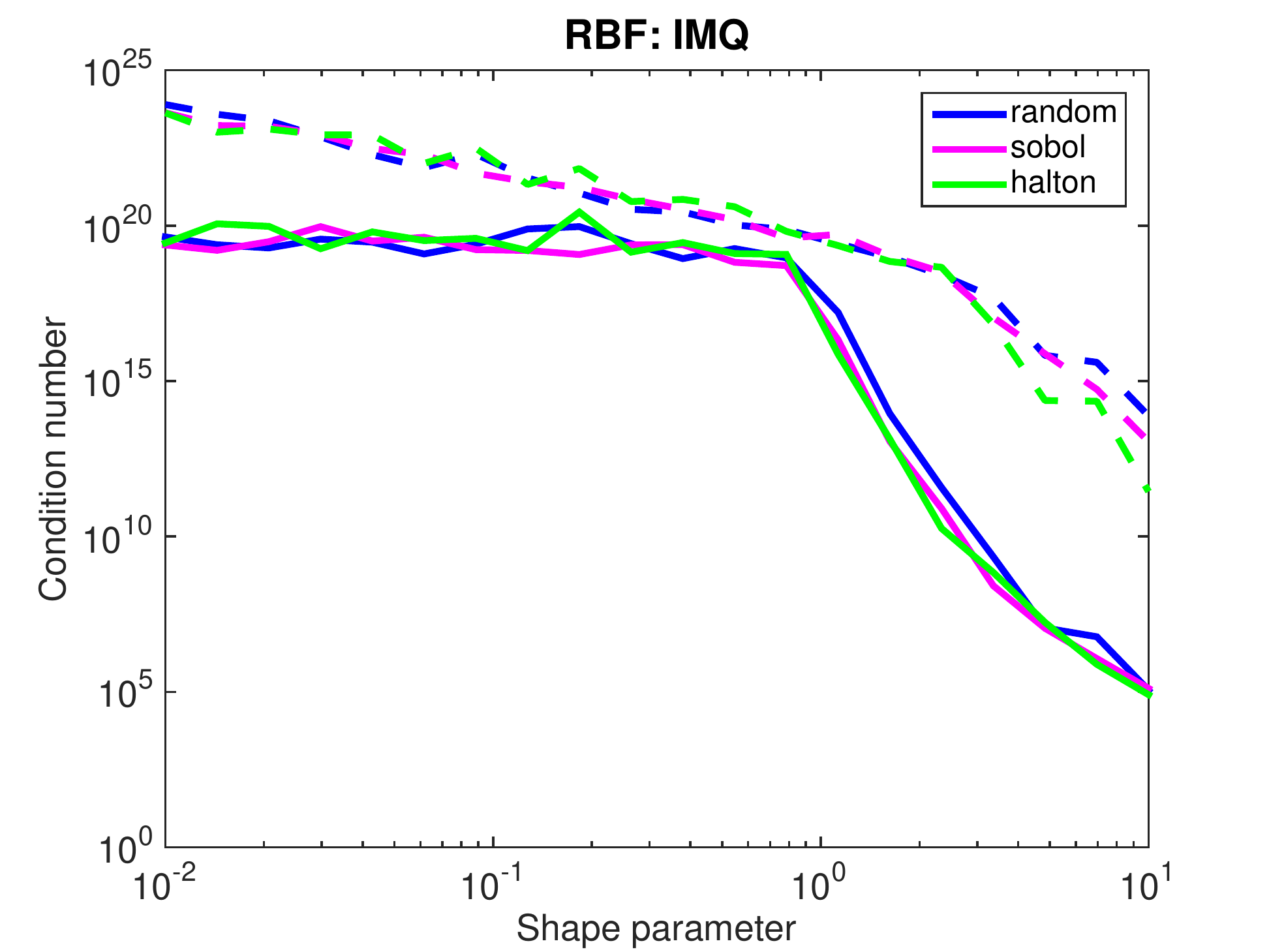}
   \includegraphics[width=4.2cm]{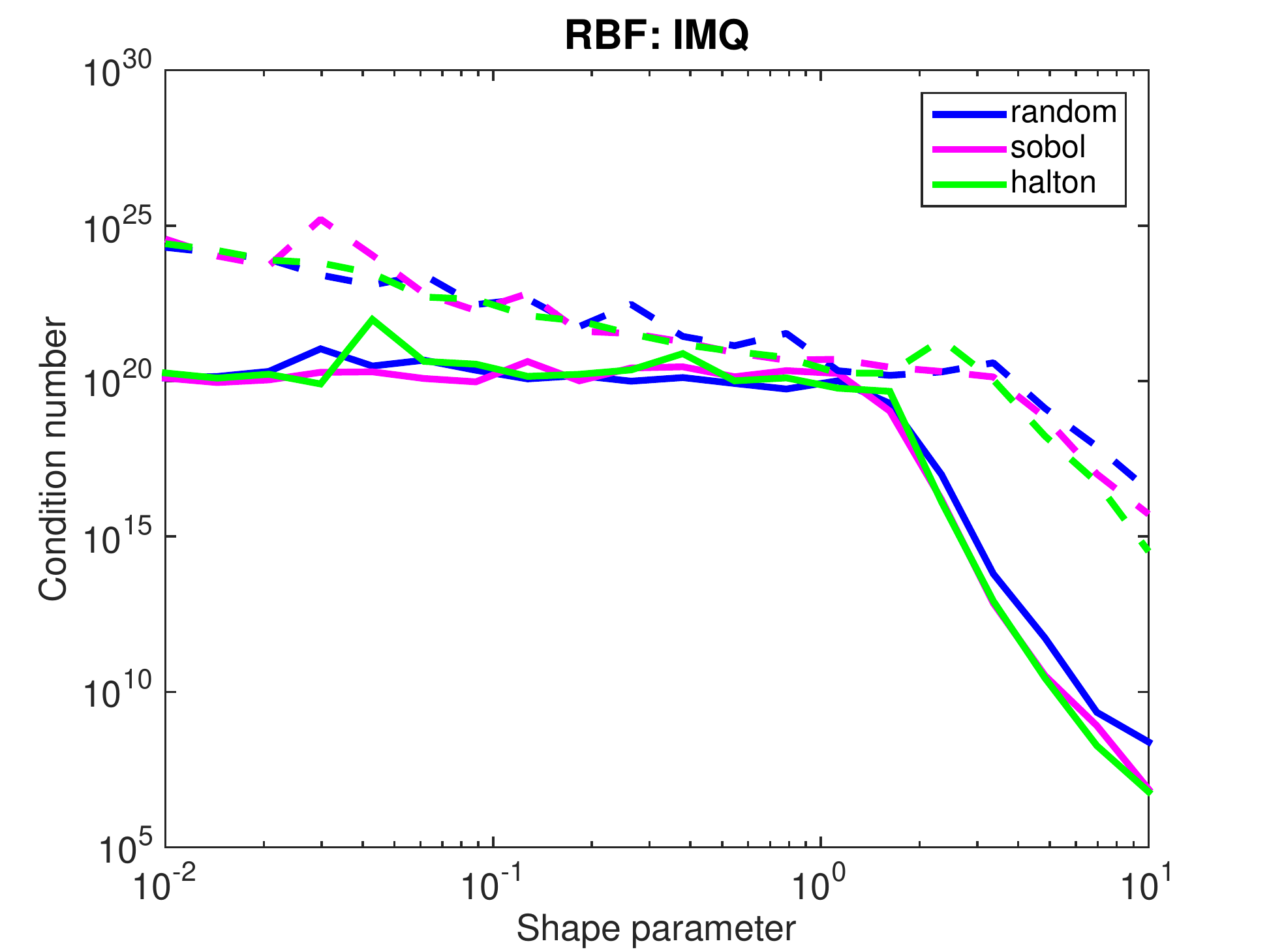}
        \includegraphics[width=4.2cm]{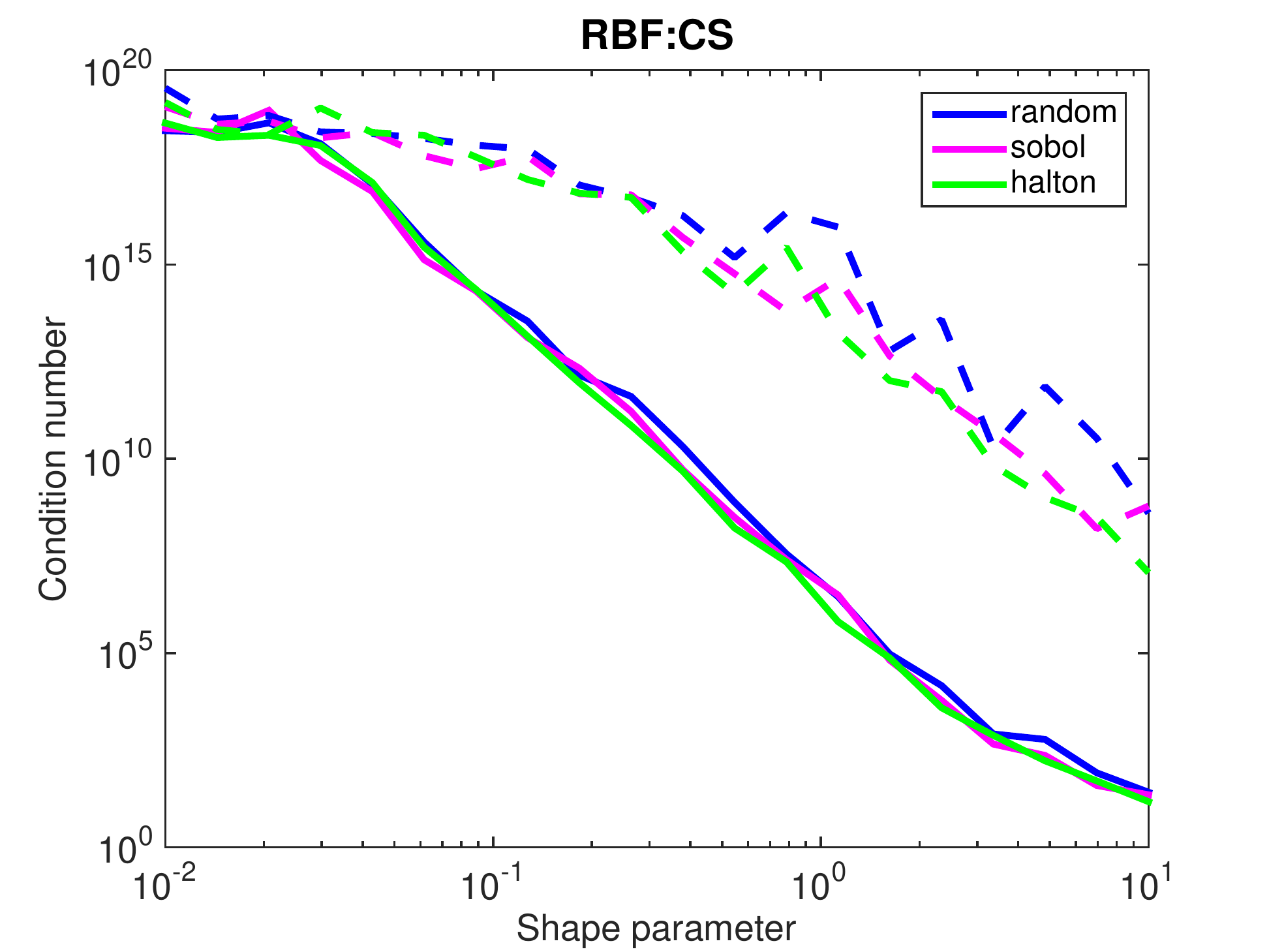}
   \includegraphics[width=4.2cm]{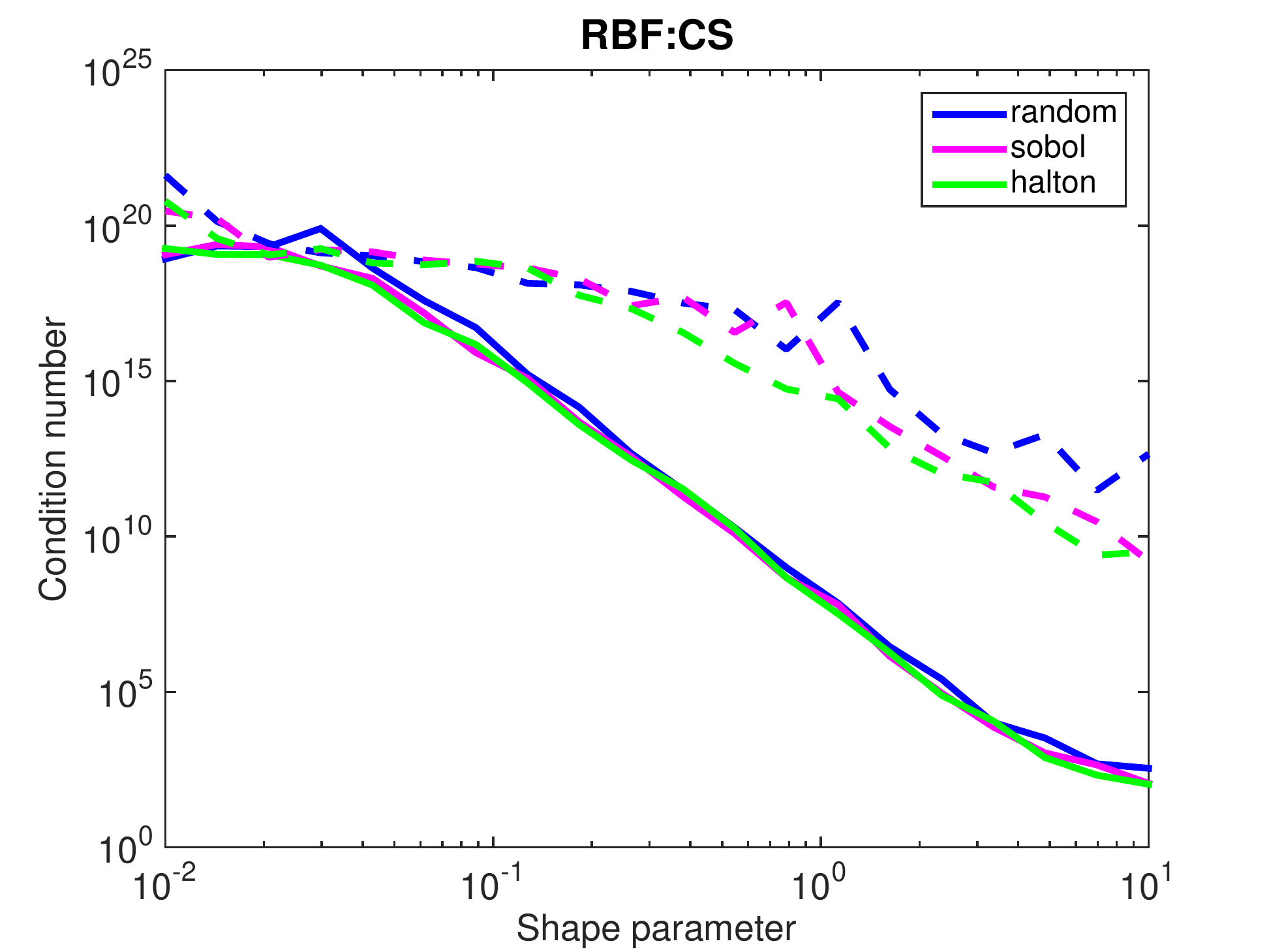}
   \includegraphics[width=4.2cm]{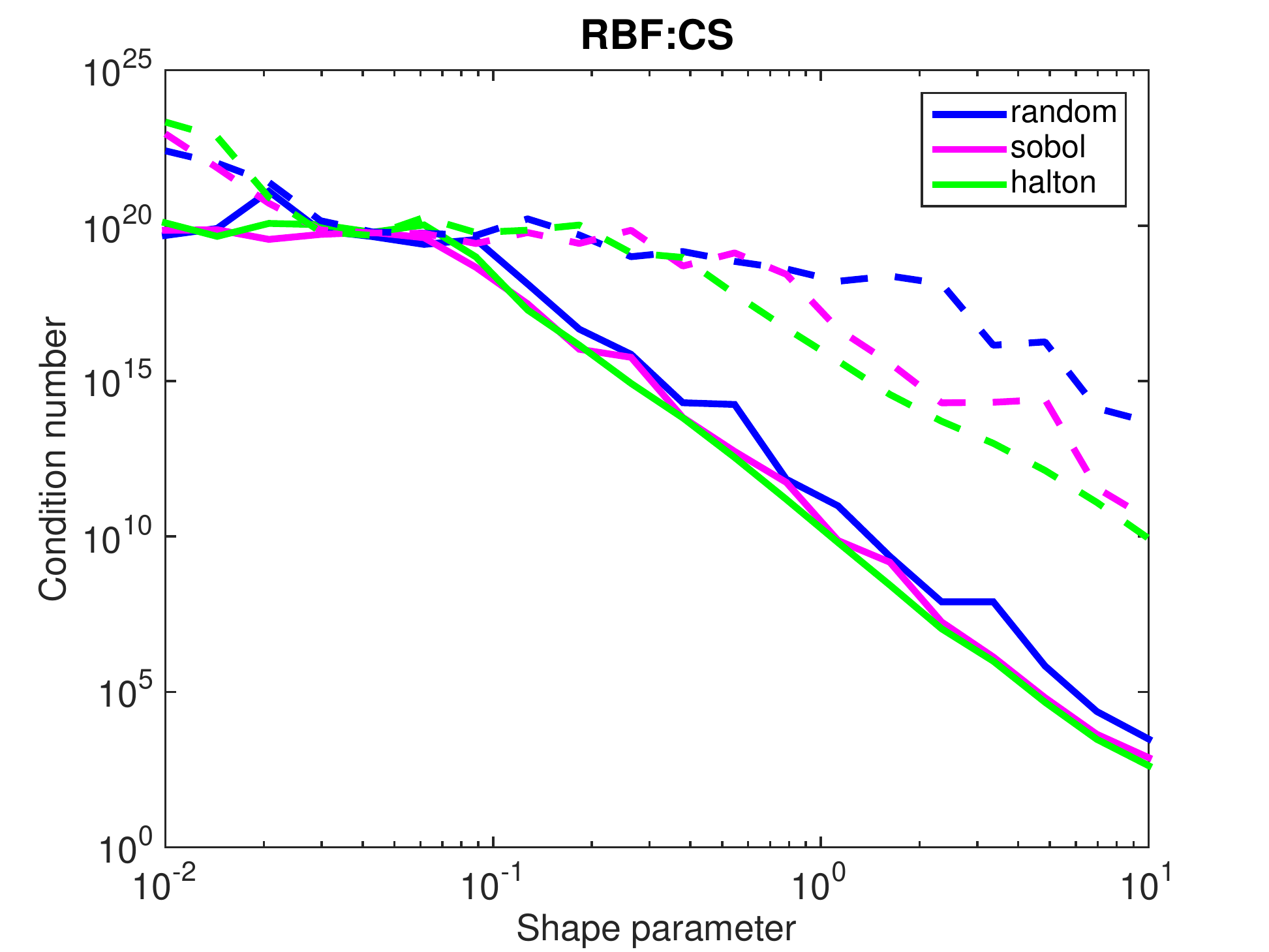}
\end{center}
  \caption{Curves of the condition numbers for the Lagrange interpolation (solid curves) and the Hermite interpolation (dashed curves) w.r.t. the shape paramter using Guassian, IMQ and CS, respectively. Left: $N=50$; Middle: $N=100$; Right: $N=300$; $d=2$.
  \label{cond:ga_imq_ep}
    }
\end{figure}

\begin{figure}[htbp]
\begin{center}
    \includegraphics[width=4.2cm]{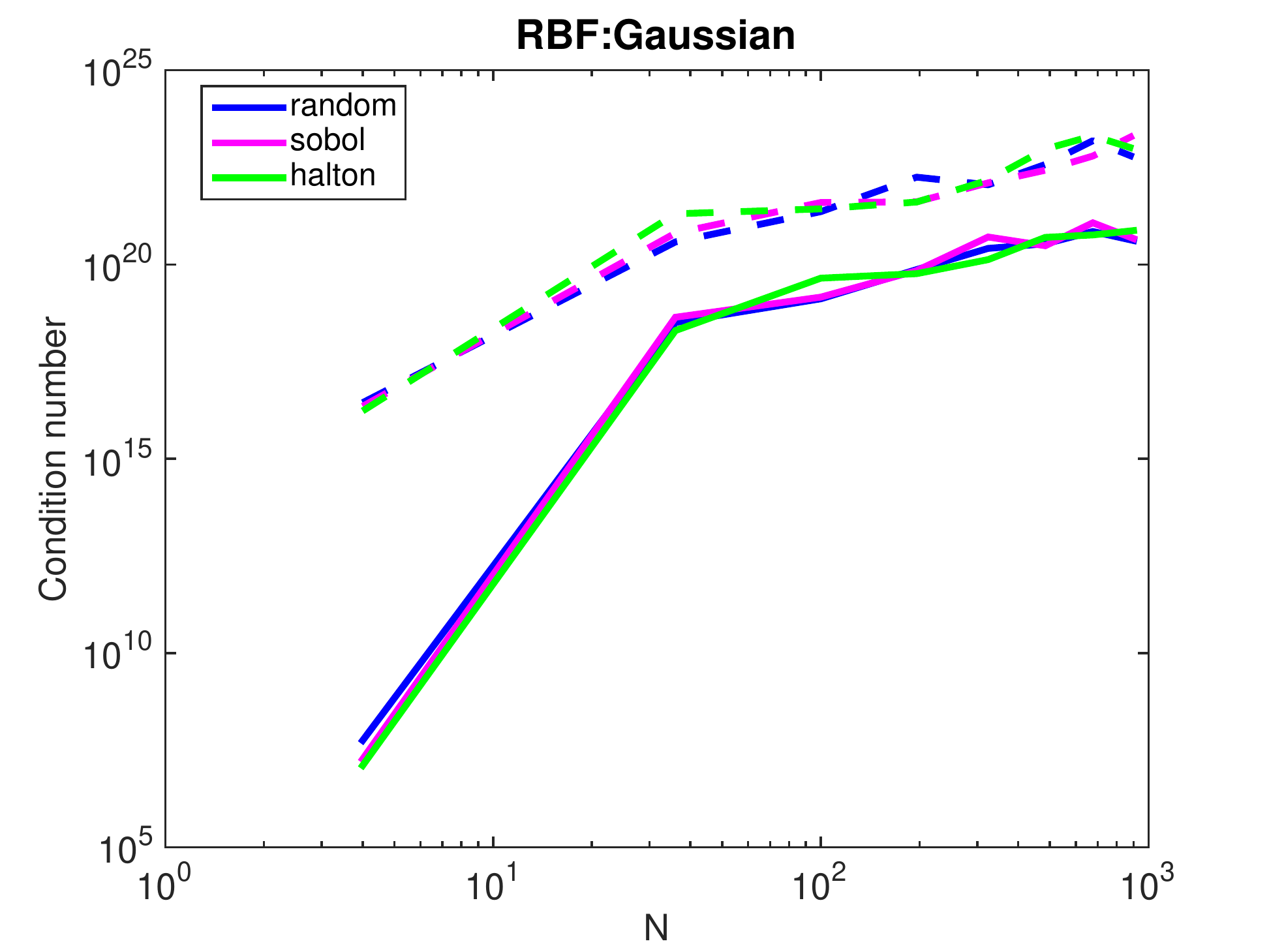}
    \includegraphics[width=4.2cm]{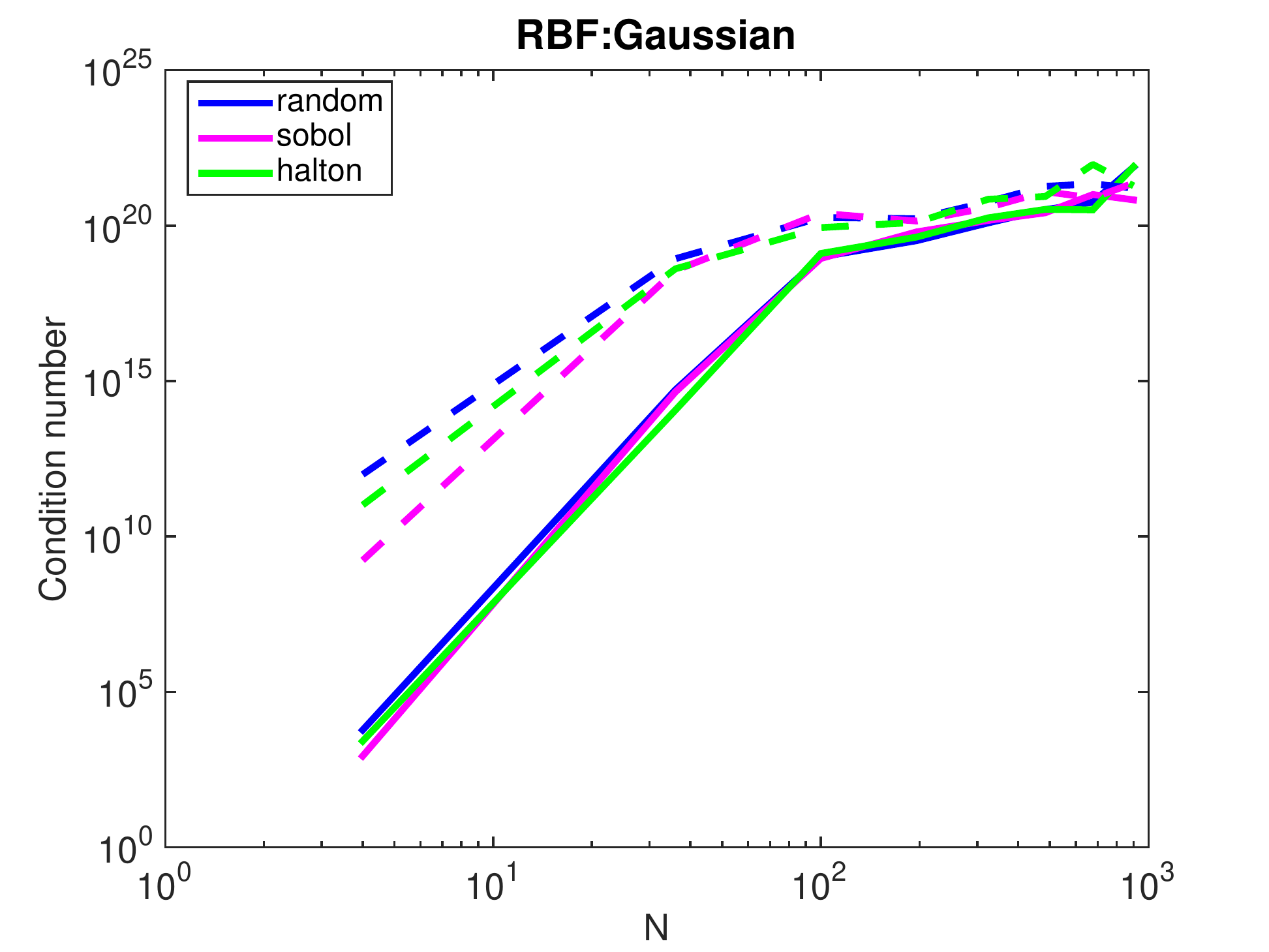}
    \includegraphics[width=4.2cm]{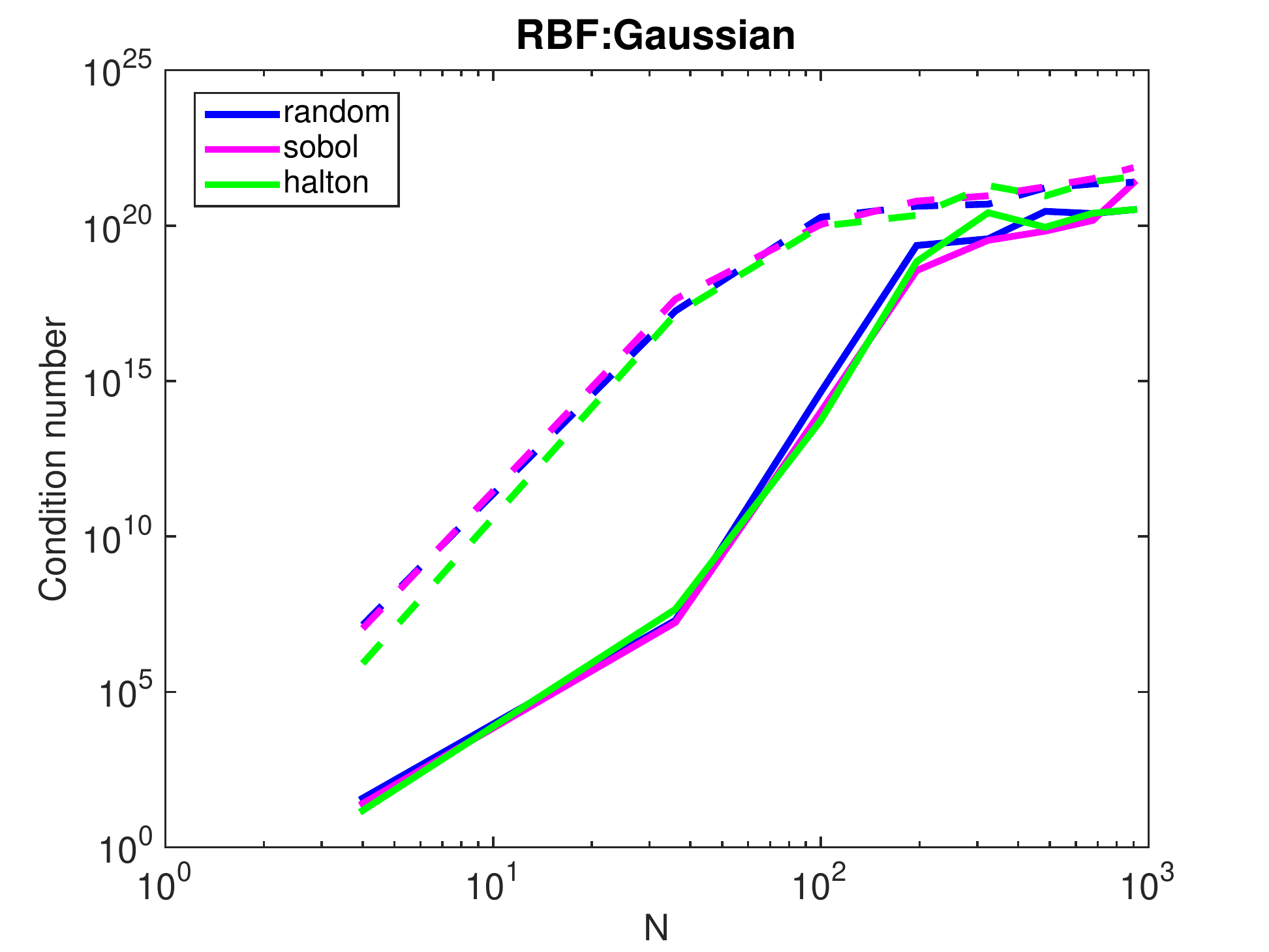}
     \includegraphics[width=4.2cm]{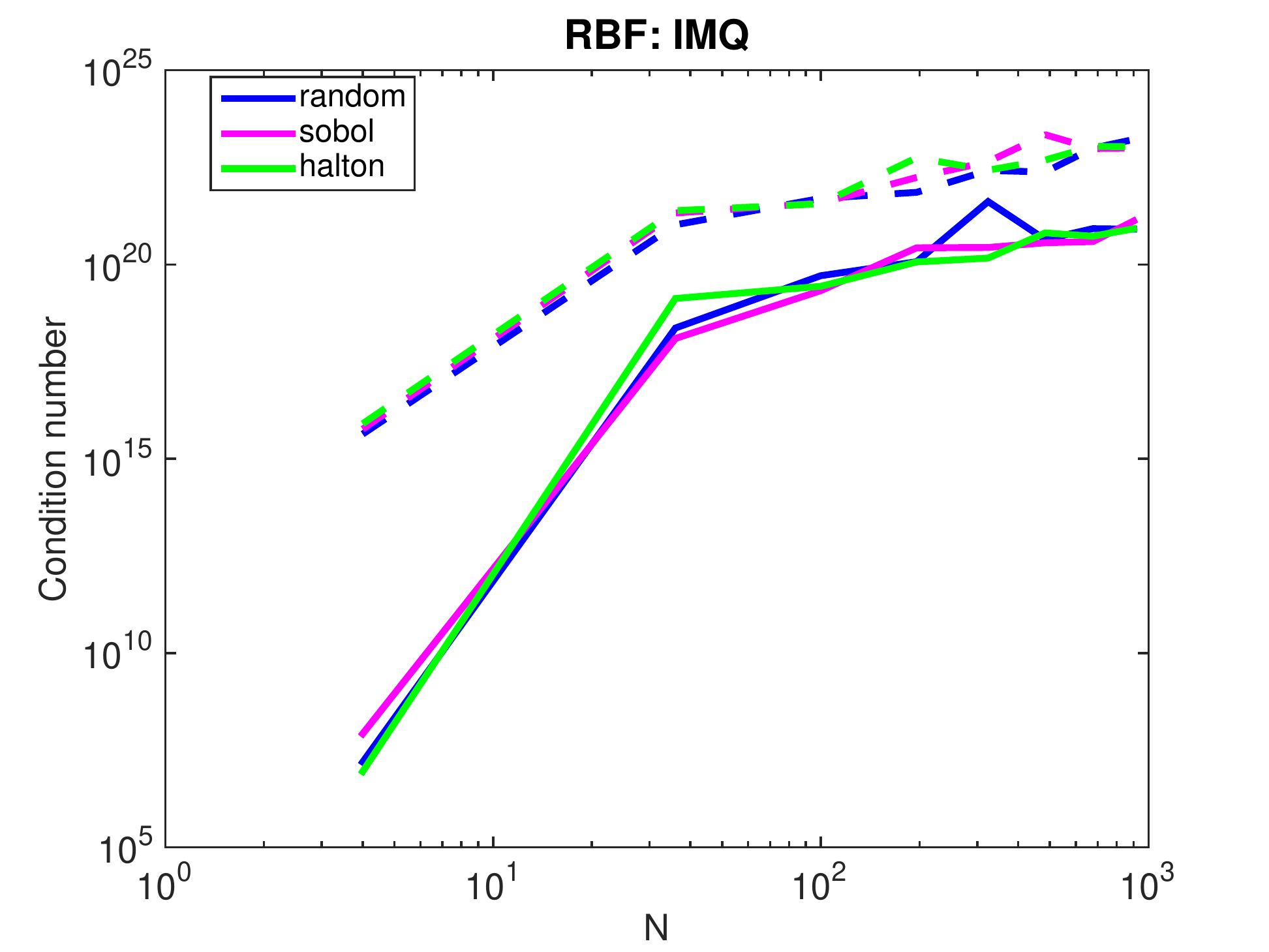}
   \includegraphics[width=4.2cm]{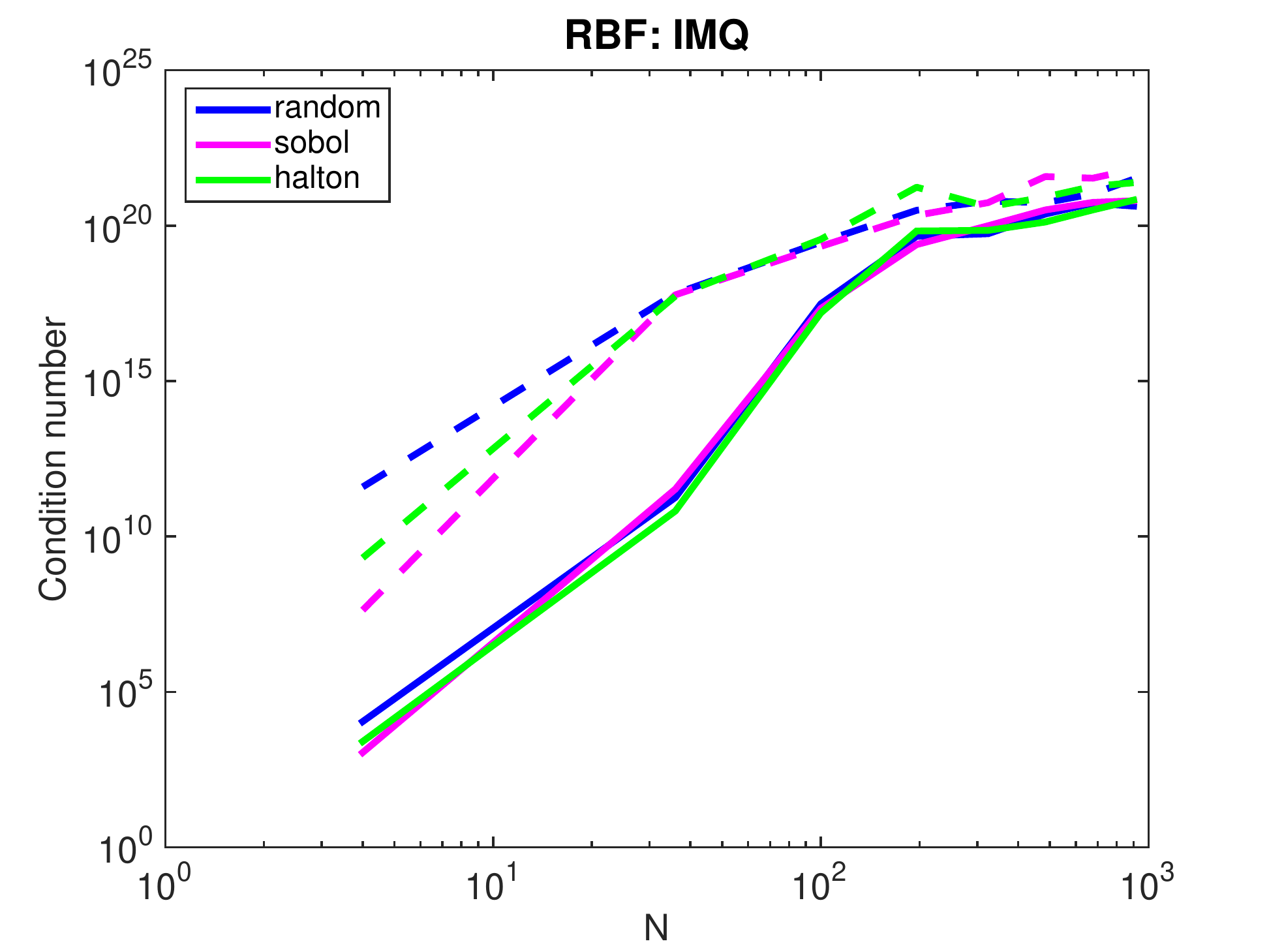}
   \includegraphics[width=4.2cm]{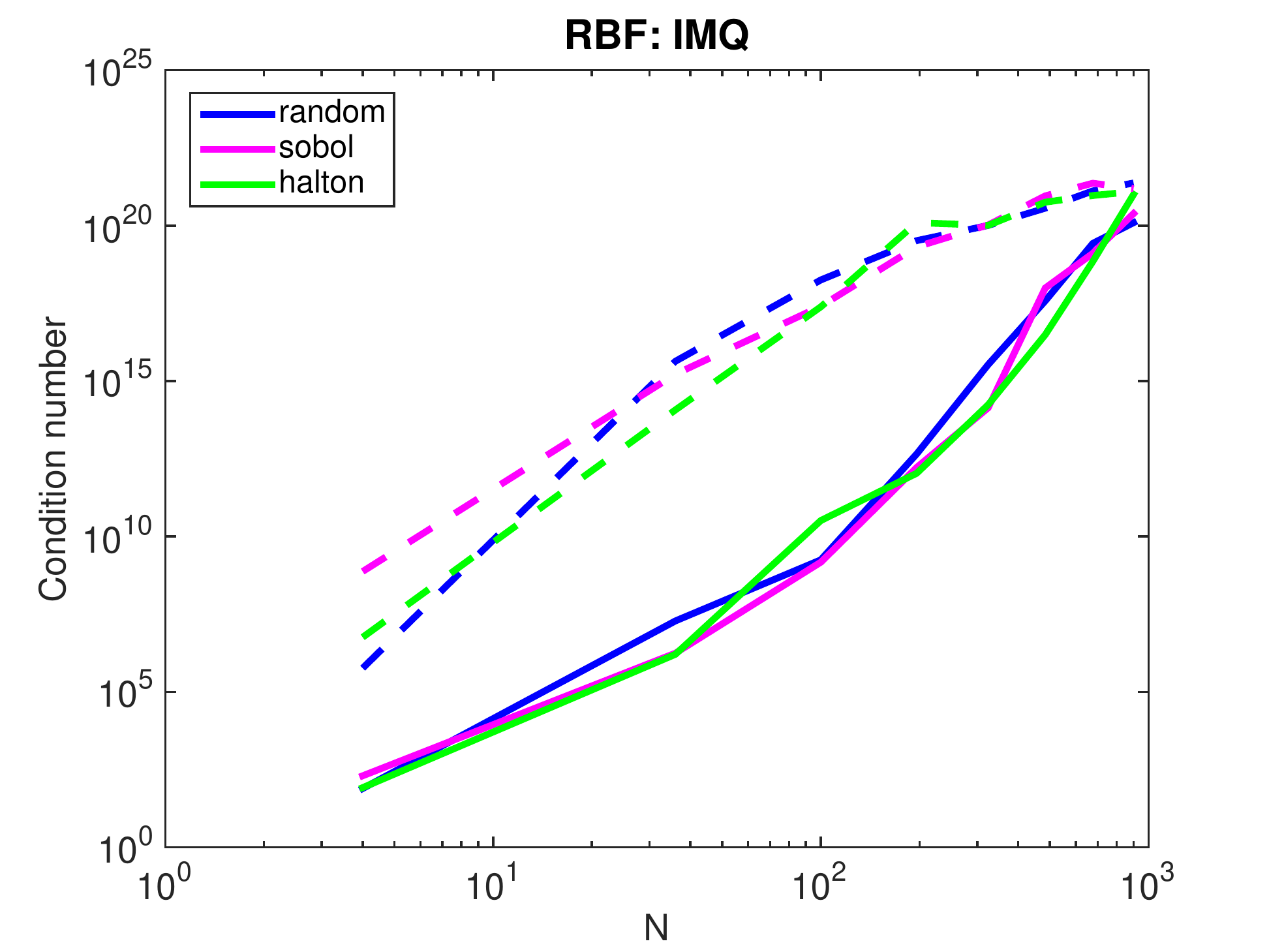}
       \includegraphics[width=4.2cm]{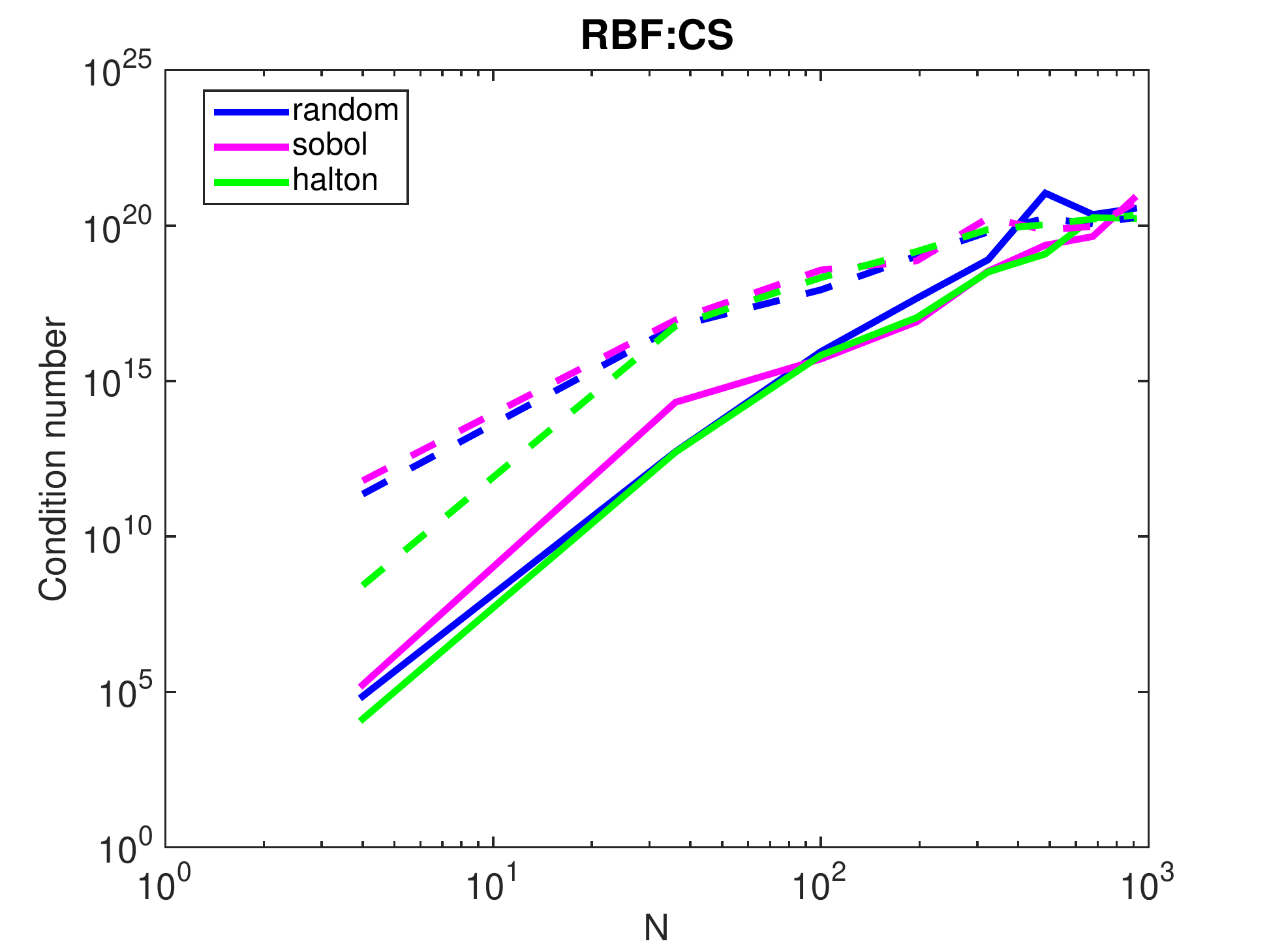}
   \includegraphics[width=4.2cm]{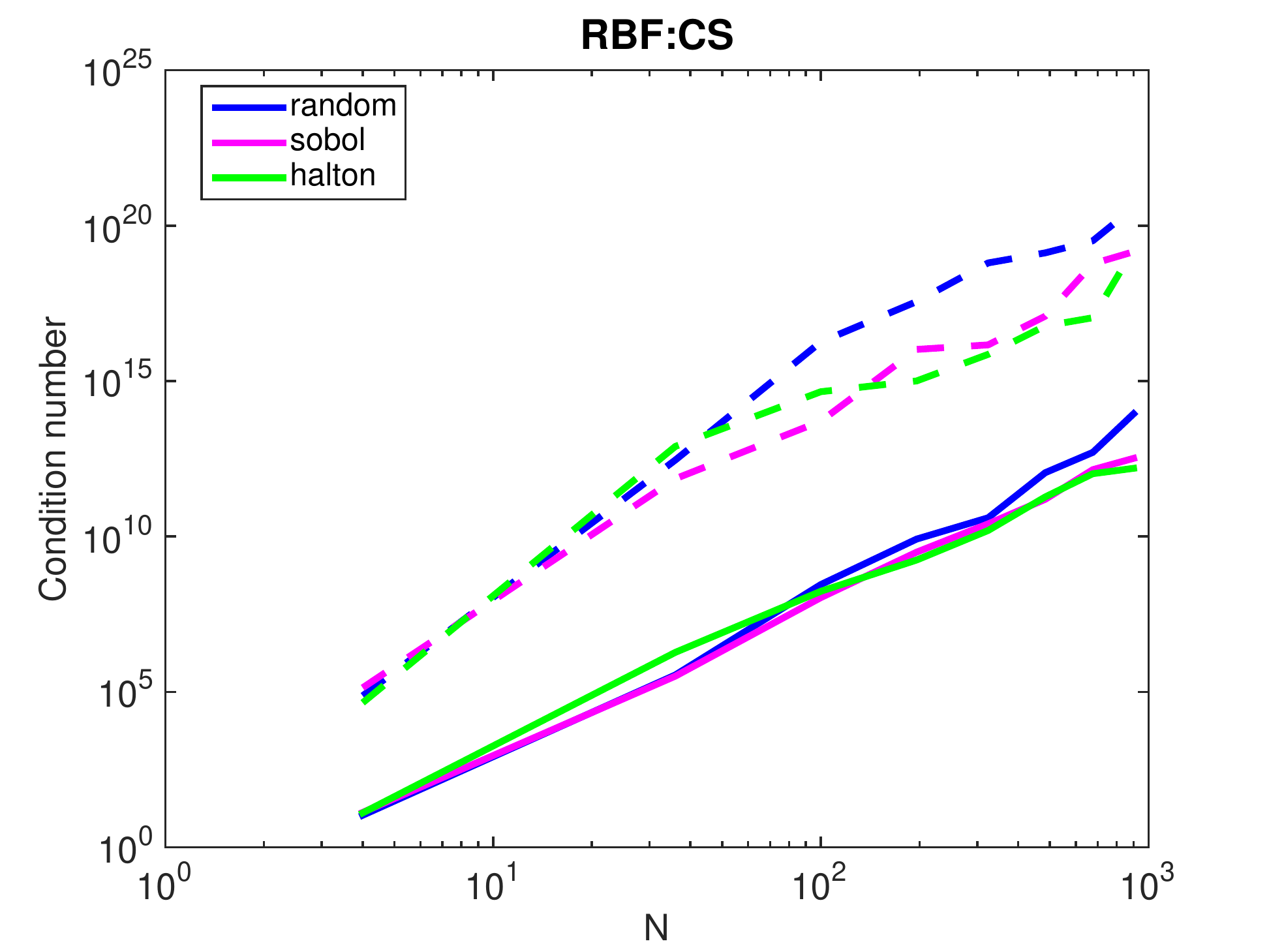}
   \includegraphics[width=4.2cm]{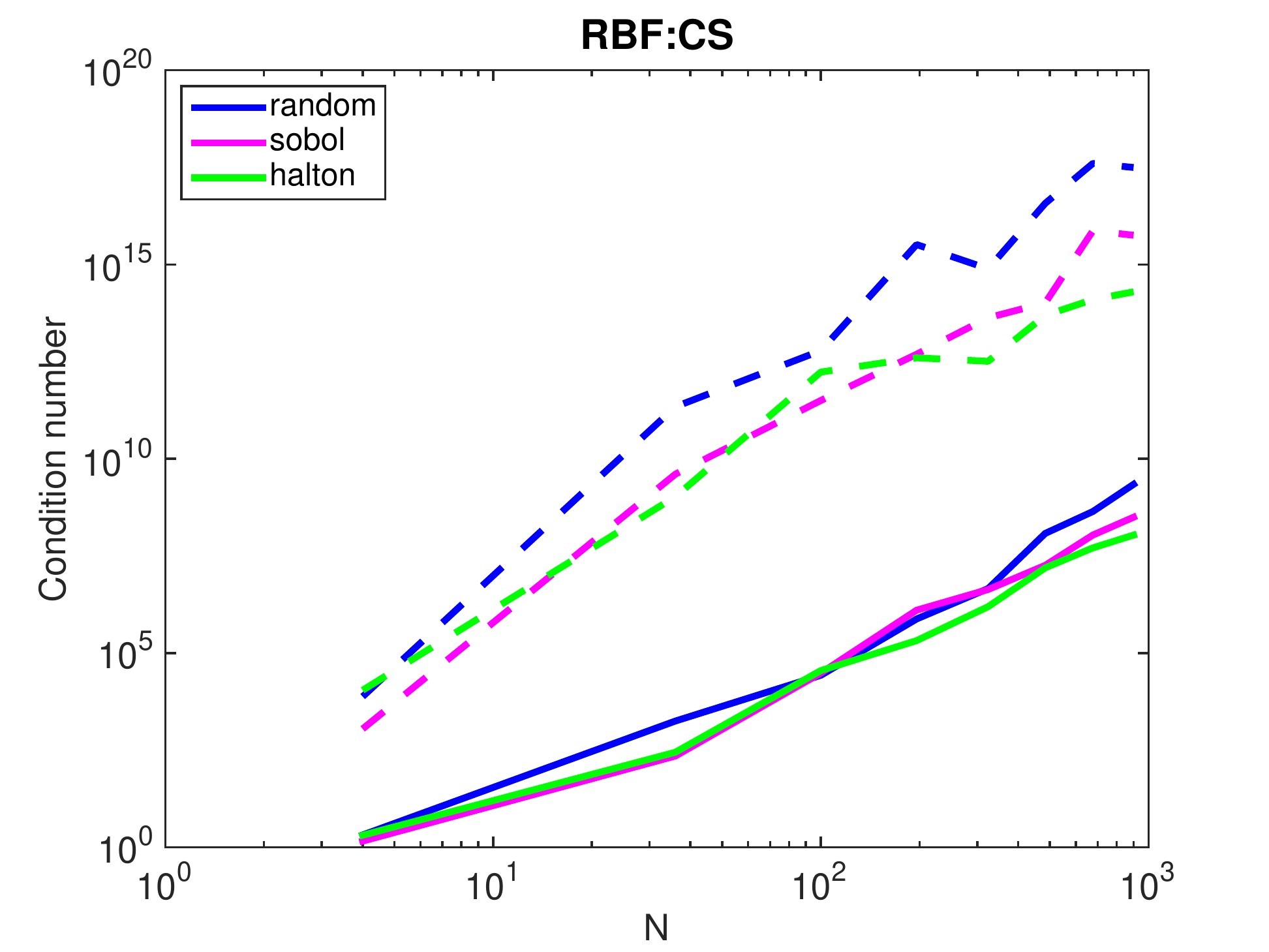}
\end{center}
  \caption{Curves of the condition numbers for the Lagrange interpolation (solid curves) and the Hermite interpolation (dashed curves) w.r.t. the number of sampling points $N$ using Guassian, IMQ and CS, respectively. Left: $\epsilon=0.1$; Middle: $\epsilon=1$; Right: $\epsilon=3$; $d=2$.
  \label{cond:ga_imq_N}
    }
\end{figure}

Note that the matrix $\mb{B}$ coincides with the joint covariance matrix of the gradient-enhanced GP emulator using the covariance kernel $K$, see \ref{app: ge-GP}.  Similar to the gradient-enhanced GP emulator, the condition numbers of the  matrix $\mb{B}$ are much larger than those of the ordinary Lagrange interpolation matrix $\mb{A}$ \cite{He+Chien2018}.  Figs. \ref{cond:ga_imq_ep} and \ref{cond:ga_imq_N} present the condition numbers of the system matrix for the design matrixes $\mb{A}$ and $\mb{B}$. Clearly, the design matrix becomes more singular as the shape parameters are closer to zero or the total number of points becomes large.

\section{Methodology} \label{sec:method}

It is obvious that design of the collocation points plays a important role in kernel interpolation. However, selecting optimal locations for trial centers in $\Xi$ is highly non-trivial. These are the motivation of the development of adaptive solvers \cite{Hon+Schaback2003, Ling+Schaback2009, Schaback+Wendland2000}. Nevertheless, it is also believed that uniformly distributed interpolation points can lead to better approximation results  \cite{Scheuerer2011}. Popular choices of such points include the Sobol' points,  the low discrepancy points  or even the random distributed points.  However, it is observed that all these choices  have the instability issue when the number of collocation points becomes large, say $N\sim \mathcal{O} (10^2)$, see Fig. \ref{cond:ga_imq_N}.

In this work, we shall propose a type of quasi-optimal collocation set for kernel interpolations, by searching from a large set of \textit{candidate} points. The idea is inspired by the computation of ``Fekete points" in polynomial approximation. 

\subsection{Fekete-type interpolation points}
We recall the kernel interpolation (\ref{interplant})
\begin{align*}
u_N(Z)=\sum^N_{j=1} c_{j}\, \Phi (\epsilon\|Z-z_j\|) = \sum_{j=1}^N d_{j}\, \ell_j \,(Z),
\end{align*}
where the latter equality writes the solution in terms of cardinal Lagrange interpolants, given by
\begin{align}\label{eq:lagrange}
  \ell_j \, (Z)&= \frac{\det \mathbf{A}(z_1,...,z_{j-1},Z,z_{j+1},...,z_N)}
  { \det \mathbf{A}(z_1,...,z_{j-1},z_j,z_{j+1},...,z_N))}, & j&=1,...,N.
\end{align}
The Lebesgue constant, corresponding to the operator norm of interpolation, is defined as
\begin{align*}
\Lambda_N := \max_{Z\in I_Z} \sum_{j=1}^N \left|\ell_j(Z)\right|.
\end{align*}
Small values of $\Lambda_N$ indicate that interpolation is a stable operation. The so called Fekete points are a configuration of $(z_1, \ldots, z_N)$ that maximizes the Vandermonde-like matrix determinant, namely,
\begin{align}\label{eq:fekete}
  \Xi^* = \max_{\Xi=\{z_1,...,z_N\}\subset I_Z} \Big\{\left| \det \mathbf{A}(z_1,...,z_N) \right| \Big\}
\end{align}
With the relation \eqref{eq:lagrange}, one observes that the configuration $\Xi^\ast$ ensures that $|\ell_j(Z)| \leq 1$ for all $Z \in I_Z$, and thus ensures $\Lambda_N \leq N$. In practice, the behavior of Lebesgue constant for Fekete points is frequently sublinear. Thus, Fekete points provide a strategy for stable interpolation. Unfortunately, there is no known way to explicitly characterize or compute these points outside of special polynomial approximation cases. In addition, computationally solving the optimization problem \eqref{eq:fekete} is a daunting task. A typical approach is to relax the optimization problem by seeking a greedy (iterative) solution, i.e.,
\begin{align*}
  z_{N+1} = \argmax_{z \in I_Z} \left| \det \mathbf{A}(z_1, \ldots, z_N) \right|.
\end{align*}
In the polynomial approximation literature such an appropriate is titled the ``approximate Fekete point" approach, and for kernel interpolation it is a power function maximization approach. We will utilize this procedure in our gradient-enhanced setting to choose centers for approximation.

\subsection{Generating quasi-optimal points via greedy selection} \label{sec:methodology}

In this section, we present a data-independent algorithm to select optimal distribution of centers.  The data-independent feature is useful for practical problems, as it allows one to determine, prior to conducting expensive simulations or experimentations, where to collect data samples.

The general goal is to achieve a set of determinant-maximizing (i.e., Fekete) points associated to the gradient-enhanced design matrix $\bs{B}$. However, direct optimization is likely to be computationally challenging since the optimization problem is non-convex, and the dimension of the optimization variables can be very large when a large interpolation grid is used. Therefore, we opt for a greedy strategy, performing the sequential optimization,
\begin{align*}
  z_{N+1} &\coloneqq \argmax_{z \in I_Z} \det \bs{B}(z_1, \ldots, z_N, z),
\end{align*}
where $(z_1, \ldots, z_N) \subset \R^d$ is a given set of $N$ centers. The procedure above defines an iterative scheme. Although in principle computing this determinant requires construction of the matrix $\bs{B}(z_1, \ldots, z_N, z)$, such an expensive construction can be avoided by exercising low-rank updates of $L U$ factorizations.
First we note that given $z_1, \ldots, z_N$, we have
\begin{align*}
  \bs{B}(z_1, \ldots, z_N,z) &\coloneqq
  \left(\begin{array}{cccc}
    \bs{A}_{0,0} & \bs{A}_{0,1} & \cdots & \bs{A}_{0,d} \\
    \bs{A}_{1,0} & \bs{A}_{1,1} & \cdots & \bs{A}_{1,d} \\
    \vdots                   & \vdots                   & \ddots & \vdots \\
    \bs{A}_{d,0} & \bs{A}_{d,1} & \cdots & \bs{A}_{d,d}
  \end{array}\right),
\end{align*}
where above $\bs{A}_{m,n} = \bs{A}_{m,n}(z_1, \ldots, z_N, z)$. Then if $\bs{P}_{N+1} \in \R^{(N+1)(d+1) \times (N+1)(d+1)}$ is a permutation matrix that maps element $(N+1) i$ to element $N(d+1) + i$ for $i = 1, \ldots, d+1$ (and retains the order for the remaining $N (d+1)$ elements in the first $N (d+1)$ entries), then
\begin{align*}
  \bs{B}(z_1, \ldots, z_N,z) = \bs{P}_{N+1} \bs{B} \bs{P}_{N+1}^{-T} &= \left(\begin{array}{cc} \bs{B}_N & \bs{W}(z) \\
  \bs{W}^T(z) & \bs{B}(z) \end{array}\right)
\end{align*}
Then the determinant can be written in terms of the Schur complement of this latter expression,
\begin{align*}
  \det \bs{B}(z_1, \ldots, z_N,z) &= \det \left( \bs{P} \bs{B} \bs{P}^{T}\right) = \det \left(\begin{array}{cc} \bs{B}_N & \bs{W} \\
                                                                                                        \bs{W}^T & \bs{B}(z) \end{array}\right) \\
                                  &= \det \bs{B}_N \det \left( \bs{B}(z) - \bs{W}^T \bs{B}_N^{-1} \bs{W} \right),
\end{align*}
where we have introduced the notation $\bs{B}_N \coloneqq \bs{B}(z_1, \ldots, z_N)$. Thus,
\begin{align*}
  \argmax_{z \in I_Z} \det \bs{B}(z_1, \ldots, z_N, z) &= \argmax_{z \in I_Z} F(z), \\
  F(z) &\coloneqq \det \left( \bs{B}(z) - \bs{W}^T \bs{B}_N^{-1} \bs{W} \right).
\end{align*}
For a fixed $z$, the major cost of evaluting $F$ involves computing the solution $\bs{X} \in \R^{N(d+1) \times (d+1)}$ to the linear system $\bs{B}_N \bs{X} = \bs{W}$. There is a secondary, comparatively negligible, cost of computing the determinant of a $(d+1) \times (d+1)$ matrix. For the linear system, we can write the solution in terms of the Cholesky factor of the permuted version of $\bs{B}_N$,
\begin{align*}
  \bs{P}_N \bs{B}_N \bs{P}_N^{T} = \bs{L}_N \bs{L}_N^T
\end{align*}
so that, if $\bs{L}_N^{-1}$ is available and stored, then evaluation of $F$ requires application of $\bs{L}_N^{-1}$, with a $\mathcal{O}(d^3 N^2)$ cost, and we can rewrite the optimization as,
\begin{align}\label{eq:z-optimization}
  z_{N+1} = \argmax_{z \in I_Z} F(z) &= \argmax_{z \in I_Z} \det \left( \bs{B}(z) - \bs{V}^T \bs{V} \right), & \bs{V} &\coloneqq \bs{L}_N^{-1} \bs{W}(z).
\end{align}
To aid in continuing the iteration, an efficient update that transforms $\bs{L}^{-1}_N$ into $\bs{L}^{-1}_{N+1}$ (the Cholesky factor for the pivoted $\bs{B}_{N+1}$) is helpful. Again exercising Schur complements, we find that
\begin{align*}
  \bs{P}_{N+1} \bs{B}_{N+1} \bs{P}_{N+1}^{T} = \bs{L}_{N+1} \bs{L}_{N+1}^T \hskip 10pt \Longrightarrow \hskip 10pt \bs{L}_{N+1} = \left(\begin{array}{cc} \bs{L}_N & \bs{0}_{N(d+1) \times (d+1)} \\
                                                                                                                            \bs{V}^T & \widetilde{\bs{L}} \end{array}\right),
\end{align*}
where $\widetilde{\bs{L}} = \widetilde{\bs{L}}(z)$ is the Cholesky factor of the Schur complement of the block $\bs{B}(z_{N+1})$ of the matrix $\bs{B}_{N+1}$,
\begin{align}\label{eq:schur-complement}
  \widetilde{\bs{L}} \widetilde{\bs{L}}^T = \bs{B}(z_{N+1}) - \bs{V}^T \bs{V}.
\end{align}
Then a computation shows that
\begin{align}\label{eq:L-update}
  \bs{L}^{-1}_{N+1} = \left(\begin{array}{cc} \bs{L}^{-1}_N & \bs{0}_{N(d+1) \times (d+1)} \\
  -\widetilde{\bs{L}}^{-1} \bs{V}^T \bs{L}_N^{-1} & \widetilde{\bs{L}}^{-1} \end{array}\right).
\end{align}
In summary, given $\bs{L}_N^{-1}$ that identifies $z_{N+1}$ from the optimization \eqref{eq:z-optimization}, then we compute $\bs{L}_{N+1}$ by (i) computing $\bs{V}(z_{N+1})$ from \eqref{eq:z-optimization}, (ii) computing the inverse of the $(d+1) \times (d+1)$ Cholesky factor $\widetilde{\bs{L}}(z_{N+1})$ of the Schur complement from \eqref{eq:schur-complement}, (iii) assembling the matrix in \eqref{eq:L-update}, which requires an additional application of $\bs{L}_N^{-1}$.

Finally, we emphasize again that in practice we do not compute the solution to \eqref{eq:z-optimization}, but rather we solve the easier optimization problem that replaces the continuum $I_Z$ by a discrete candidate set.

\section{Numerical tests} \label{sec:tests}

In this section we use several numerical tests to demonstrate the benefit of the quasi-optimal interpolation points. In order to implement our proposed method,  we first pick $M$ random centers on $I_Z$ as a candidate set, then we select $N$ optimal points from Section \ref{sec:methodology} as  RBF centers. In the following examples, we choose $M=10^4$.

For comparisons against other methods, we employ the following sampling methods.
\begin{itemize}
\item Random sampling: $N$ samples are taken randomly in $I_Z$;
\item Sobol sampling:  choose $N$ values of the Sobol sequence;
\item Halton sampling: choose  $N$ values of the Halton sequence.
\end{itemize}

To measure the performance of an approximation, we will use the $l_2$ error (RMSE). Specifically given a set of $Q=1000$ random samples $\Xi_{test}=\{z_i\}^Q_{i=1}\in I_Z$ and samples of the true function $u(z_i)$ and the kernel approximation $s(z_i)$ we compute
\begin{equation}
E_{l_2}=\Big(\frac{1}{Q}\sum^Q_{i=1}|s(z_i)-u(z_i)|^2\Big)^{1/2}.
\end{equation}

In what follows we will use the term  ``Cholesky" for our proposed algorithm, ``random" to denote random sampling points, ``sobol" to denote Sobol points, and ``halton" to denote Halton points. Furthermore,  in all examples that follow we perform 50 trials of each procedure and report the mean results along with $20\%$ and $80\%$ quantiles.  For CS RBFs, we use the following  Wendlands compactly supported function $ \Phi_{d}(r)$,
which has the form
$$
    \Phi_{d}(r) = (1-r)_{+}^{l+3}\Big((l^3+9l^2+23l+15)r^3+(6l^2+36l+45)r^2+(15l+45)r+15\Big)/15
$$
with $ l=\left \lfloor \frac{d}{2} \right \rfloor+4.$

\subsection{Kernel interpolation for UQ}
In this section, some numerical examples are shown to demonstrate the stability and convergent behaviors of our proposed algorithms for UQ.
\subsubsection{Matrix stability}
\begin{figure}[htbp]
\begin{center}
    \includegraphics[width=6cm]{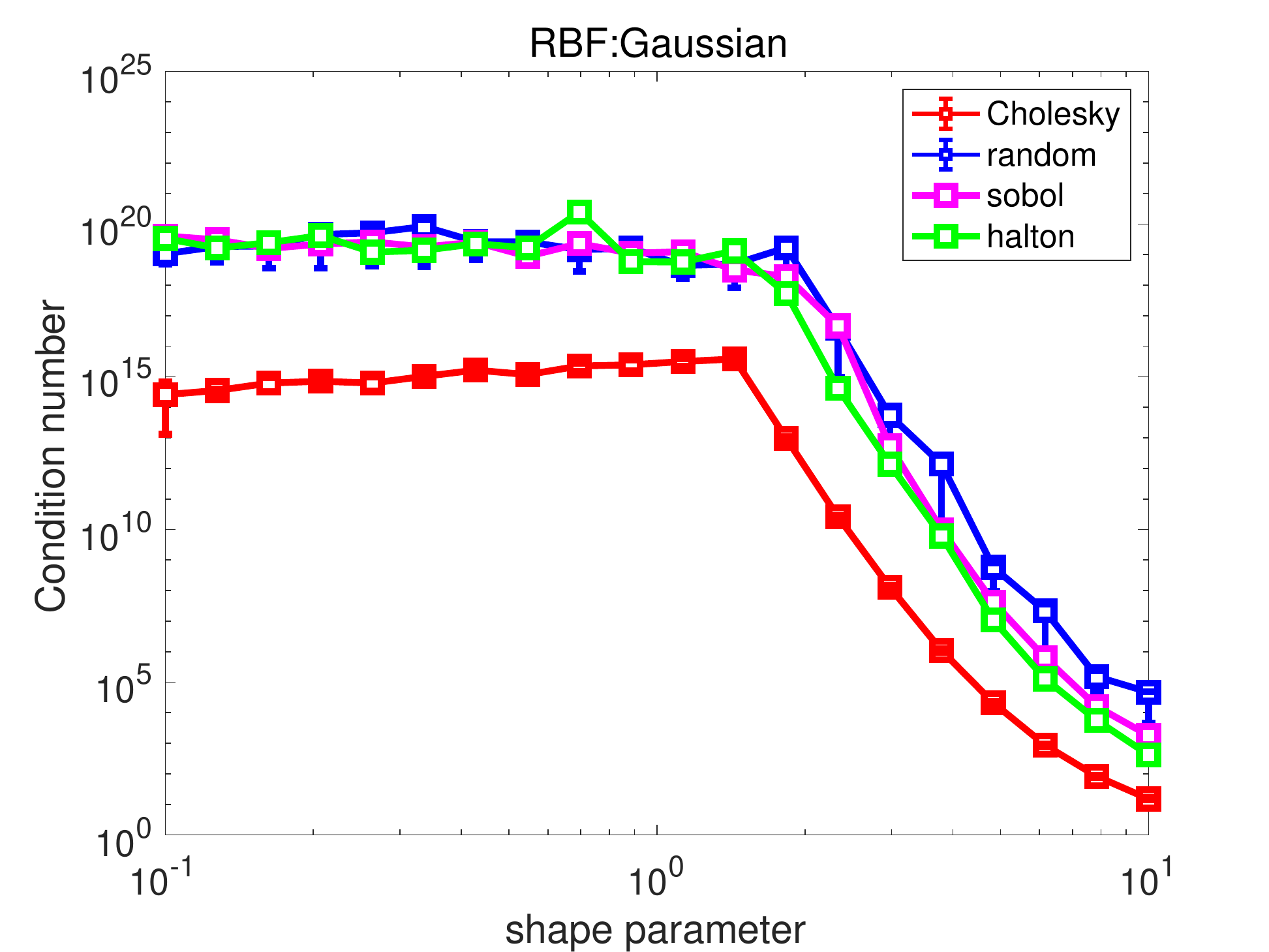}
    \includegraphics[width=6cm]{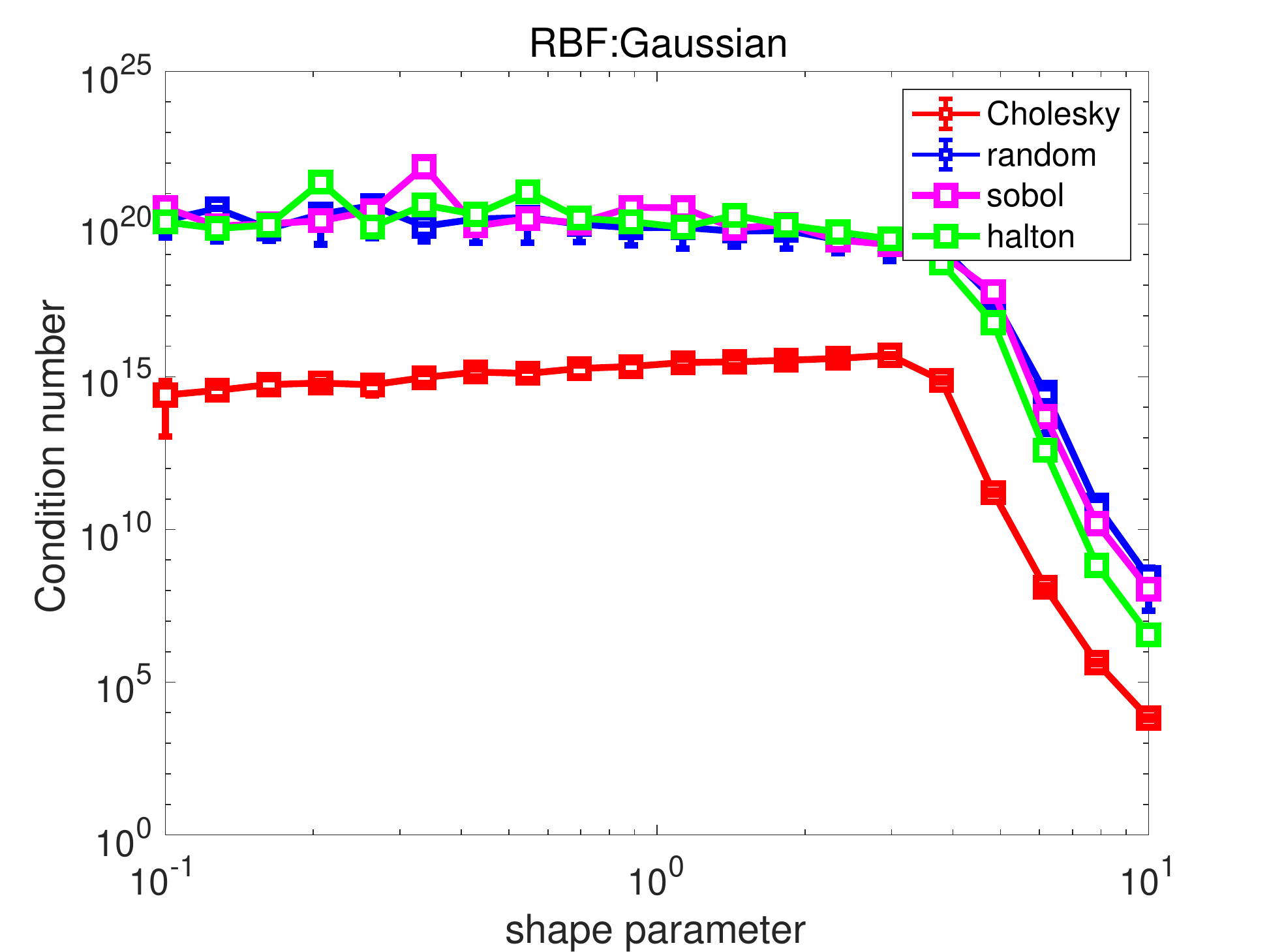}
   \includegraphics[width=6cm]{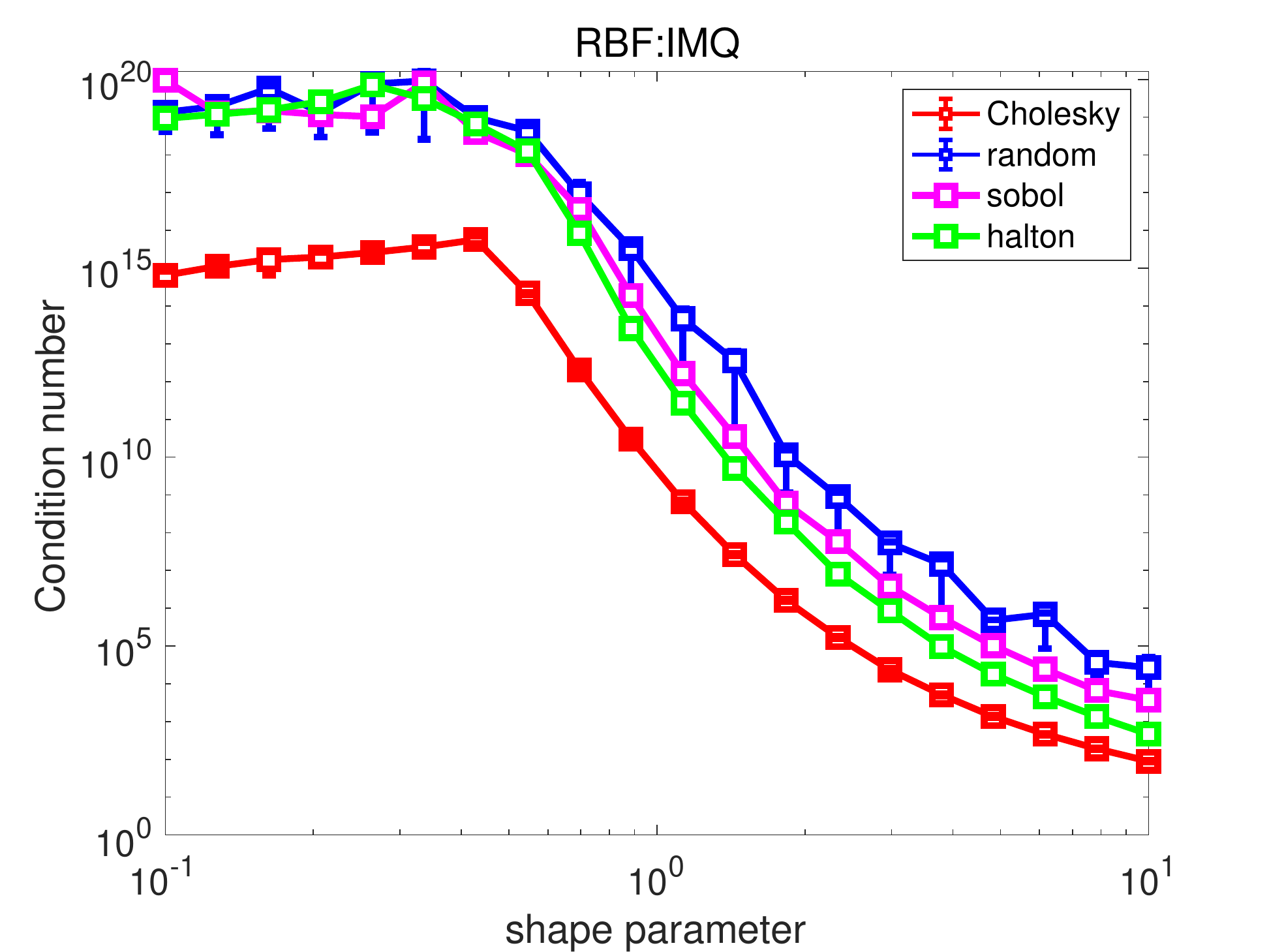}
   \includegraphics[width=6cm]{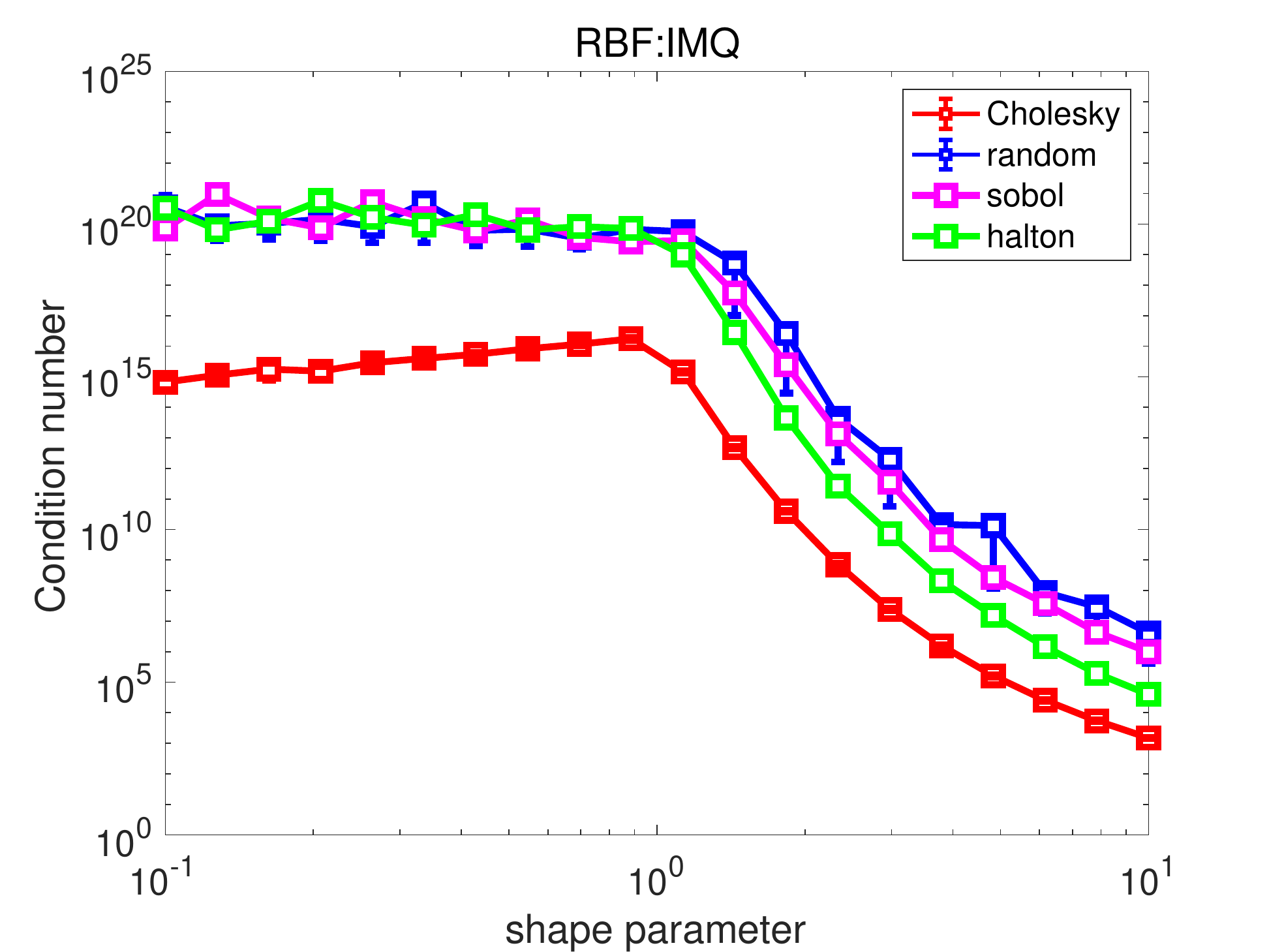}
      \includegraphics[width=6cm]{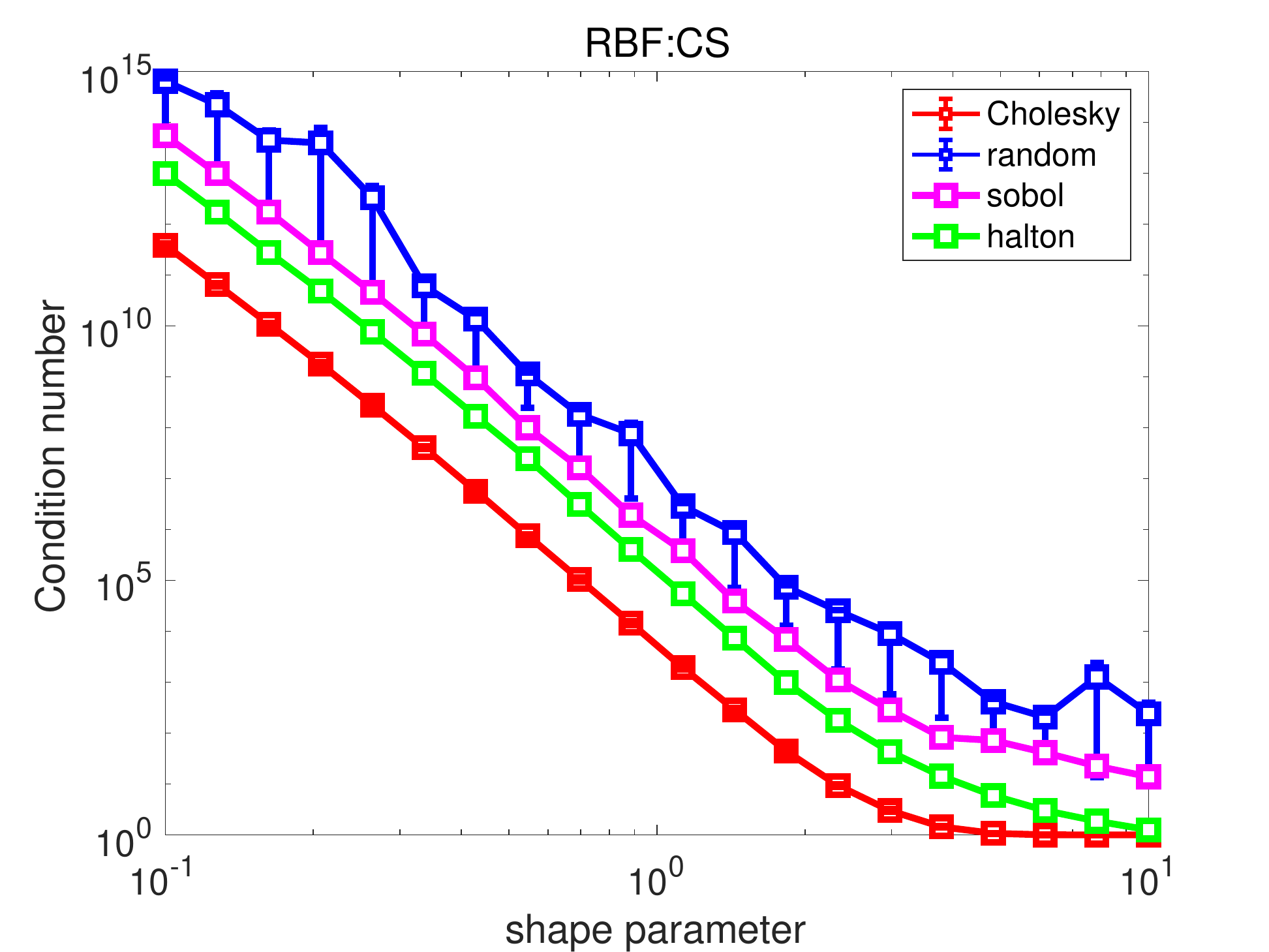}
   \includegraphics[width=6cm]{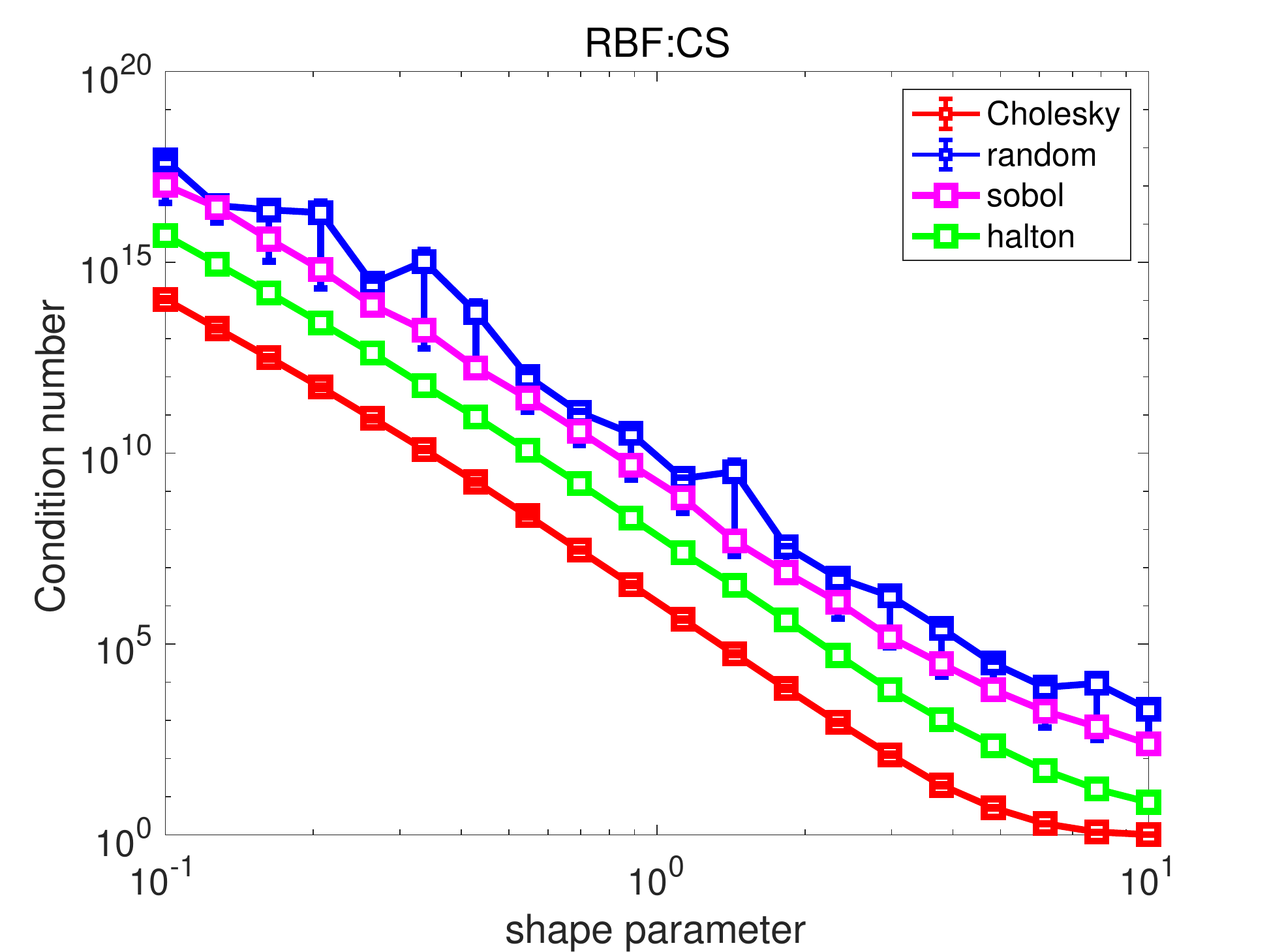}
\end{center}
  \caption{Condition numbers with respect to shape parameters using Gaussian, IMQ and CS. Left: $N=100$;  Right: $N=300$; $d=2.$
  \label{cond:ga_imq_ep_chol}
    }
\end{figure}

\begin{figure}[htbp]
\begin{center}
    \includegraphics[width=6cm]{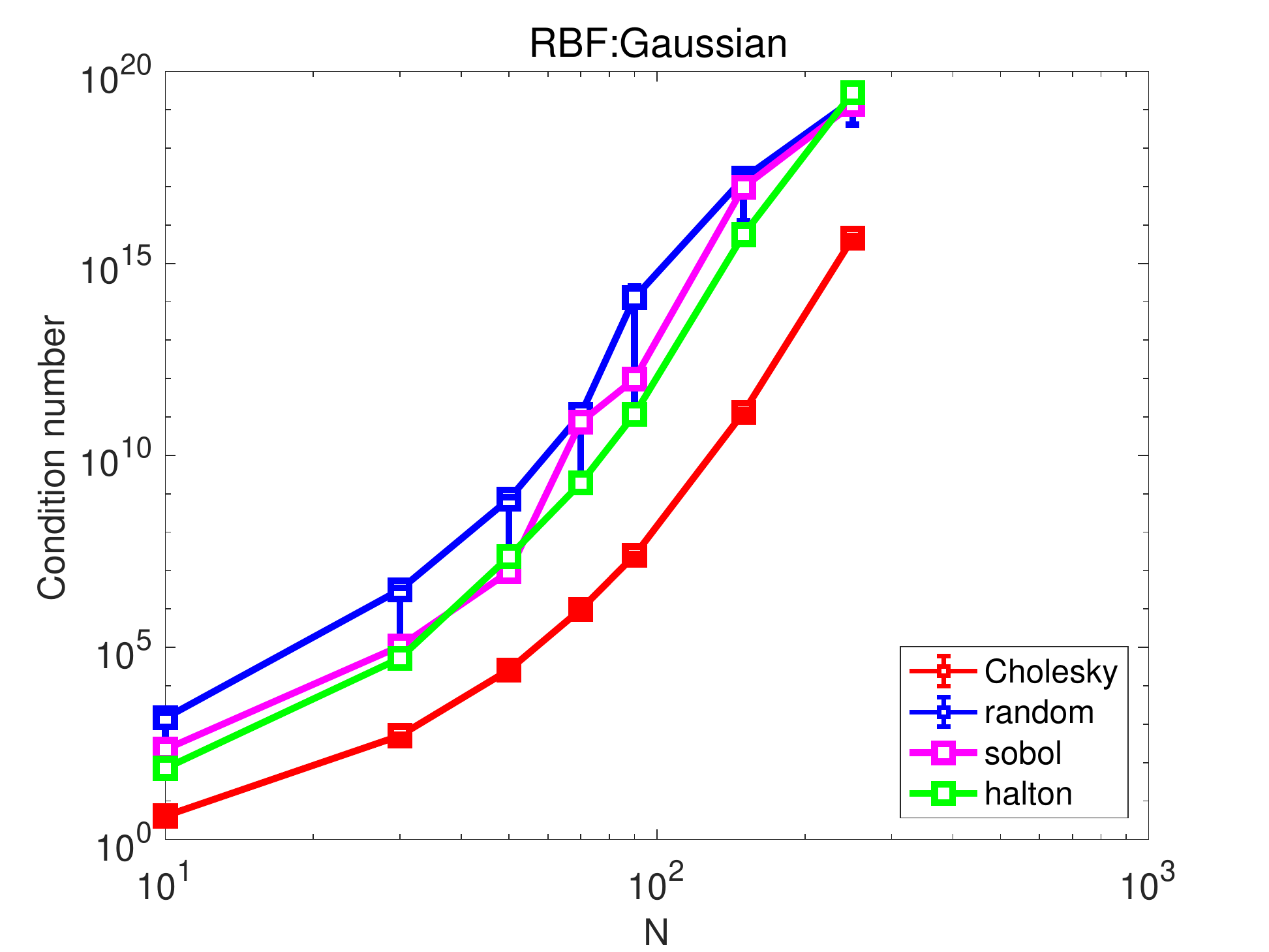}
    \includegraphics[width=6cm]{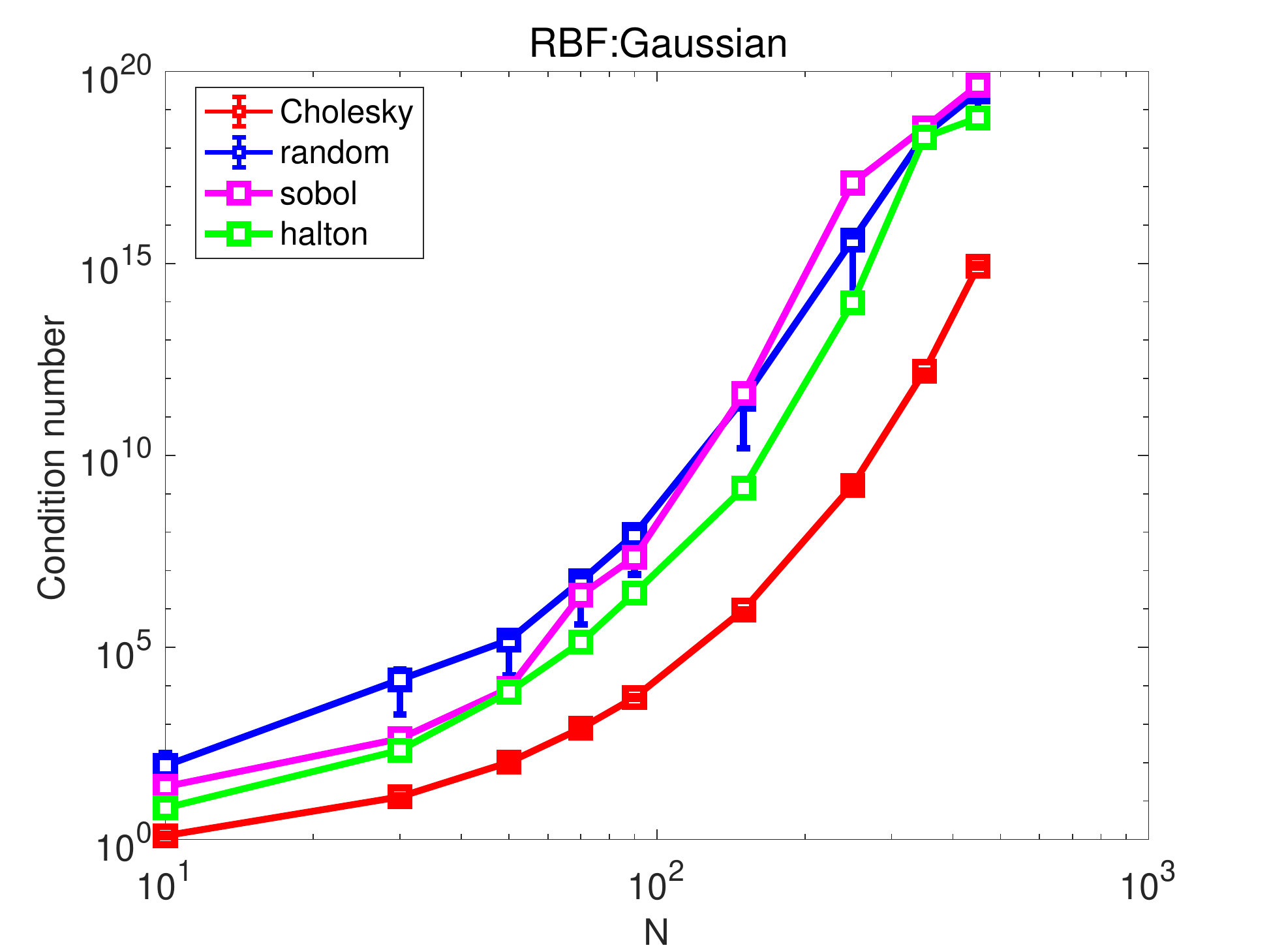}
    \includegraphics[width=6cm]{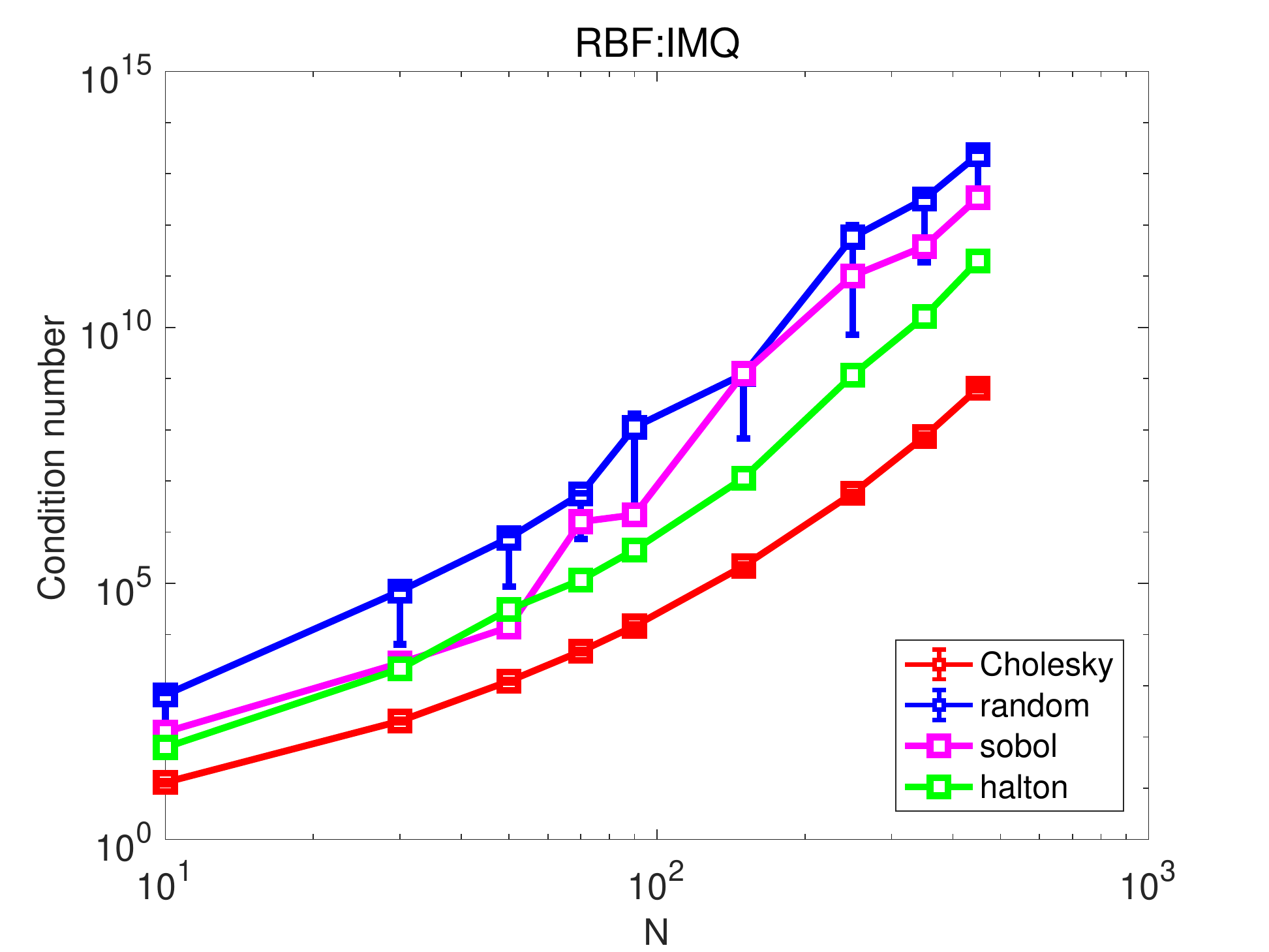}
     \includegraphics[width=6cm]{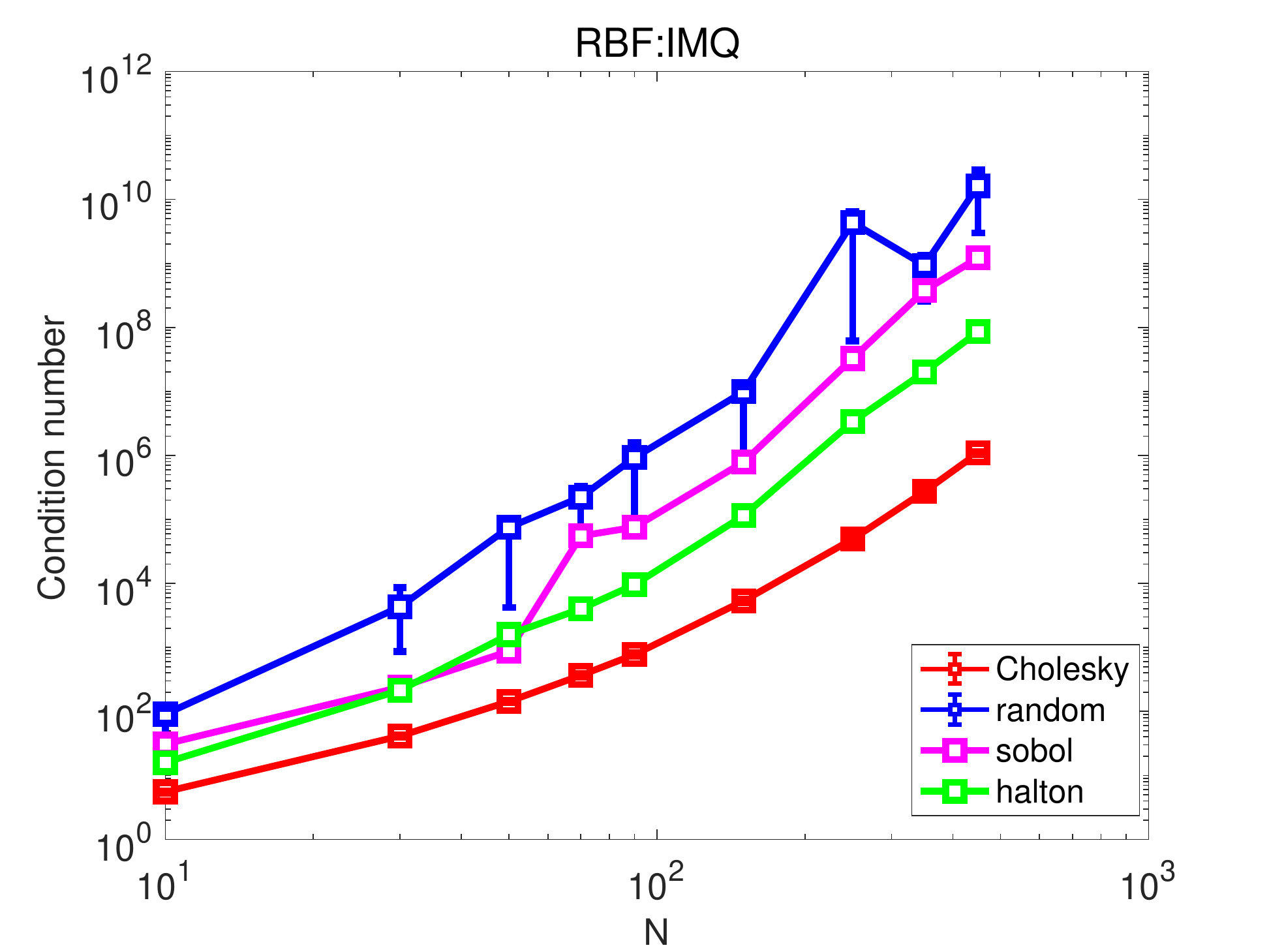}
        \includegraphics[width=6cm]{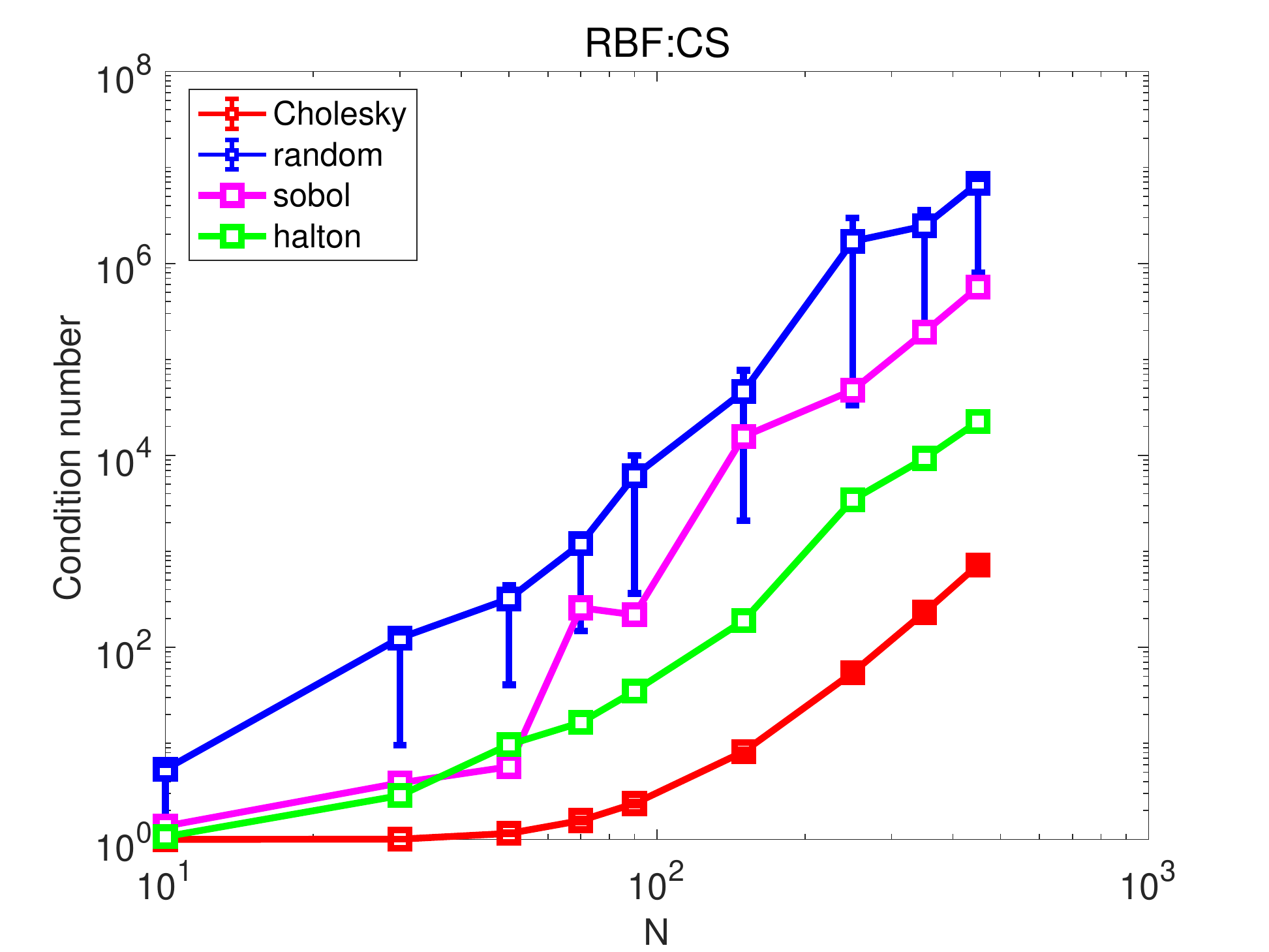}
     \includegraphics[width=6cm]{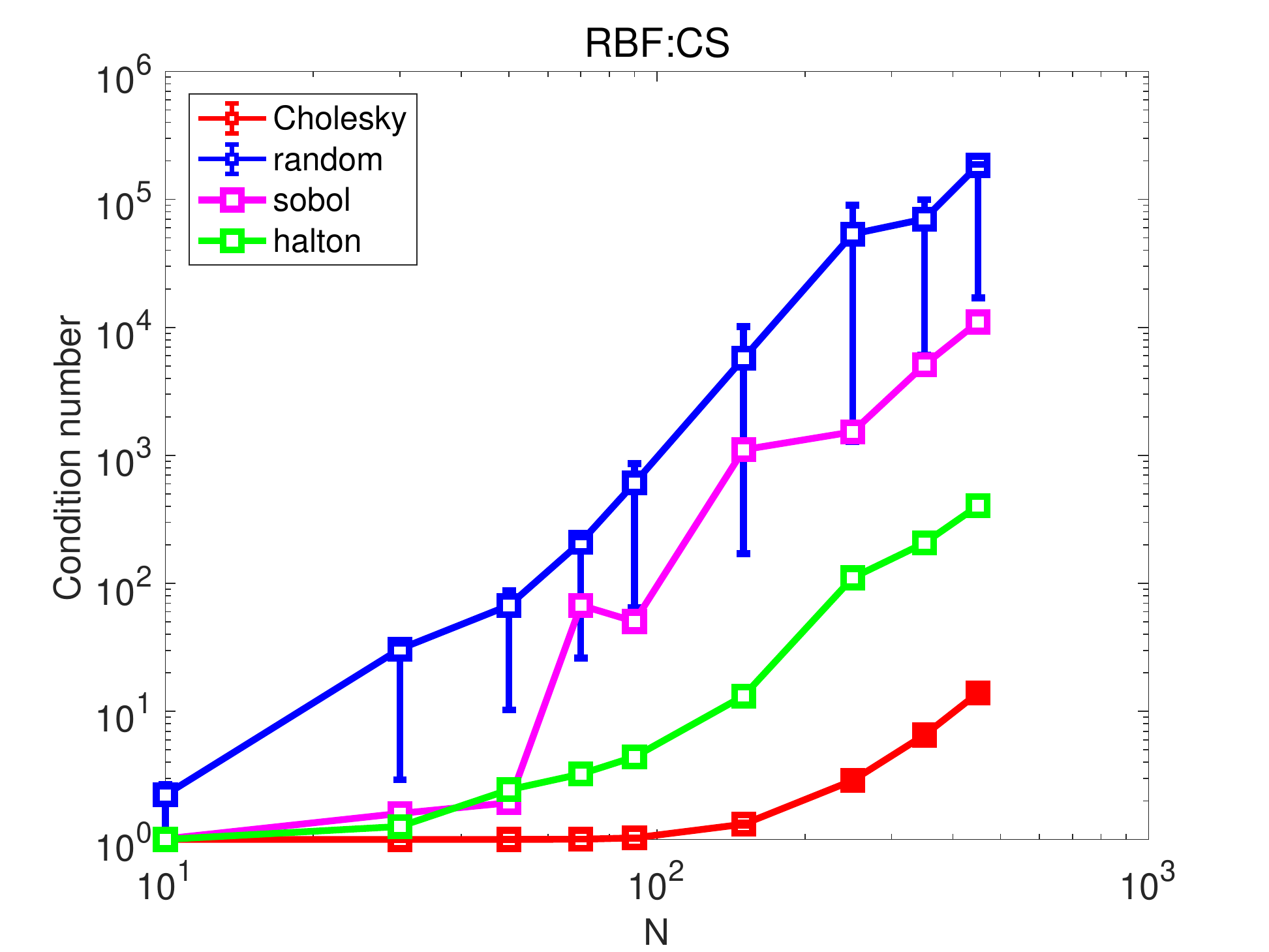}
\end{center}
  \caption{Condition numbers with respect to shape parameters with respect to the number of sample points $N$ using Gaussian, IMQ and CS. Left: $\epsilon=3$; Right: $\epsilon=5$; $d=2$.
  \label{cond:ga_imq_N_chol}
    }
\end{figure}

\begin{figure}[htbp]
\begin{center}
    \includegraphics[width=6cm]{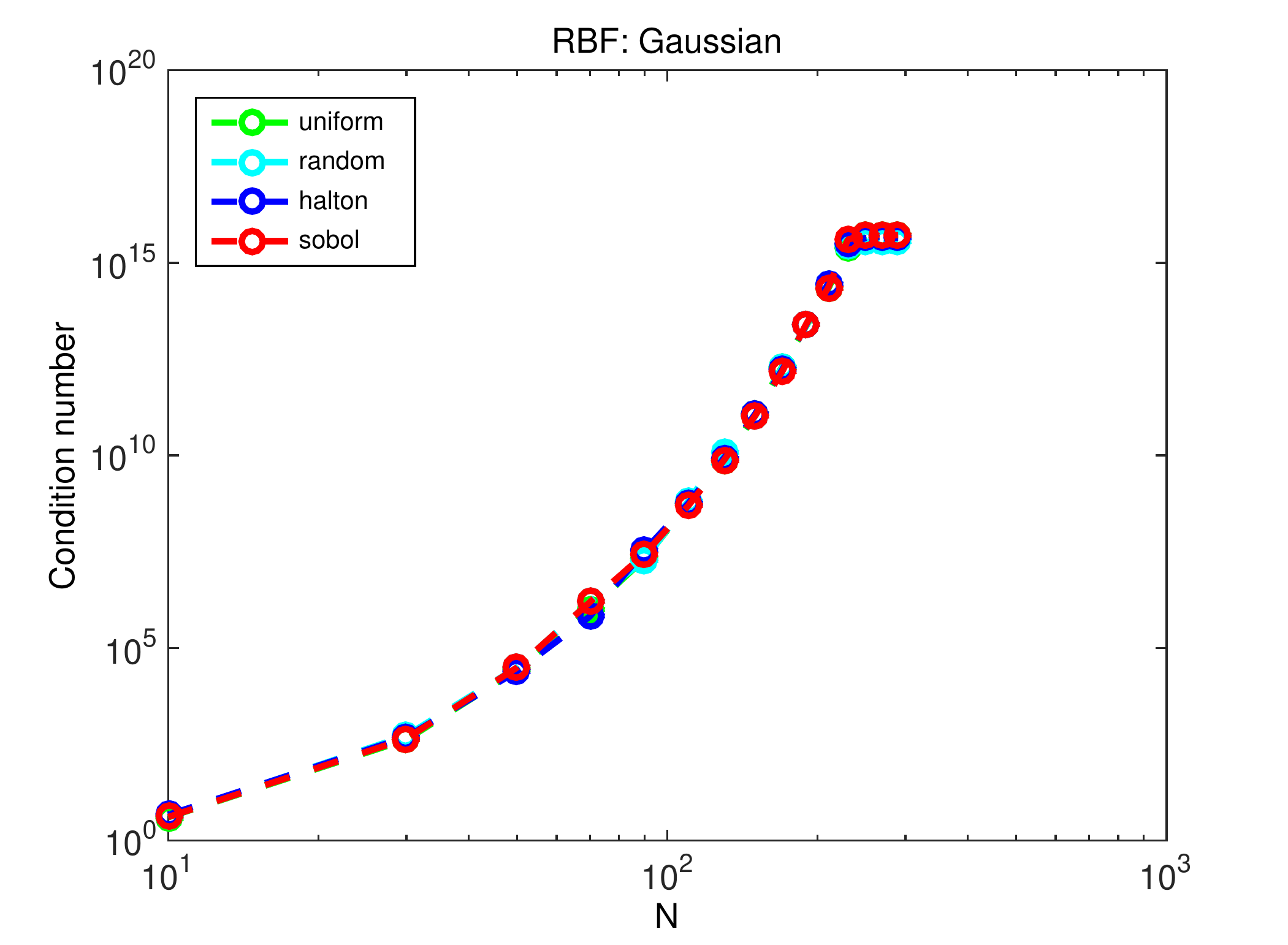}
     \includegraphics[width=6cm]{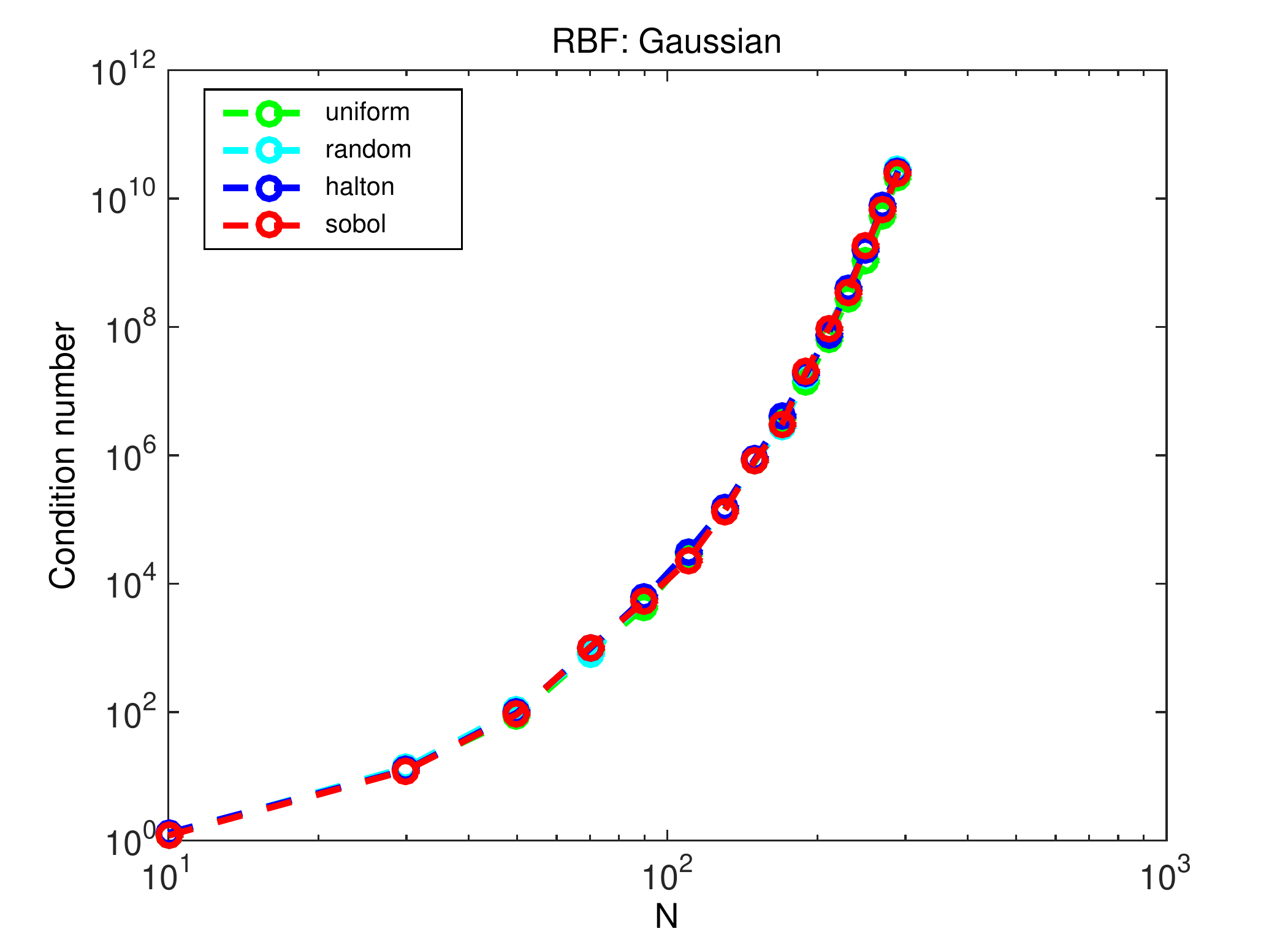}
       \includegraphics[width=6cm]{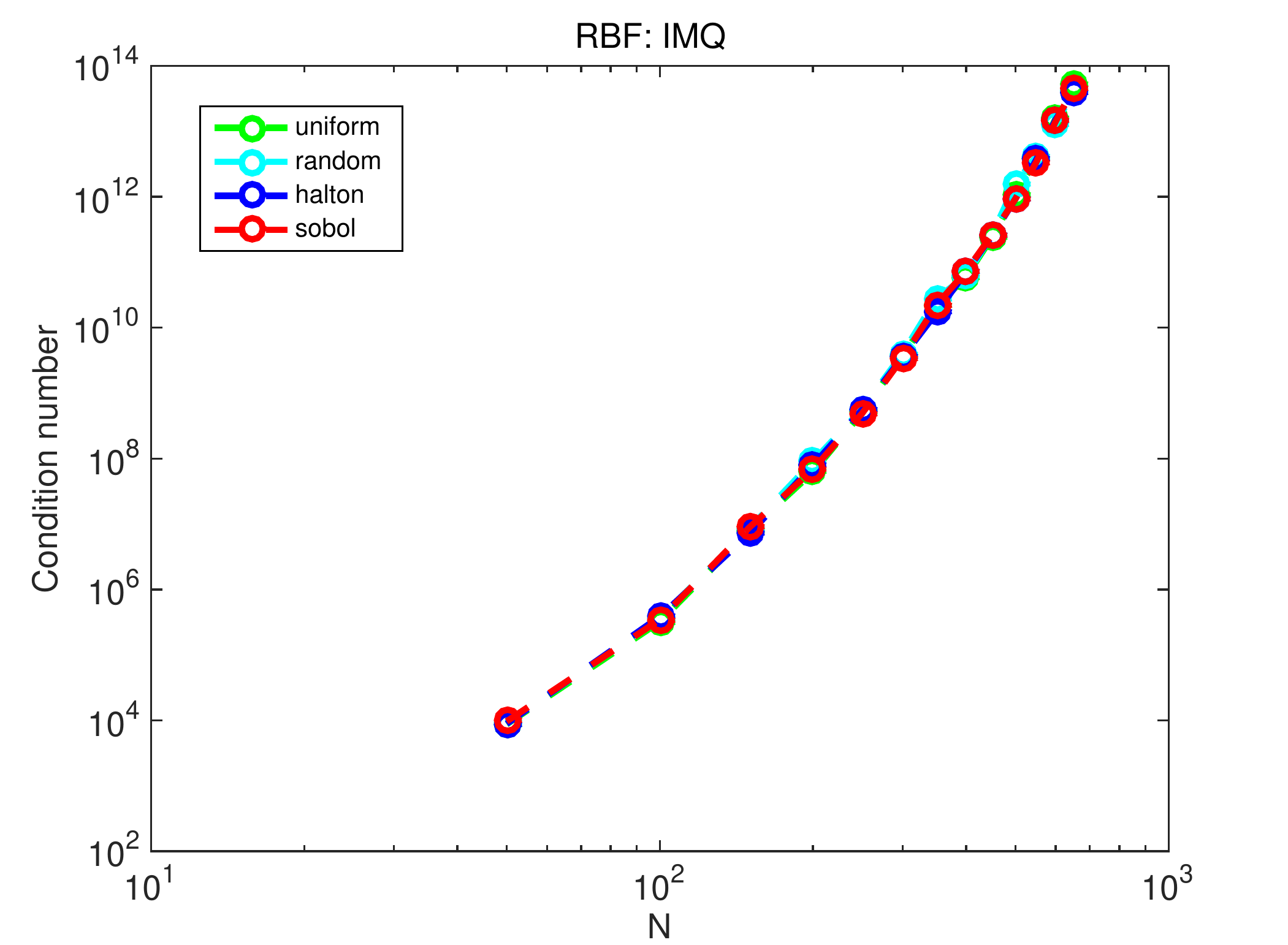}
 \includegraphics[width=6cm]{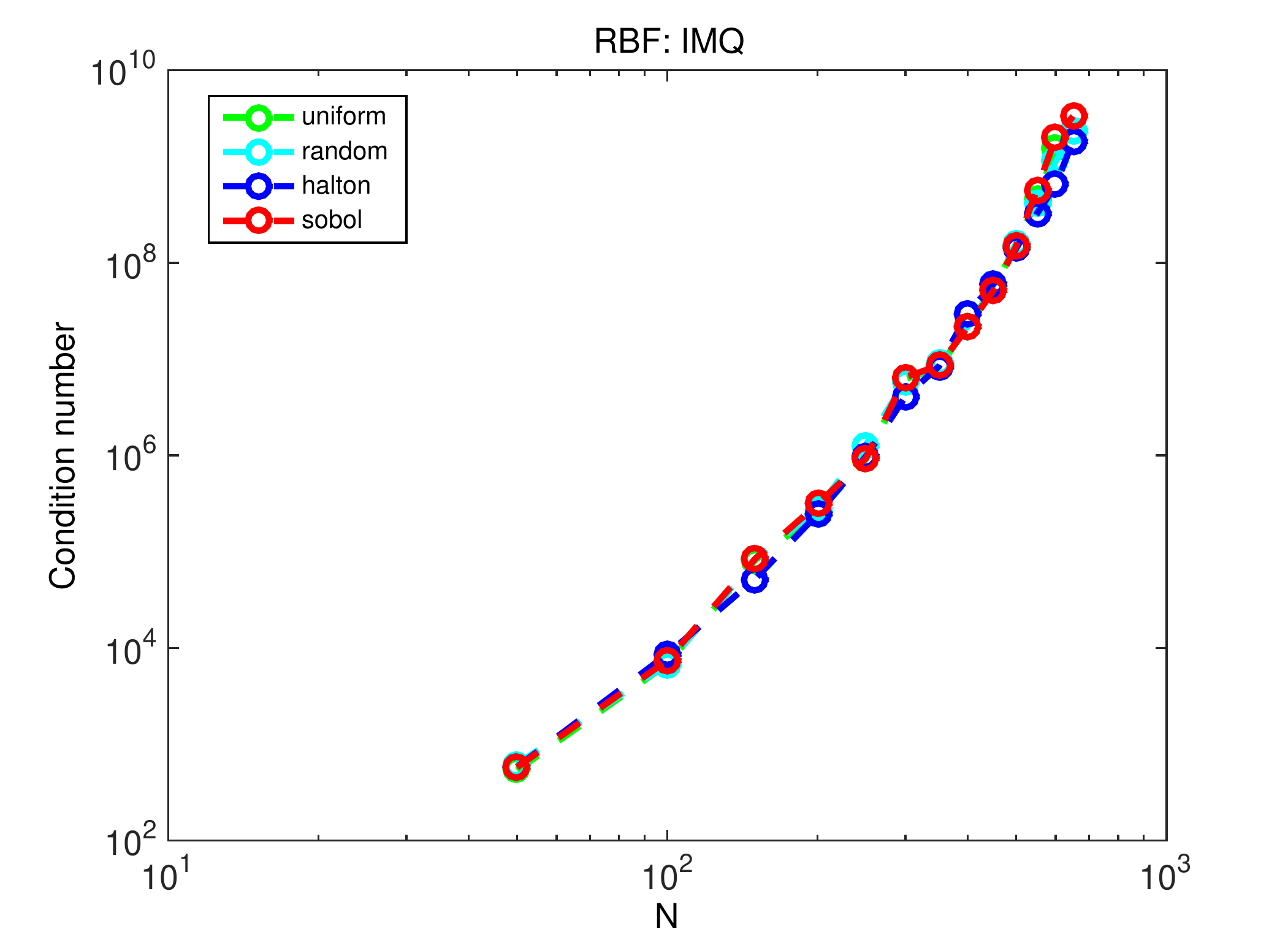}
     \includegraphics[width=6cm]{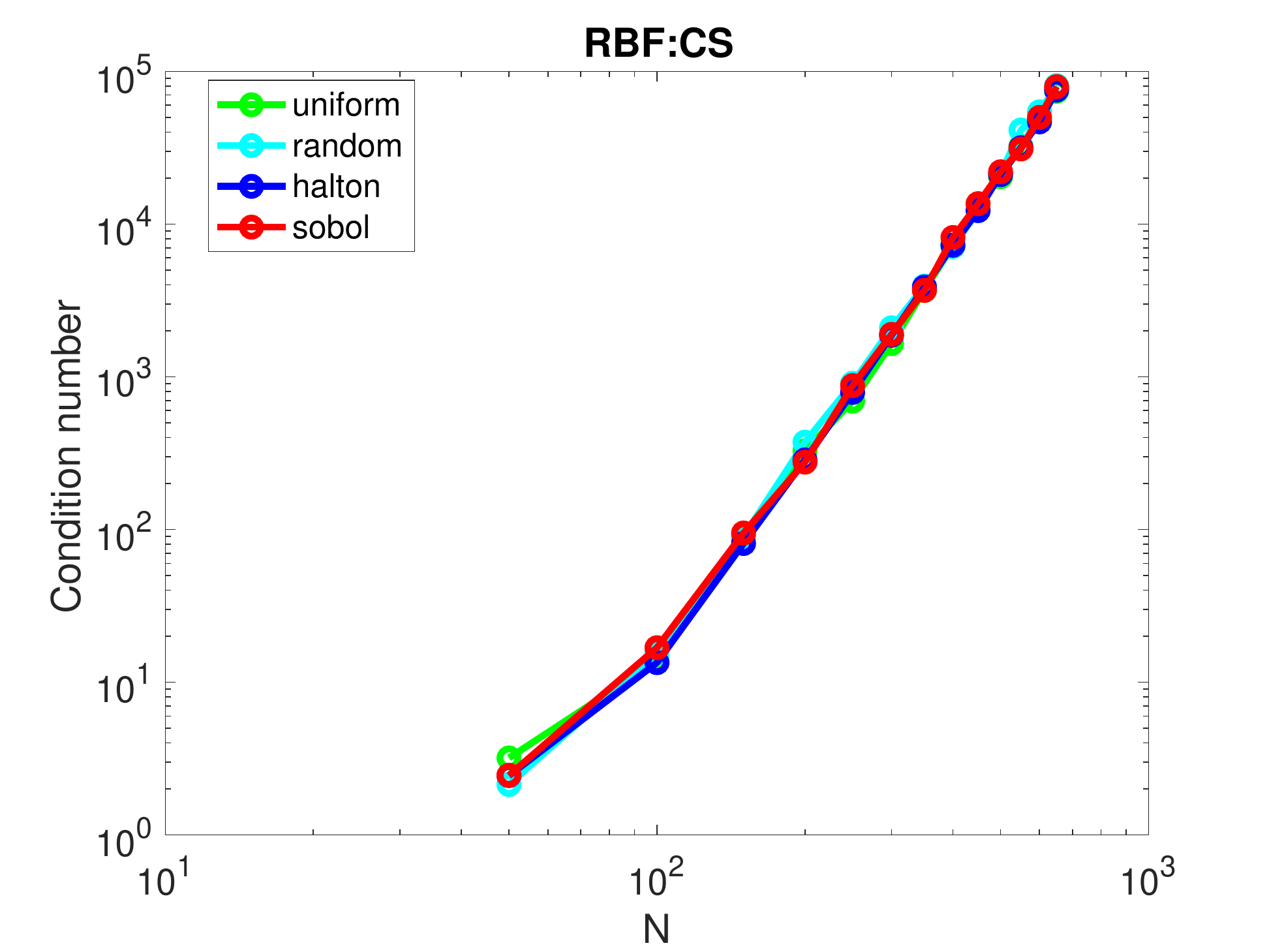}
 \includegraphics[width=6cm]{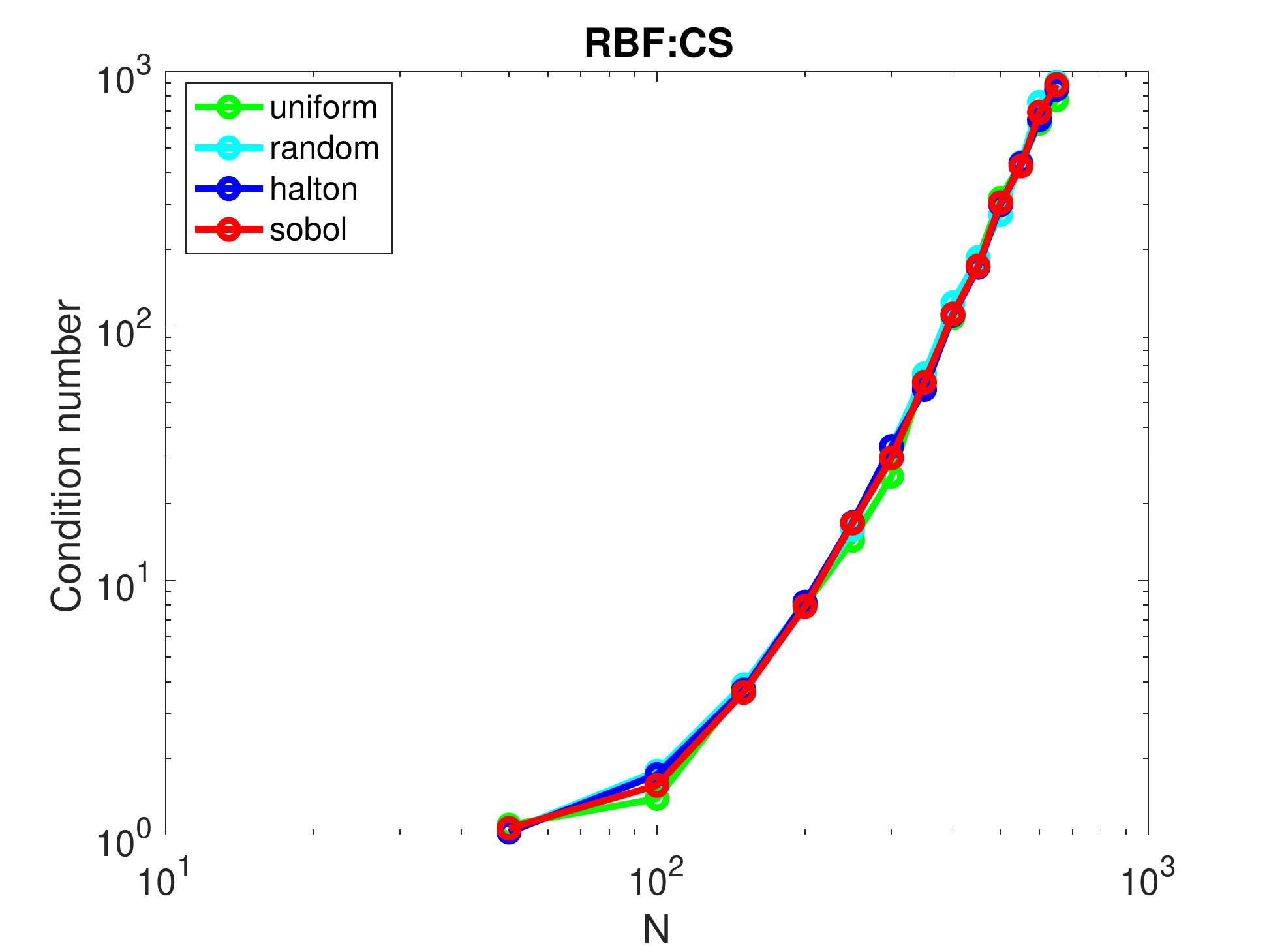}
\end{center}
  \caption{Condition number with respect to the number of sample points $N$ for different choice of the candidates. Left: $\epsilon=3$; Right: $\epsilon=5$.
  \label{fig:d2_gau_N_ca}
    }
\end{figure}

We first investigate the condition number $Cond(\mb{A})=\frac{\sigma_{max}(\mb{A})}{\sigma_{min}(\mb{A})}$ of the design matrix $\mb{A}$.   In  Fig. \ref{cond:ga_imq_ep_chol}, we show the condition number with respect to the different shape parameters.   In  Fig. \ref{cond:ga_imq_N_chol}, we show the condition number with respect to the number of the selected RBF centers $N$.  For both Gaussian and IMQ, our algorithm is much more stable compared to the other sampling methods.  In fact, the other choice of the candidates points do not dramatically affect the performance of the proposed method.  Fig.\ref{fig:d2_gau_N_ca}  show the corresponding results. In the figure, ``uniform" utilizes the evenly spaced points in $I_Z$.

\subsubsection{Numerical accuracy}
In this section, we will compare Cholseky, random, sobol and halton in terms of their ability to approximate test function.

We first consider the benchmark Franke's function $u(z), z=(z_1,z_2)\in[0,1]^2$, given by
 \begin{eqnarray} \label{frankefun}
 \begin{split}
   u(z)=& \frac{3}{4}e^{-((9 z^{(1)}-2)^2+(9 z^{(2)}-2)^2)/4}+\frac{3}{4}e^{-(9 z^{(1)} +1)^2/49-(9 z^{(2)} + 1)^2)/10}\\
   &+\frac{1}{2}e^{-((9 z^{(1)} -7)^2+(9 z^{(2)} -3)^2)/4}-\frac{1}{5}e^{-(9 z^{(1)} -4)^2-(9 z^{(2)} -7)^2}.
\end{split}
 \end{eqnarray}

\begin{figure}[htbp]
\begin{center}
   \includegraphics[width=6cm]{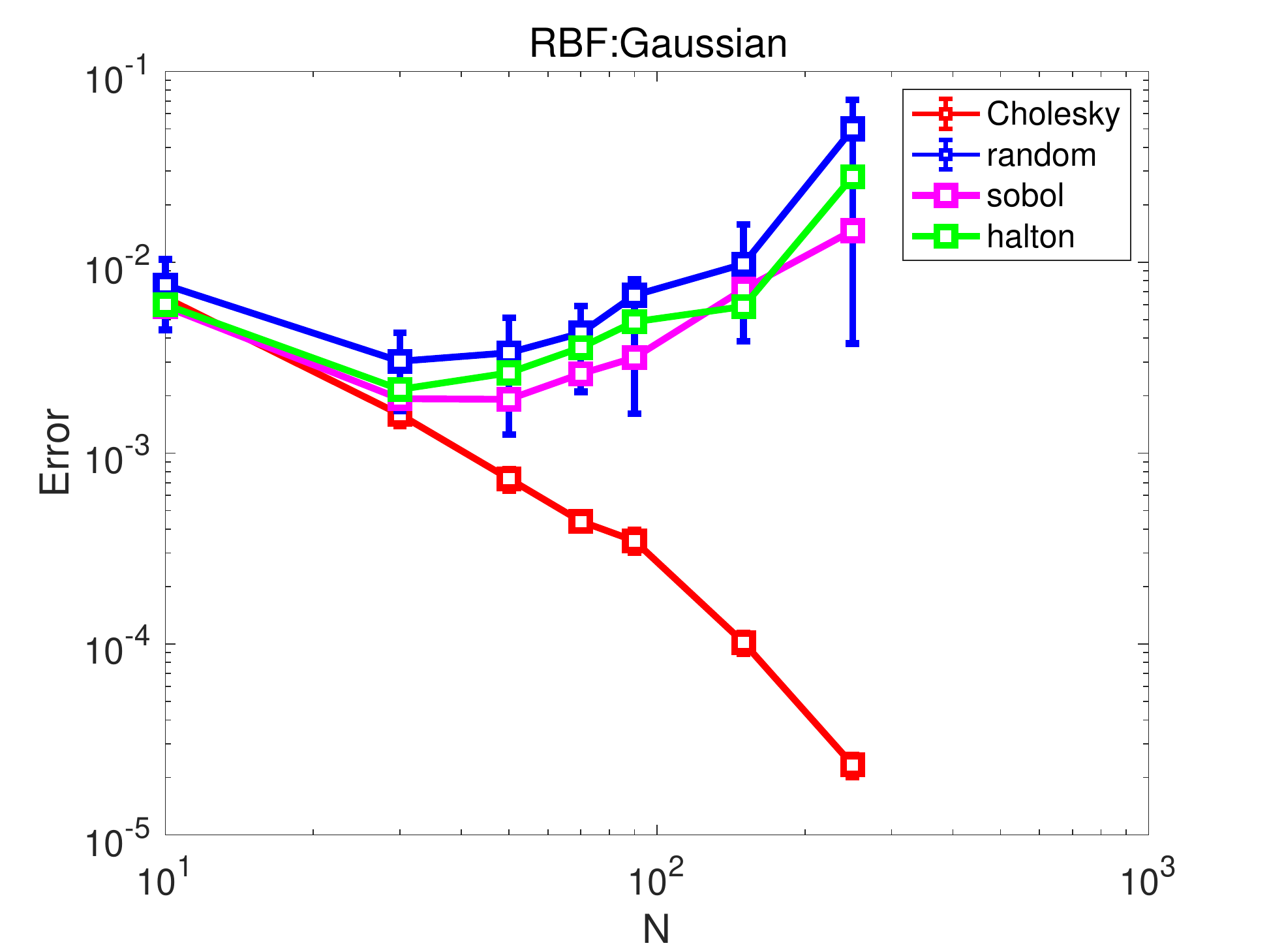}
   \includegraphics[width=6cm]{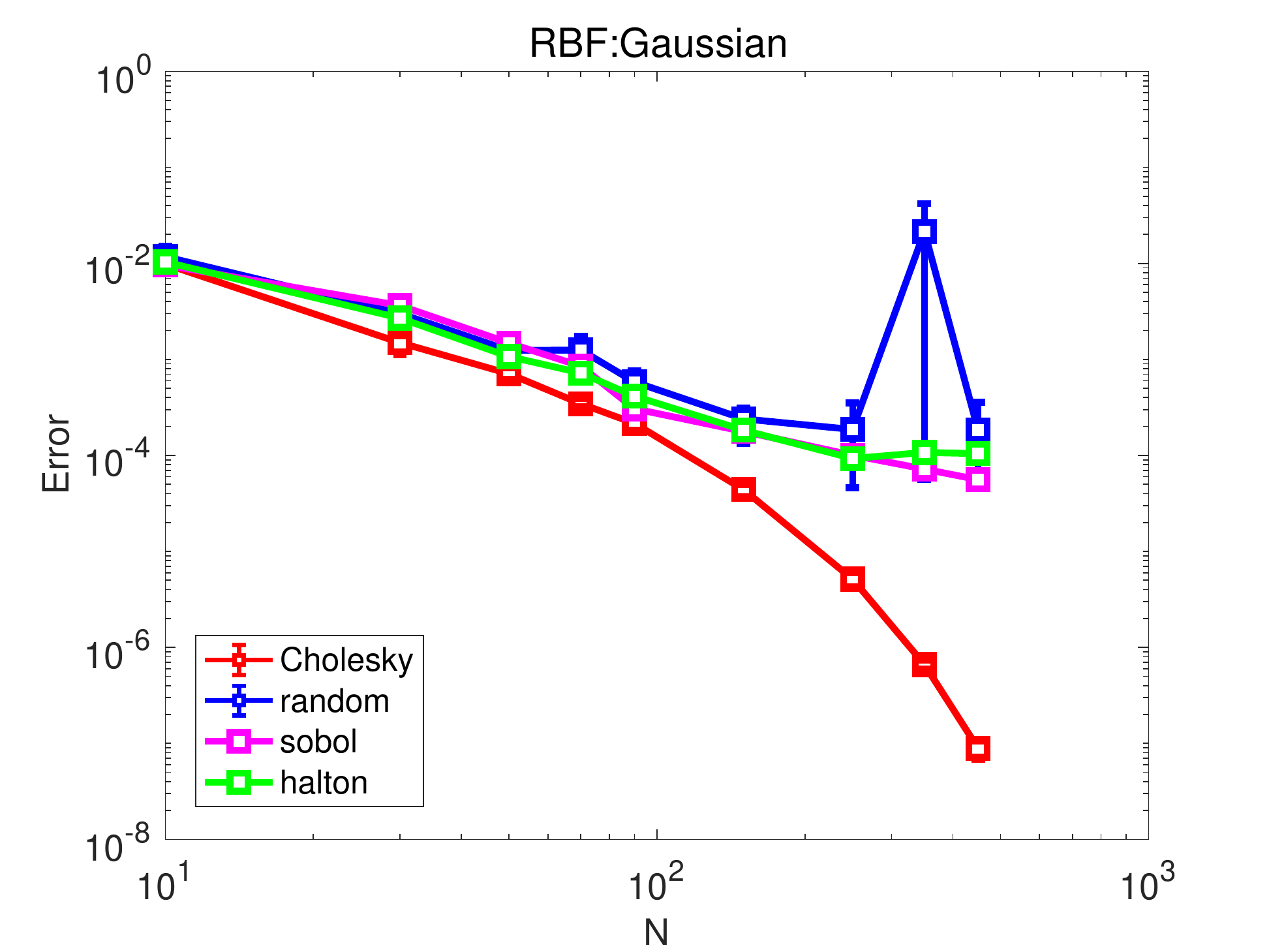}
        \includegraphics[width=6cm]{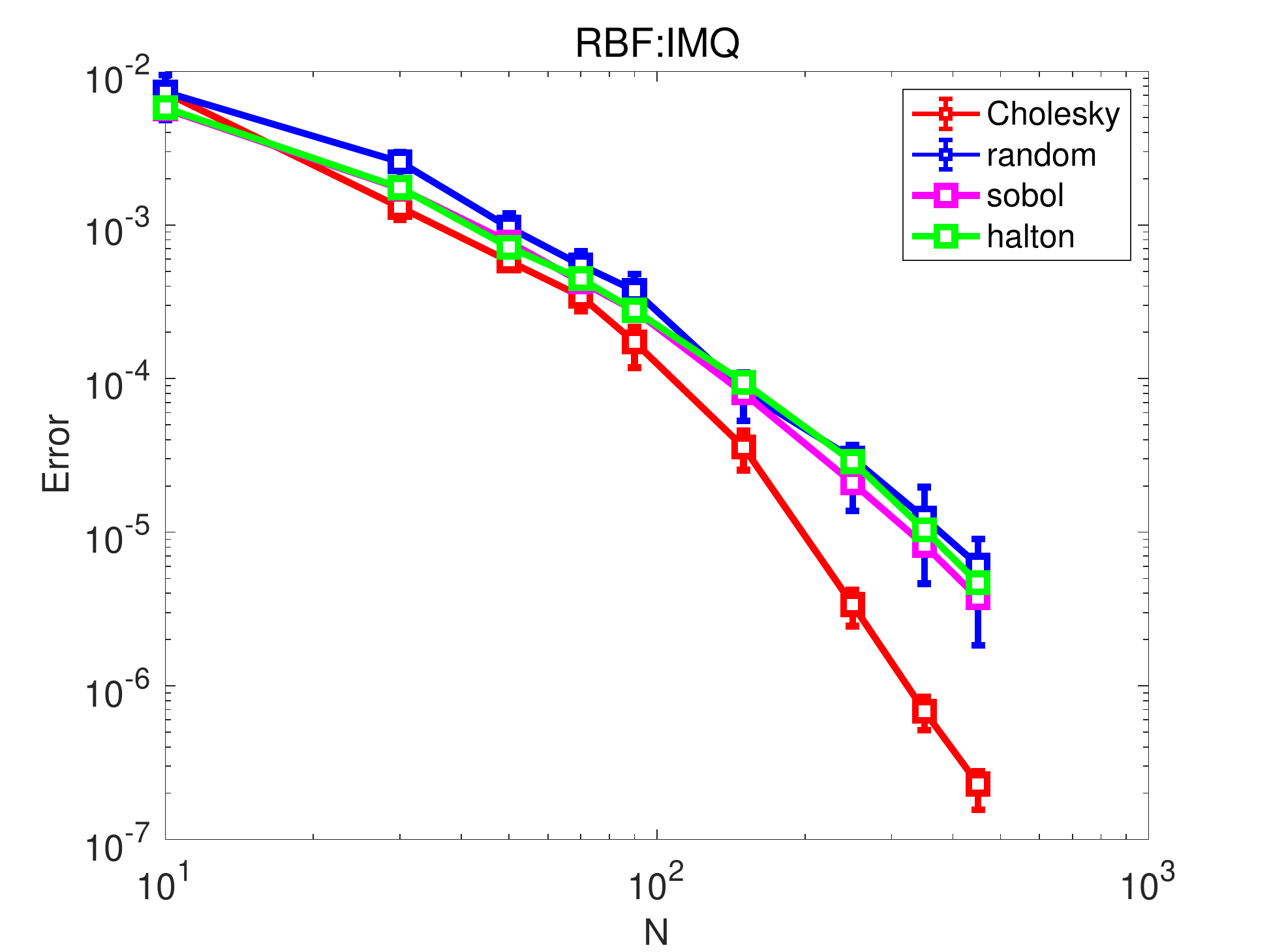}
  \includegraphics[width=6cm]{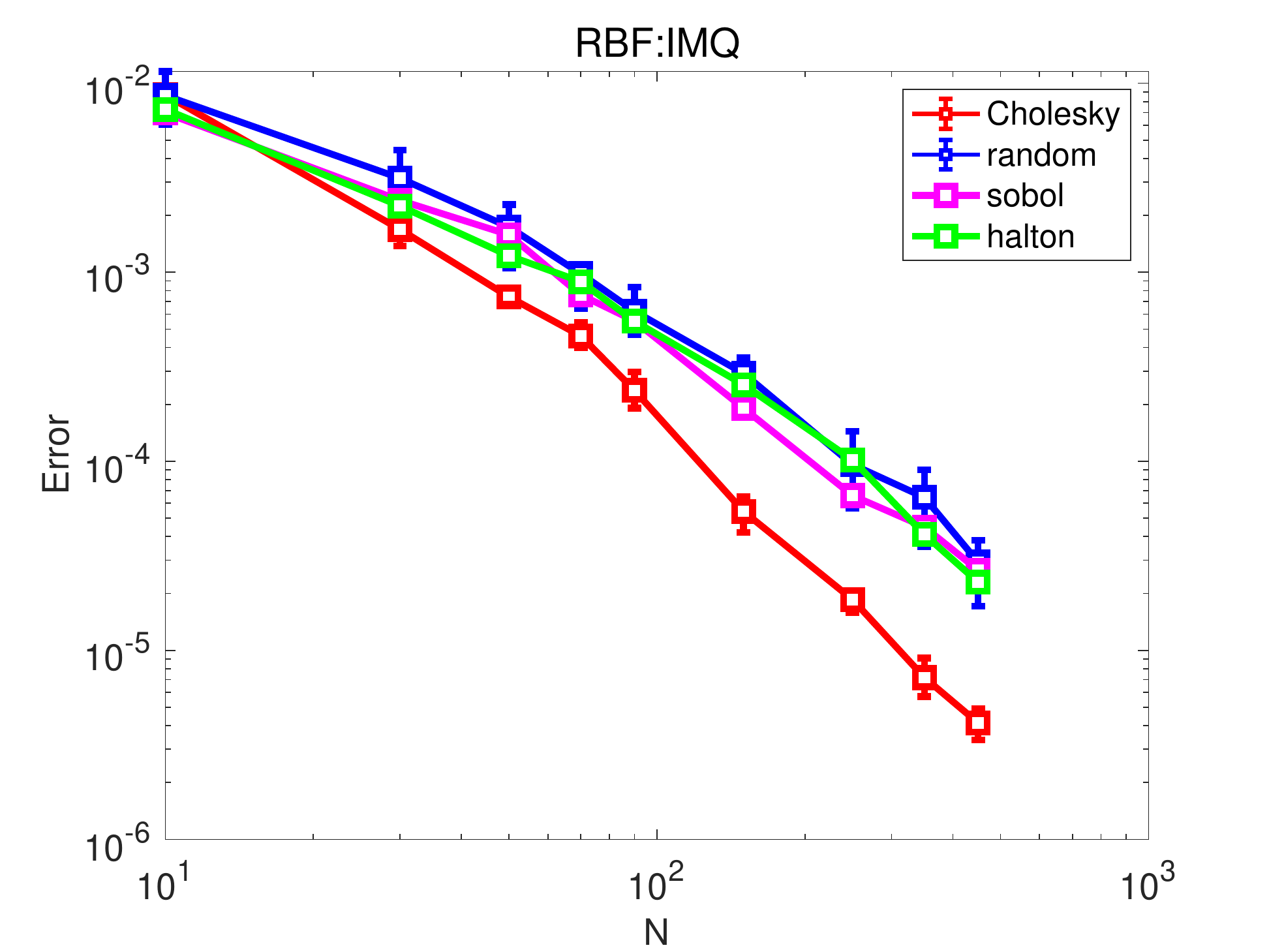}
          \includegraphics[width=6cm]{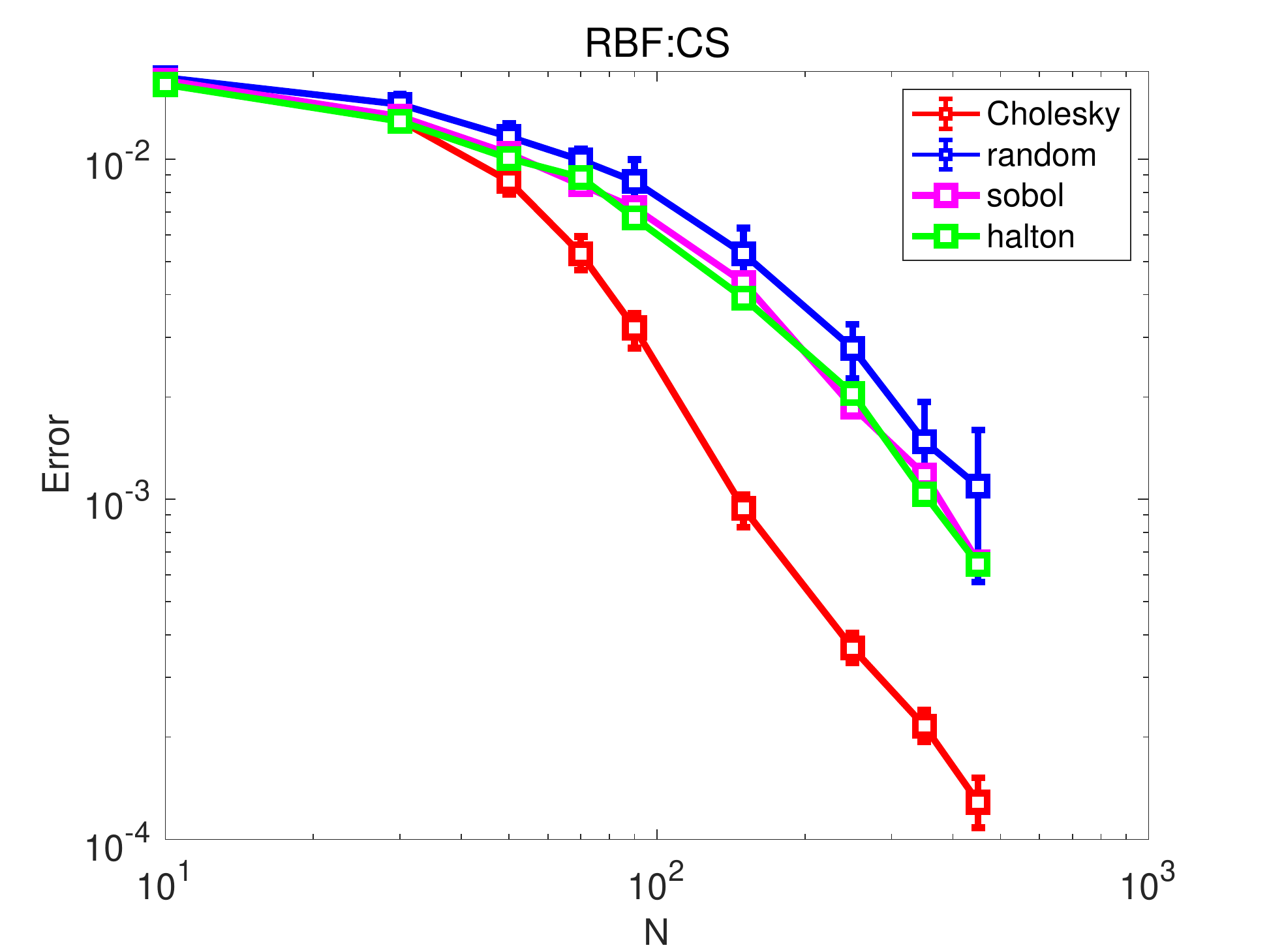}
  \includegraphics[width=6cm]{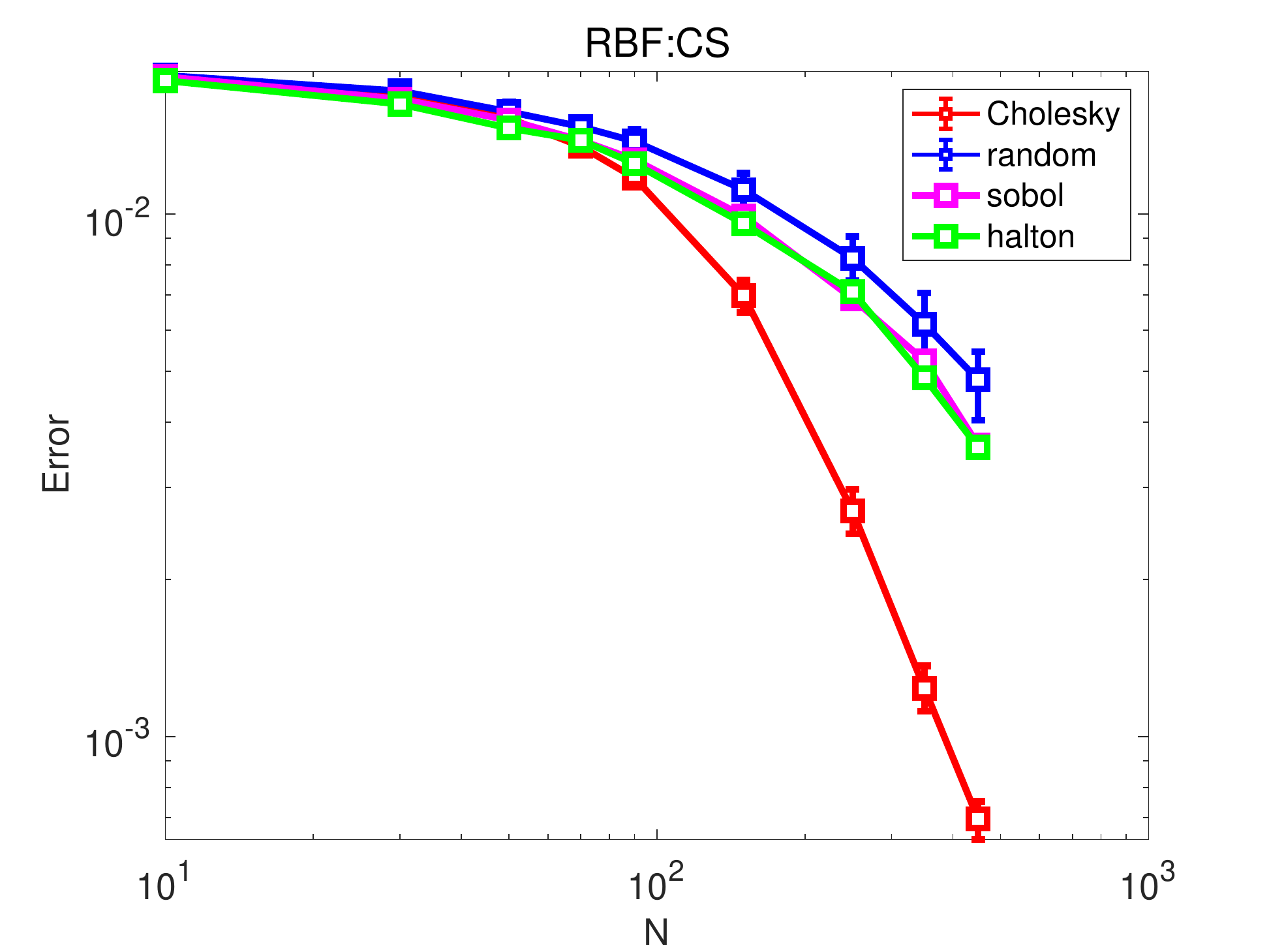}
\end{center}
  \caption{ Approximation error with respect to the number of sample points $N$ using Gaussian, IMQ and CS. Left: $\epsilon=3$; Right: $\epsilon=5$.
  \label{fig:d2_gau_N}
    }
\end{figure}

\begin{figure}[htbp]
\begin{center}
   \includegraphics[width=6cm]{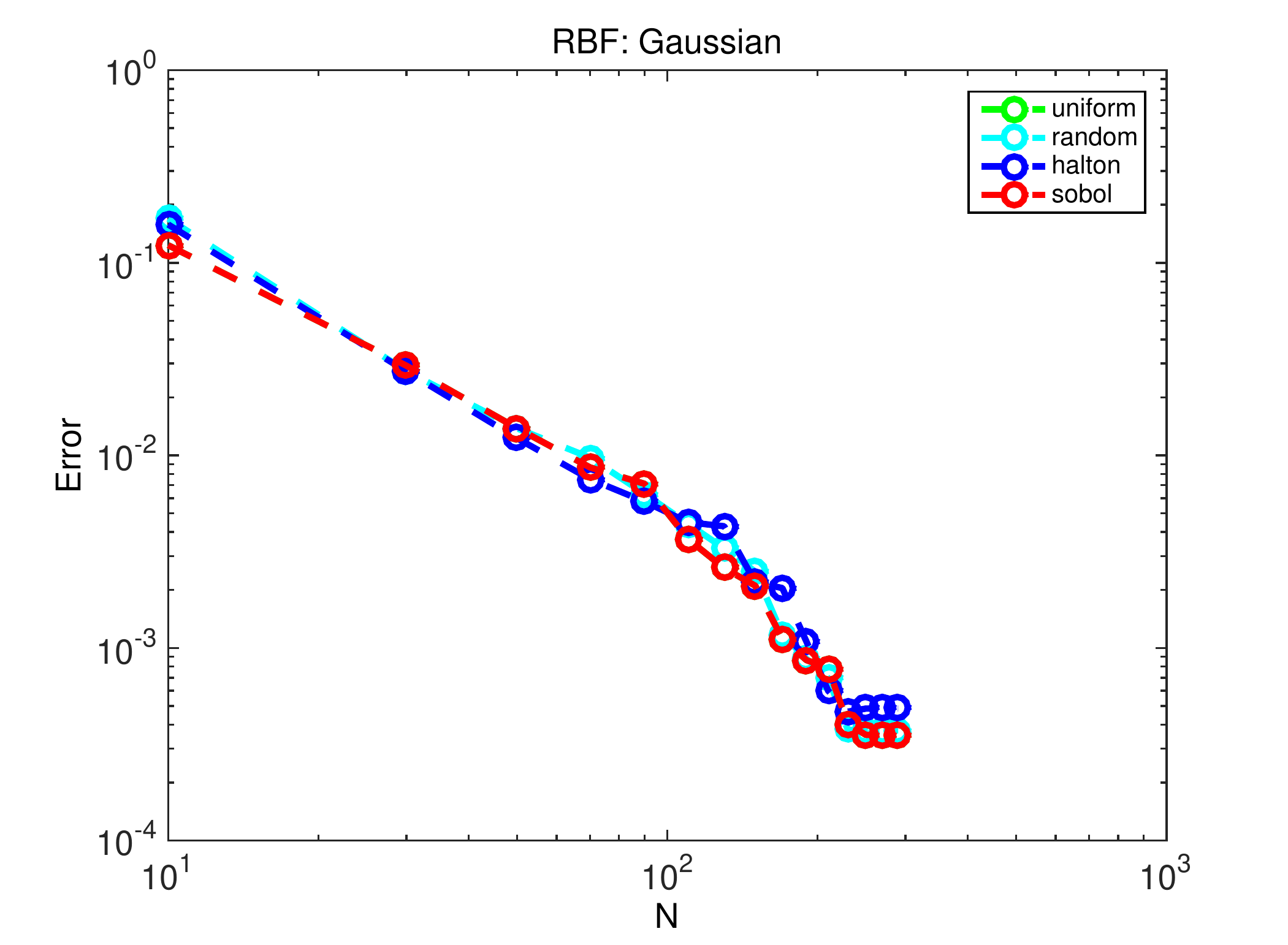}
   \includegraphics[width=6cm]{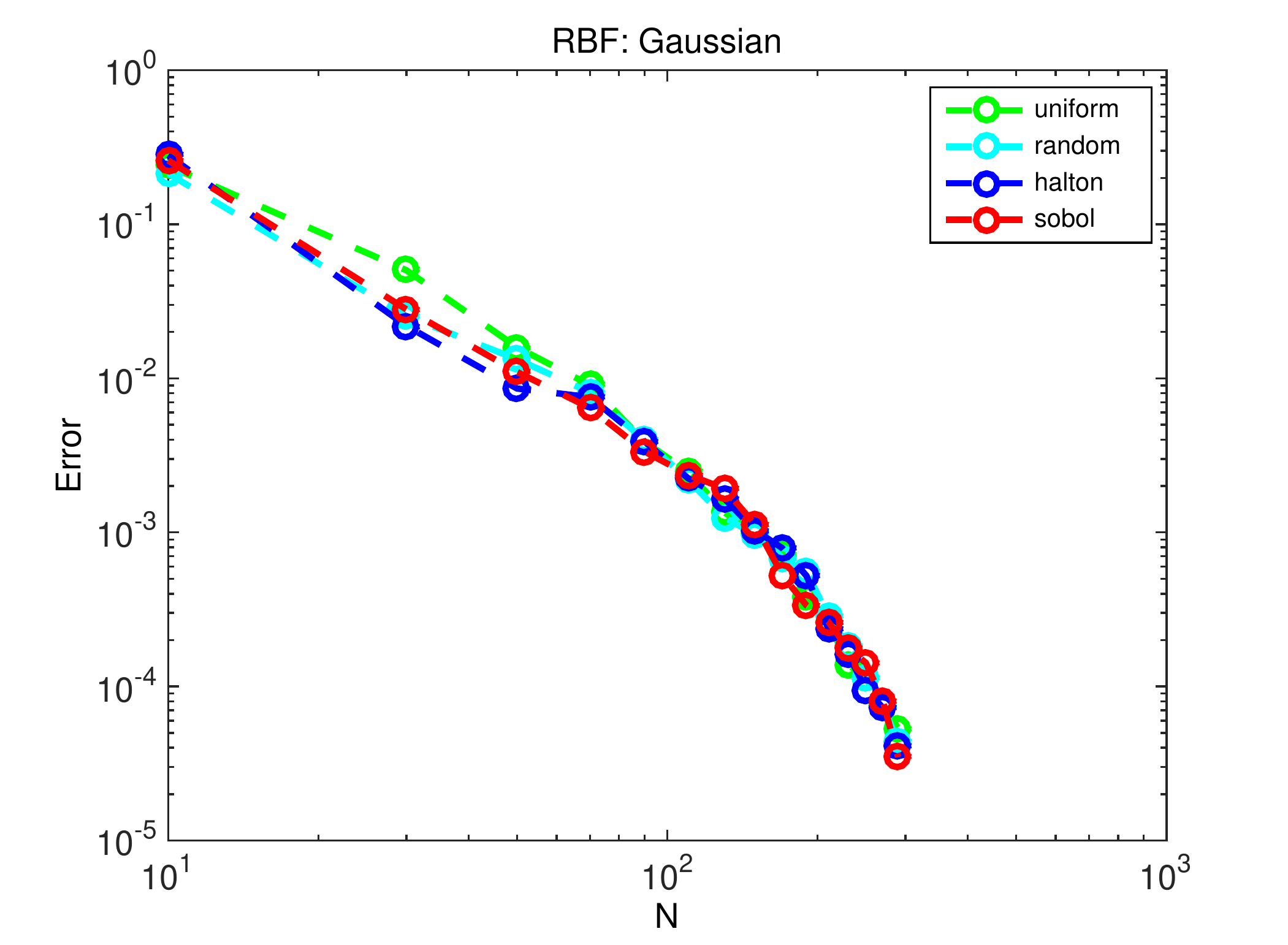}
  \includegraphics[width=6cm]{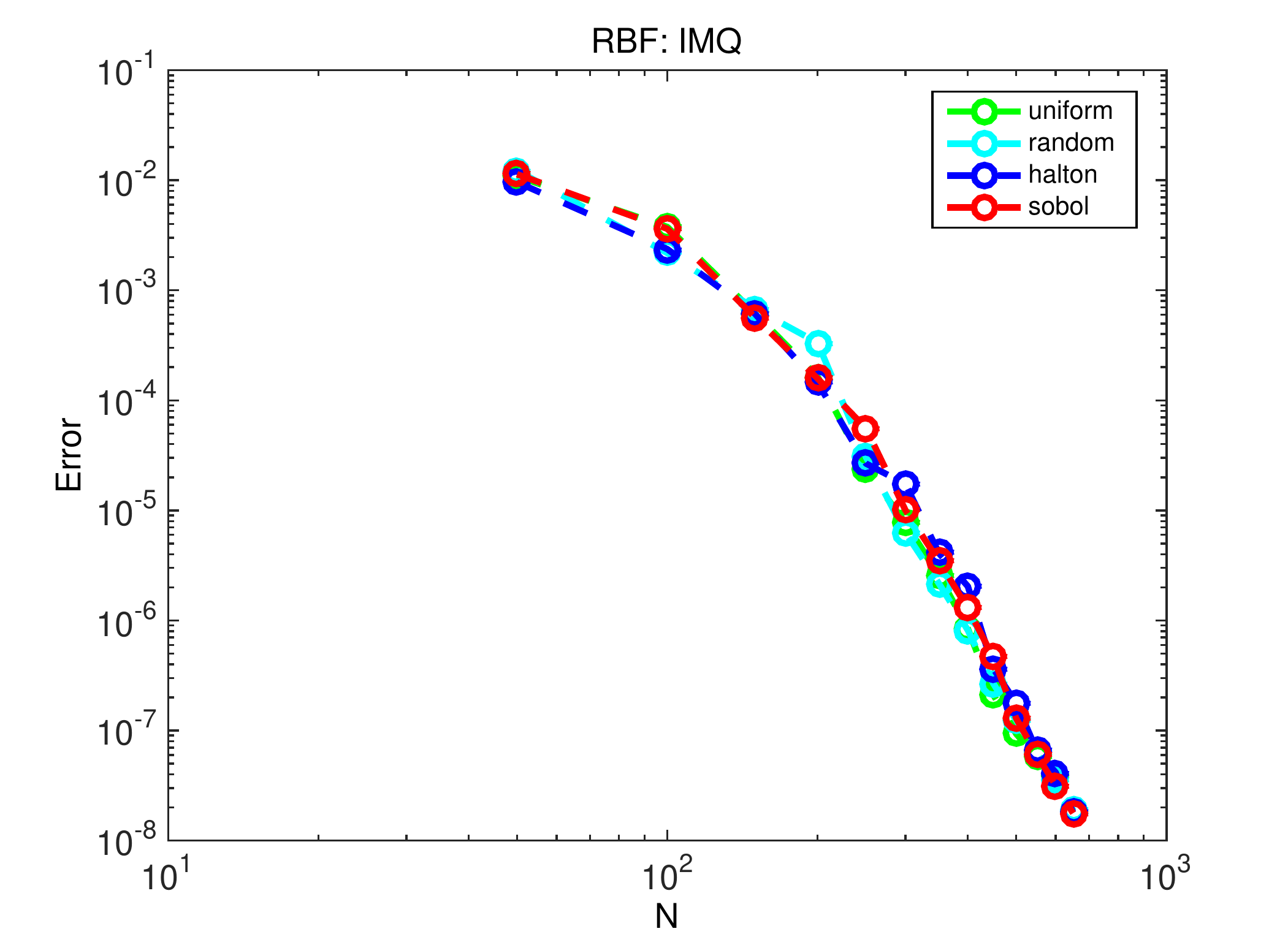}
  \includegraphics[width=6cm]{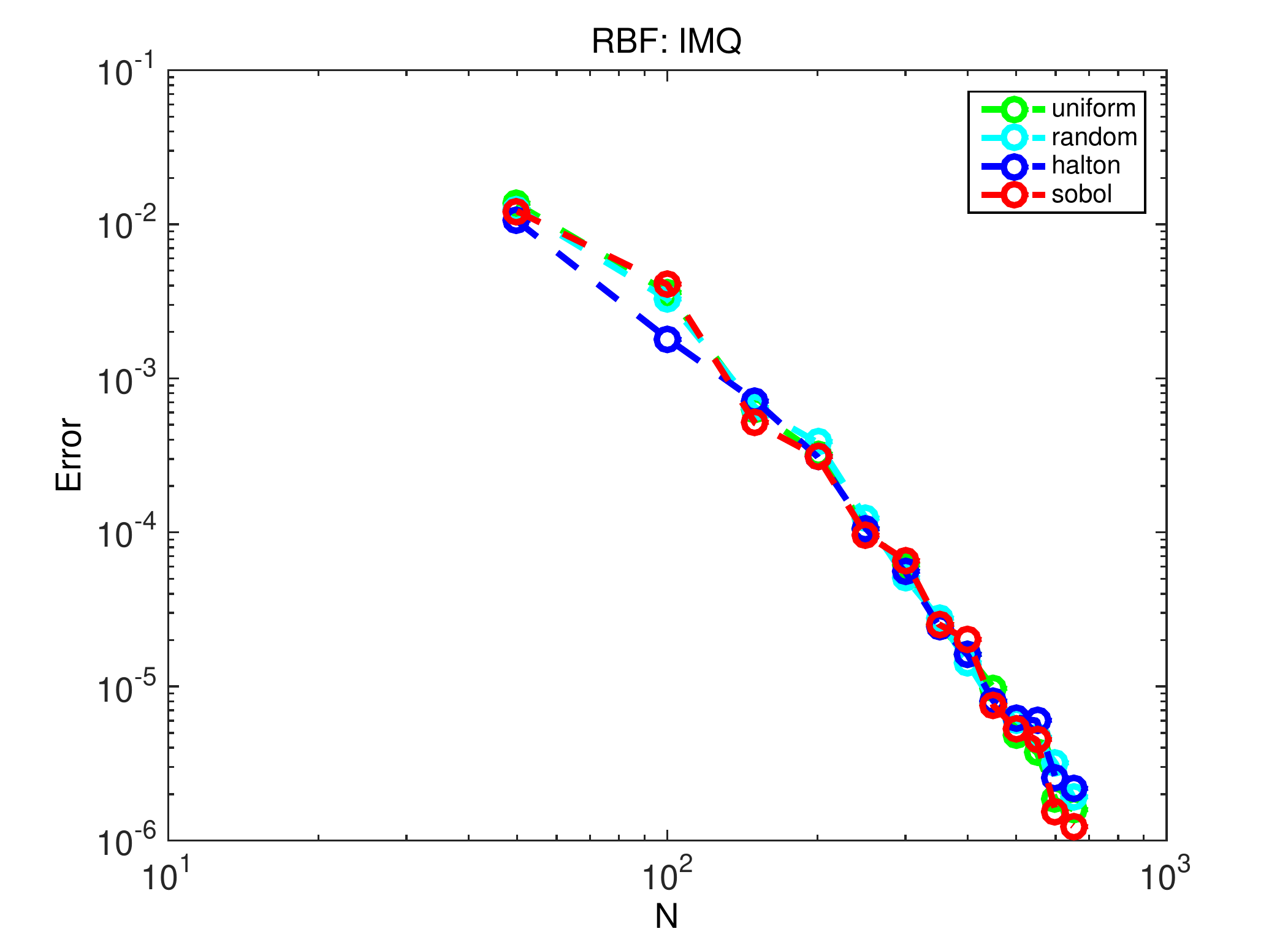}
  \includegraphics[width=6cm]{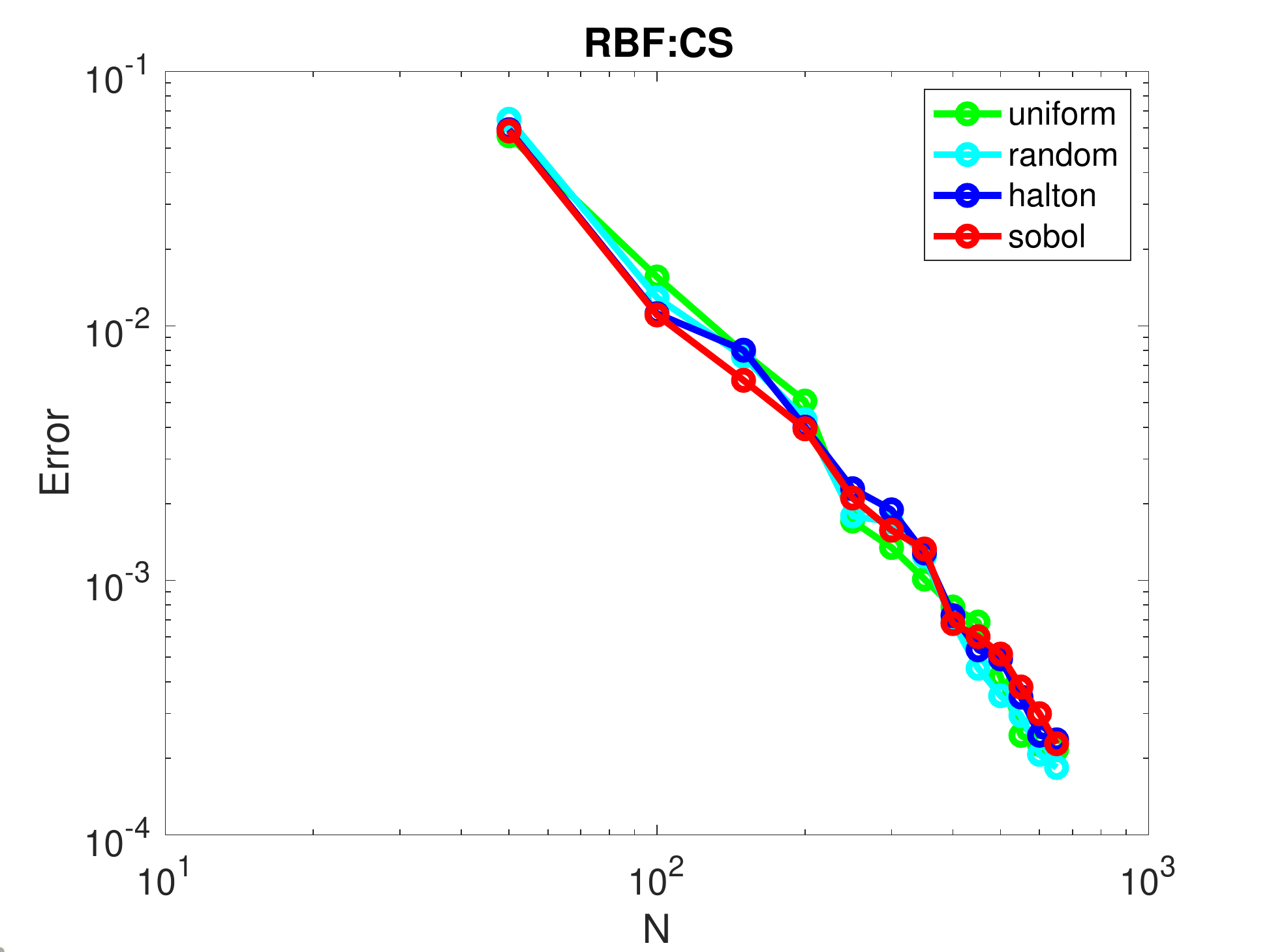}
  \includegraphics[width=6cm]{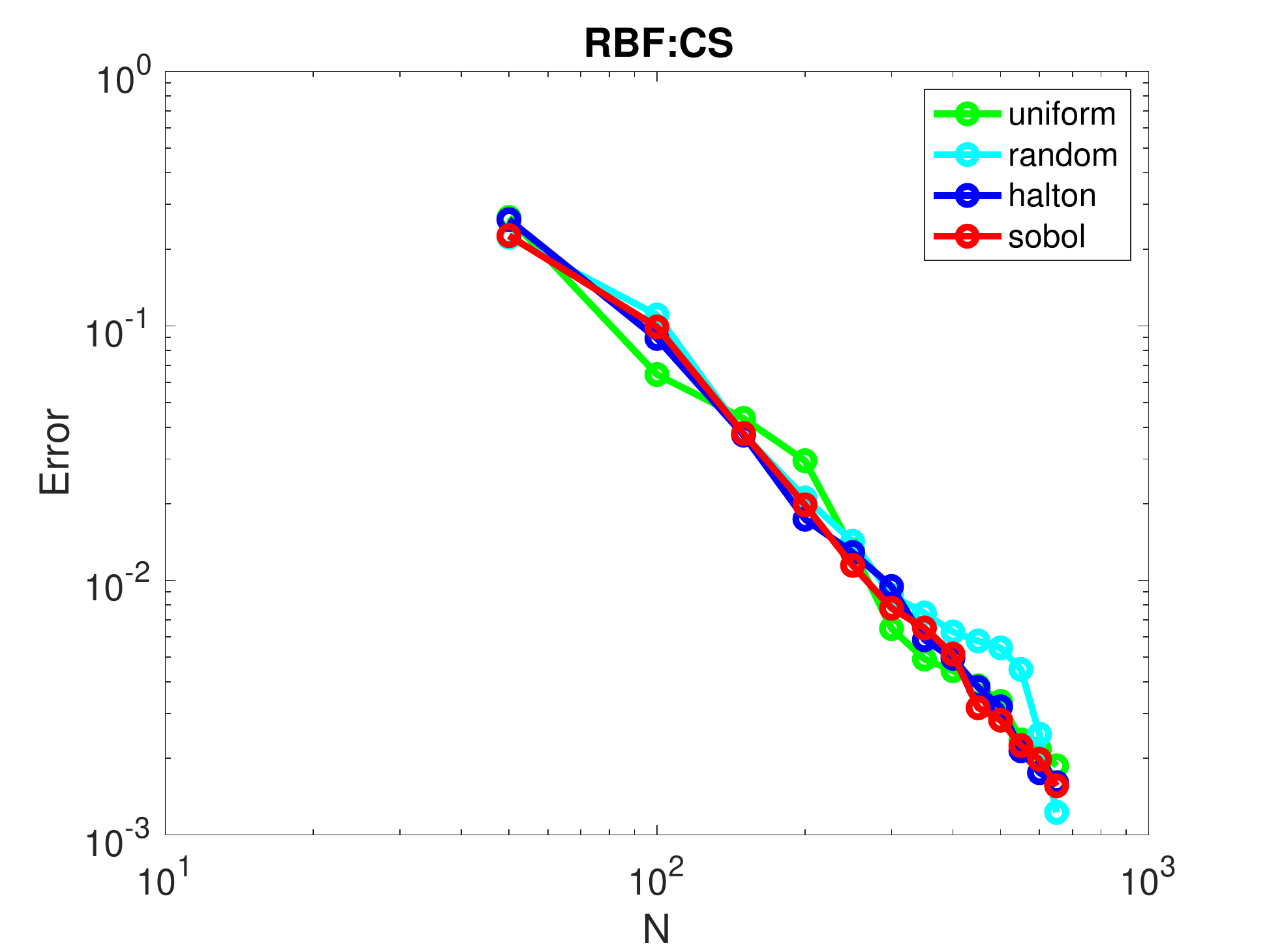}
\end{center}
  \caption{Numerical results with respect to the number of sample points $N$ for different choice of the candidates: Left: $\epsilon=3$; Right: $\epsilon=5$.
    \label{fig:d2_imq_N_ca}
    }
\end{figure}

The corresponding numerical errors for the test function in (\ref{frankefun}) are shown in Fig. \ref{fig:d2_gau_N}. We observe that  our proposed algorithm  produces superior results than the other sampling methods.   Moreover, the new algorithm shows clear convergence patterns as $N$ increases.  This is considered to be an improvement over the other sampling methods: it is clear in the error profile of the other sampling methods that providing more sampling points does not always result in better accuracy. Again, the  choice of the candidates points do not dramatically affect the performance of the proposed method.  Fig.\ref{fig:d2_imq_N_ca}  show the corresponding results.

\begin{figure}[htbp]
\begin{center}
  \includegraphics[width=6cm]{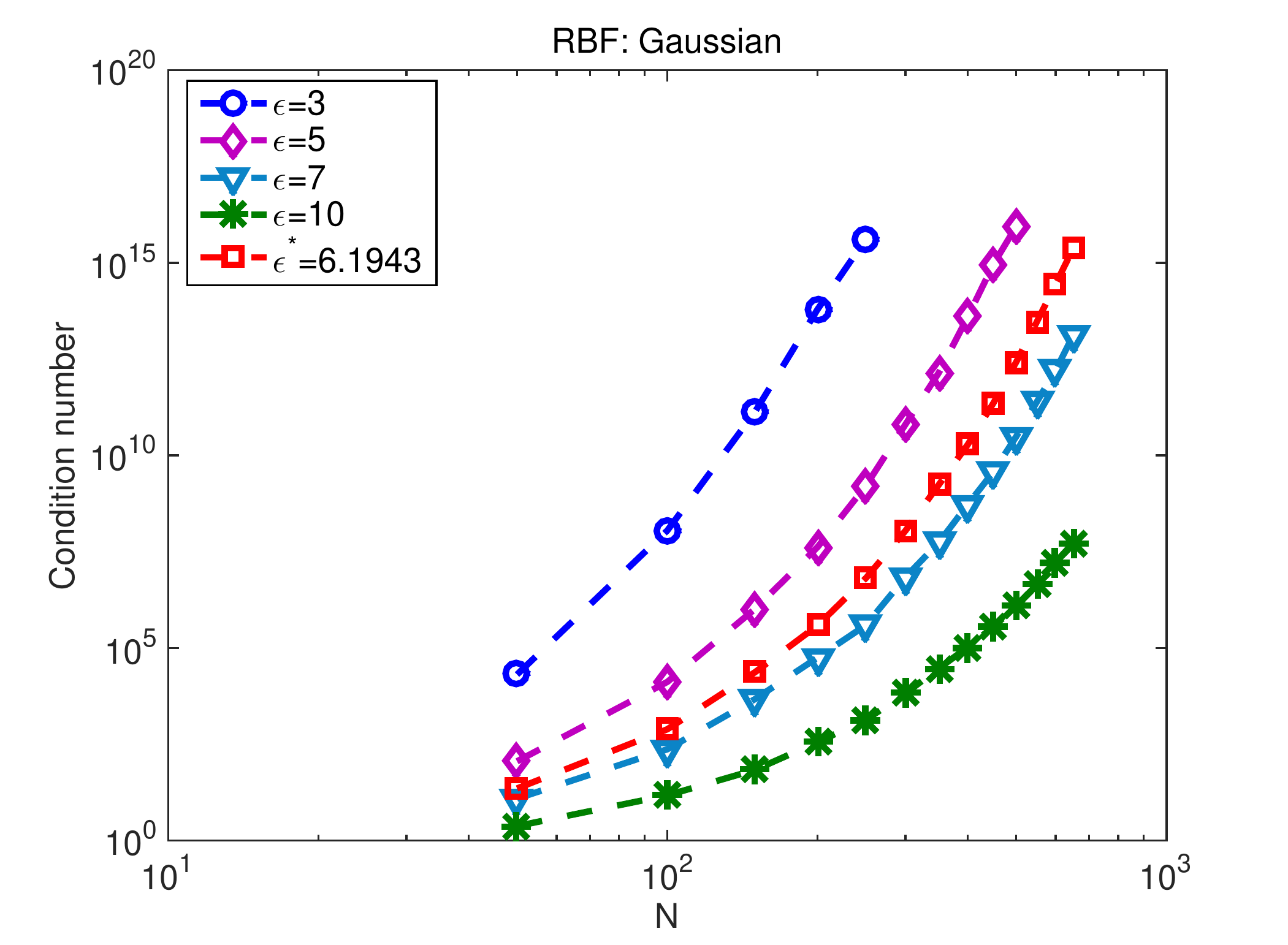}
  \includegraphics[width=6cm]{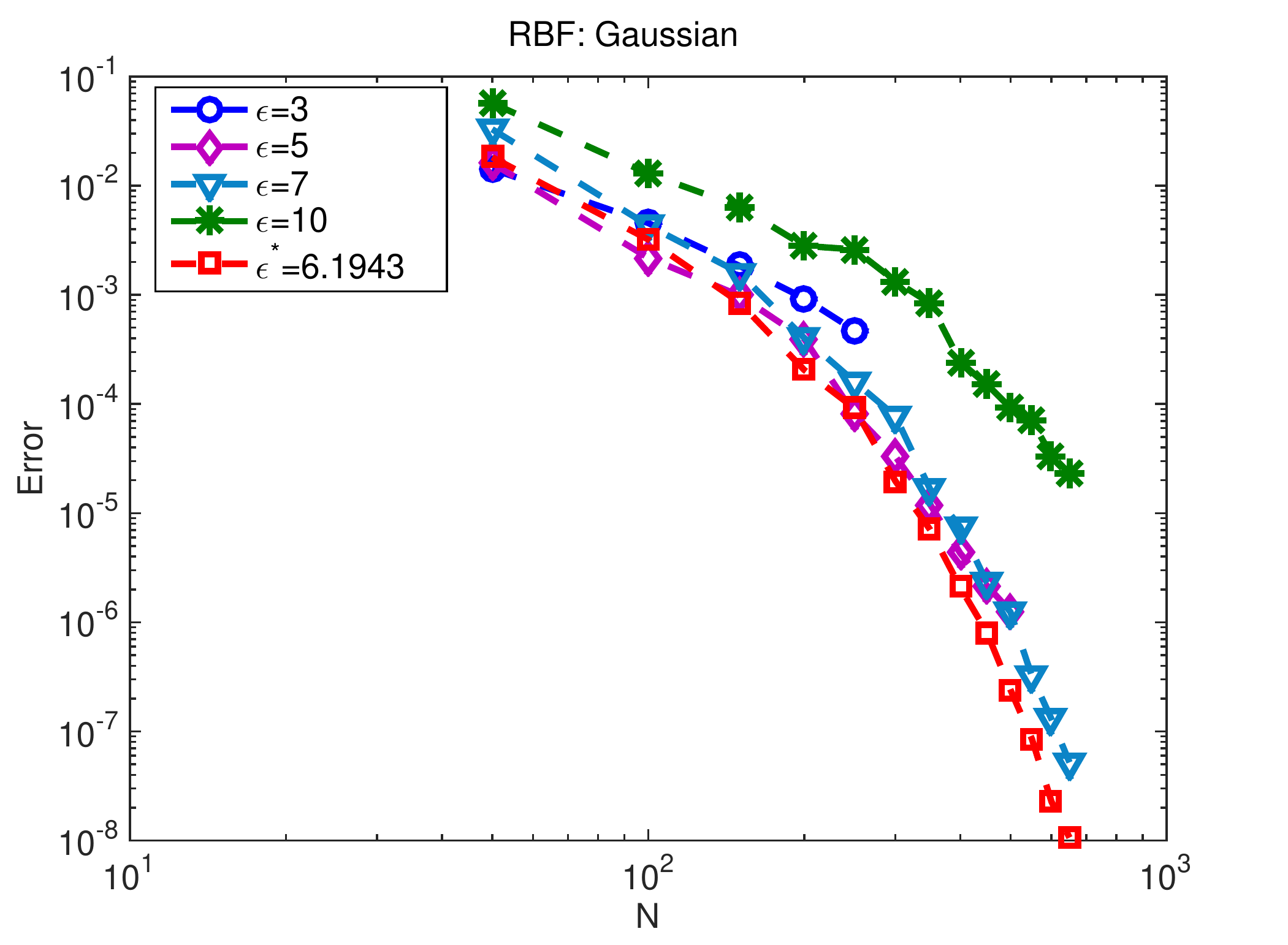}
     \includegraphics[width=6cm]{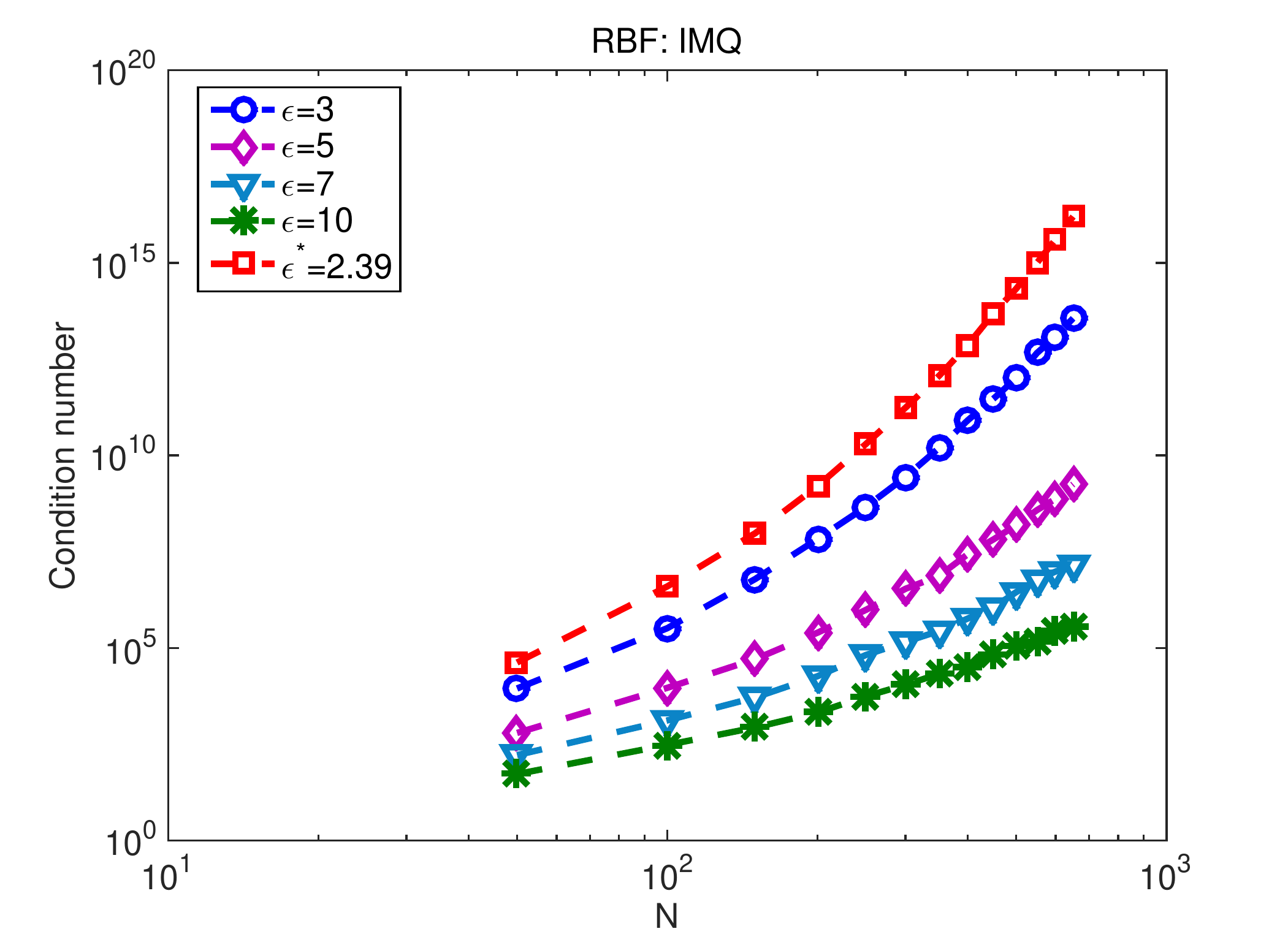}
   \includegraphics[width=6cm]{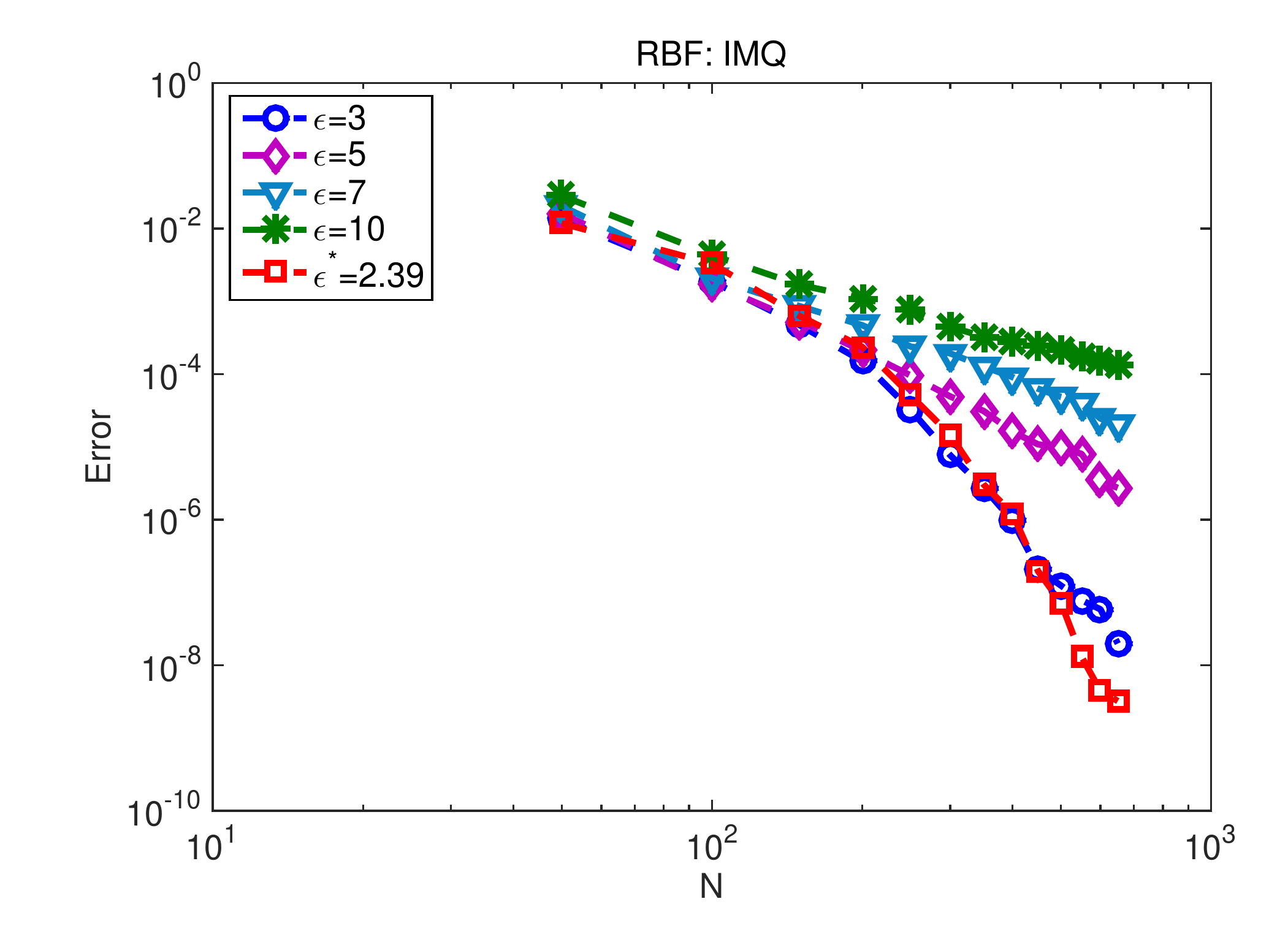}
     \includegraphics[width=6cm]{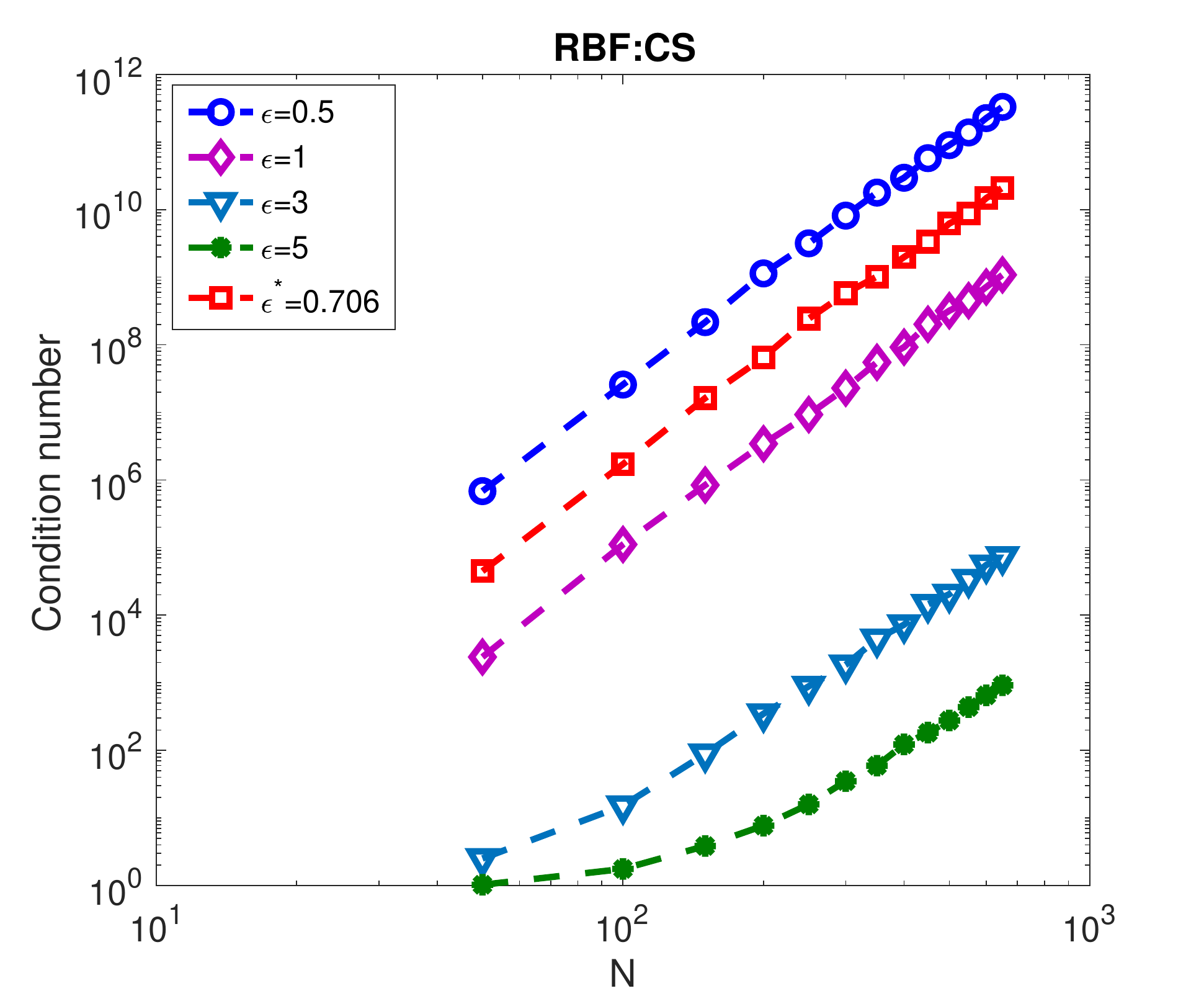}
   \includegraphics[width=6cm]{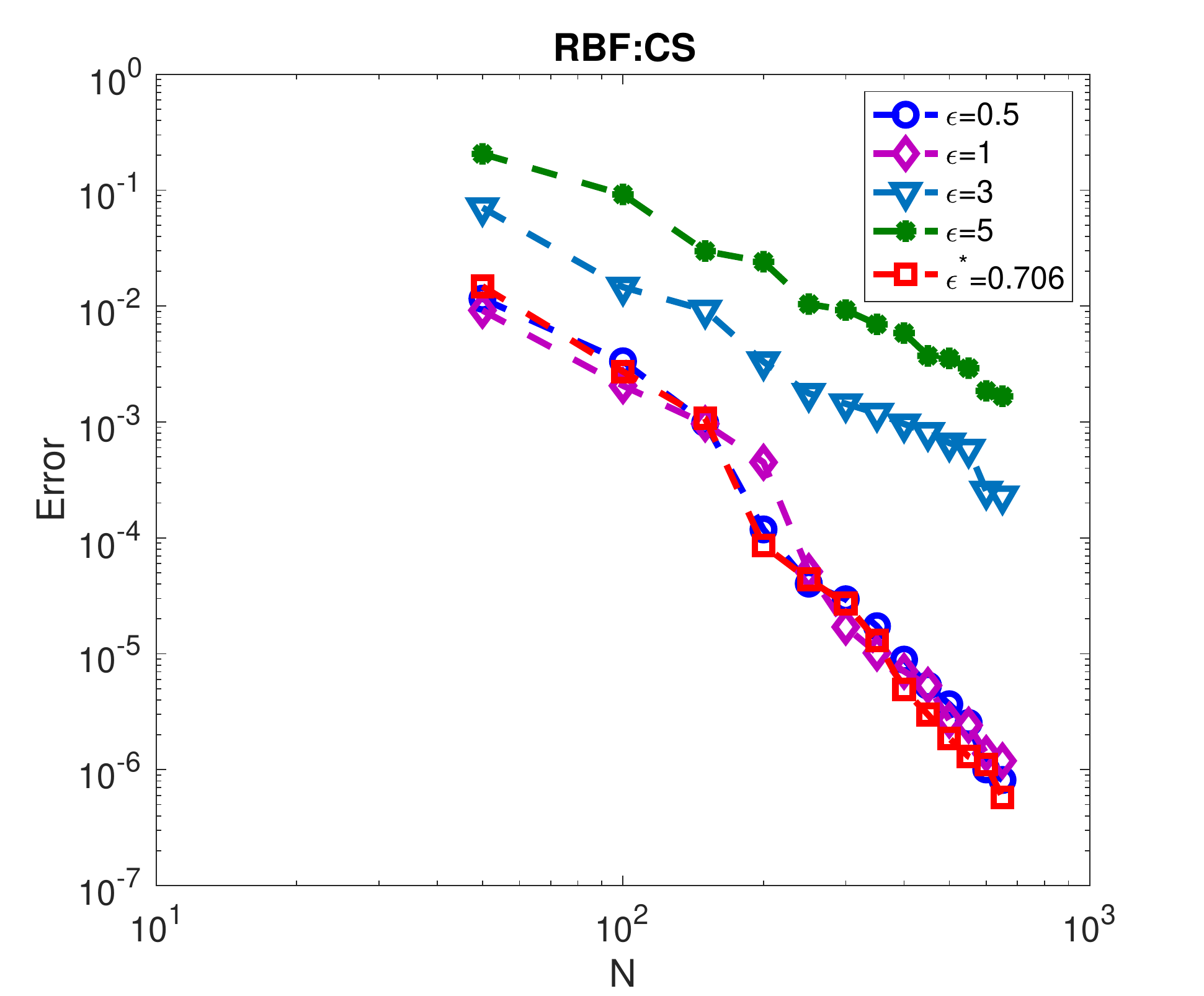}
\end{center}
  \caption{Numerical results with respect to the number of sample points $N$ for different  shape parameters. Left: condition number; Right: errors.
   \label{fig:d2_gau_c}
    }
\end{figure}

In Fig.\ref{fig:d2_gau_c} we show the numerical results as a function of the sample points, $N$, obtained using the Cholesky algorithm  with various values of the shape parameter $\epsilon$.  The left plot shows the condition number while the right plot shows the corresponding errors. It can be seen from this figure that as $\epsilon$ gets smaller, the resulting matrix becomes ill-conditioned.
The numerical results when the LOOCV algorithm is used are also presented in Fig. \ref{fig:d2_gau_c}. From these figures, we observe that the Cholesky algorithm using LOOCV to choose the good value of shape parameter can provide a very good numerical results.

We now consider the stochastic elliptic equation, one of the most used
benchmark problems in UQ, in one spatial dimension,
 \begin{equation} \label{spde}
 -\frac{d}{dx} [\kappa(x,z)\frac{du}{dx}(x,z)]=f, \quad  (x,z)\in (0,1) \times \mathbb{R}^d,
 \end{equation}
 with boundary conditions
 \begin{eqnarray*}
 u(0,z)=0,  \quad u(1,z)=0,
 \end{eqnarray*}
and $f=2$. The random diffusivity takes the following form
 \begin{equation}\label{randdiff}
 \kappa(x,z)=1+\sigma \sum^d_{k=1} \frac{1}{k^2 \pi^2} \cos(2\pi k x)z^{(k)},
 \end{equation}
 where $z=(z^{(1)},\dots,z^{(d)})$ is a random vector with independent
 and identically distributed components.

\begin{figure}[htbp]
\begin{center}
   \includegraphics[width=4.4cm]{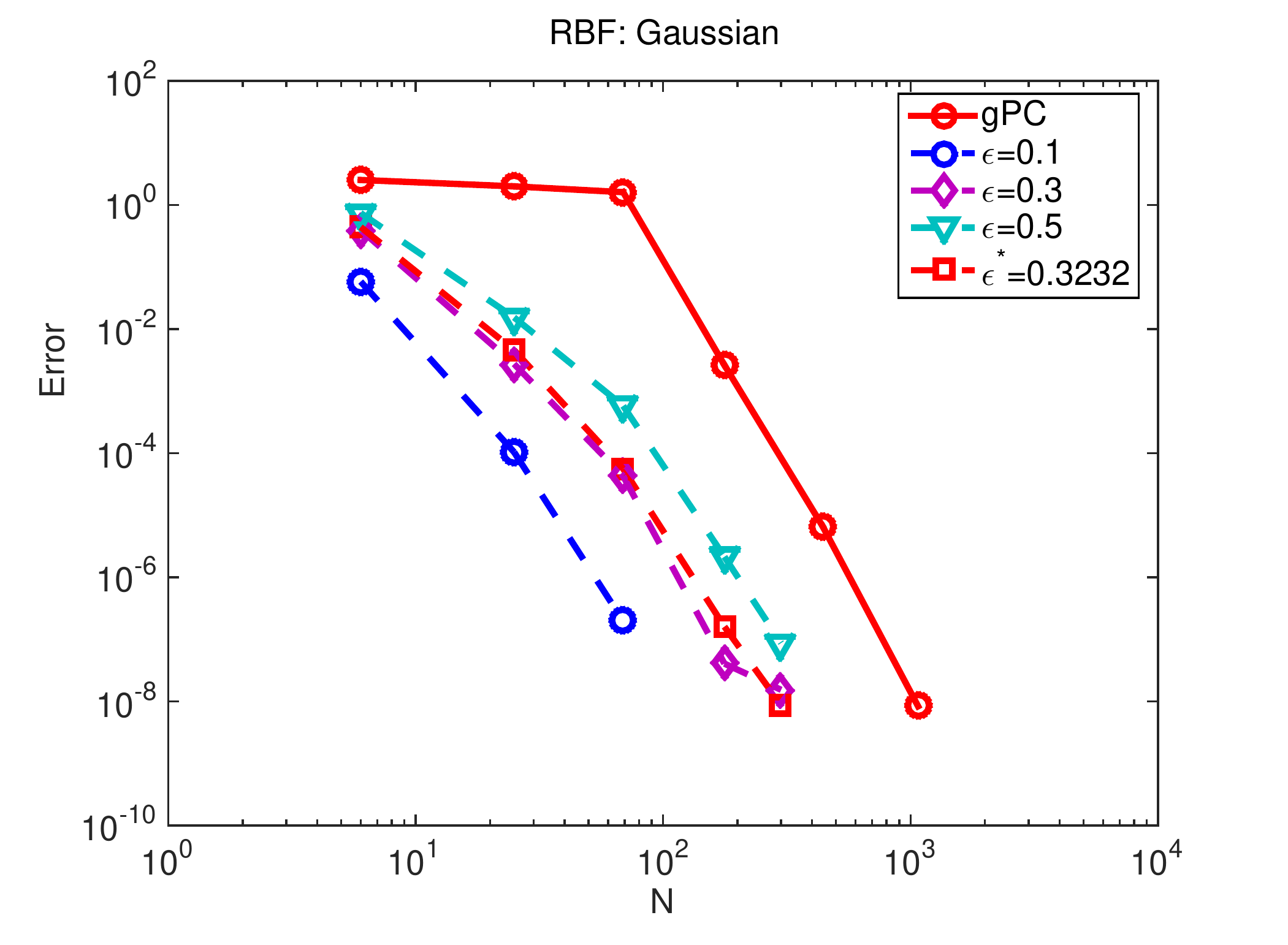}
   \includegraphics[width=4.4cm]{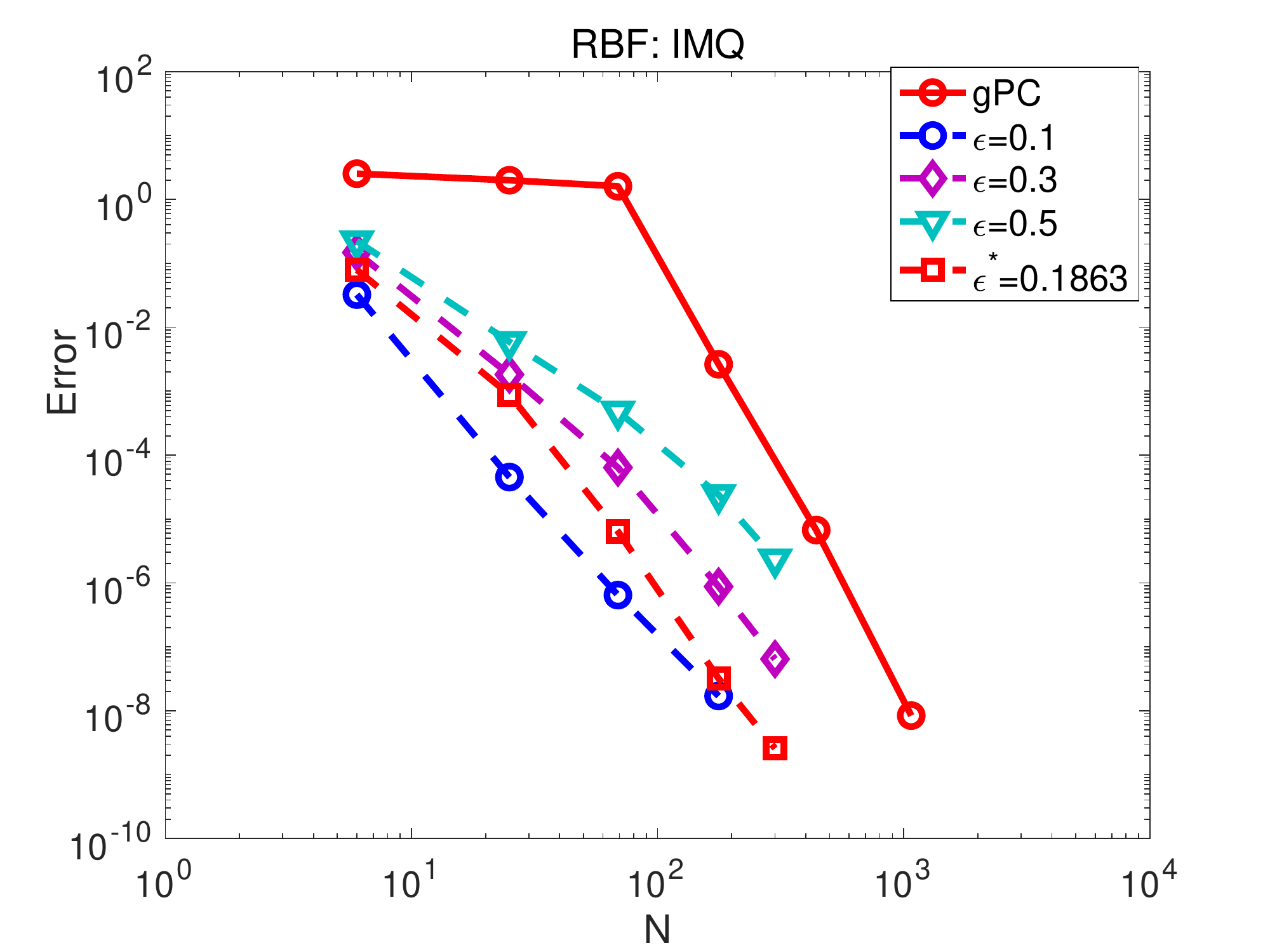}
   \includegraphics[width=4.4cm]{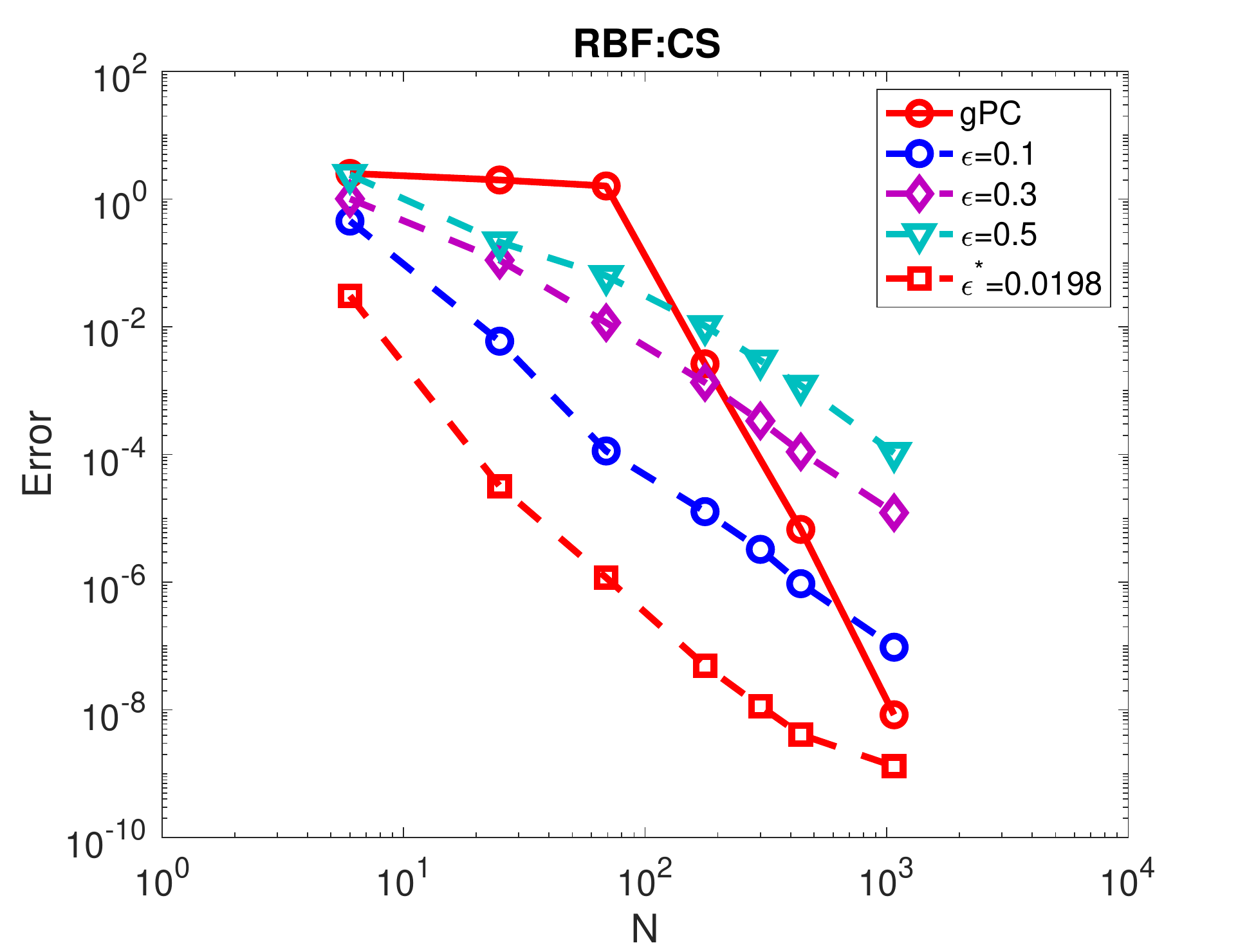}
\end{center}
  \caption{ Approximation error against the number of sample points $N$ (d=3). Left: Gaussian; Middle: IMQ; Right: CS.
    \label{fig:d3_diff}
  }
  \end{figure}

\begin{figure}[htbp]
\begin{center}
   \includegraphics[width=4.4cm]{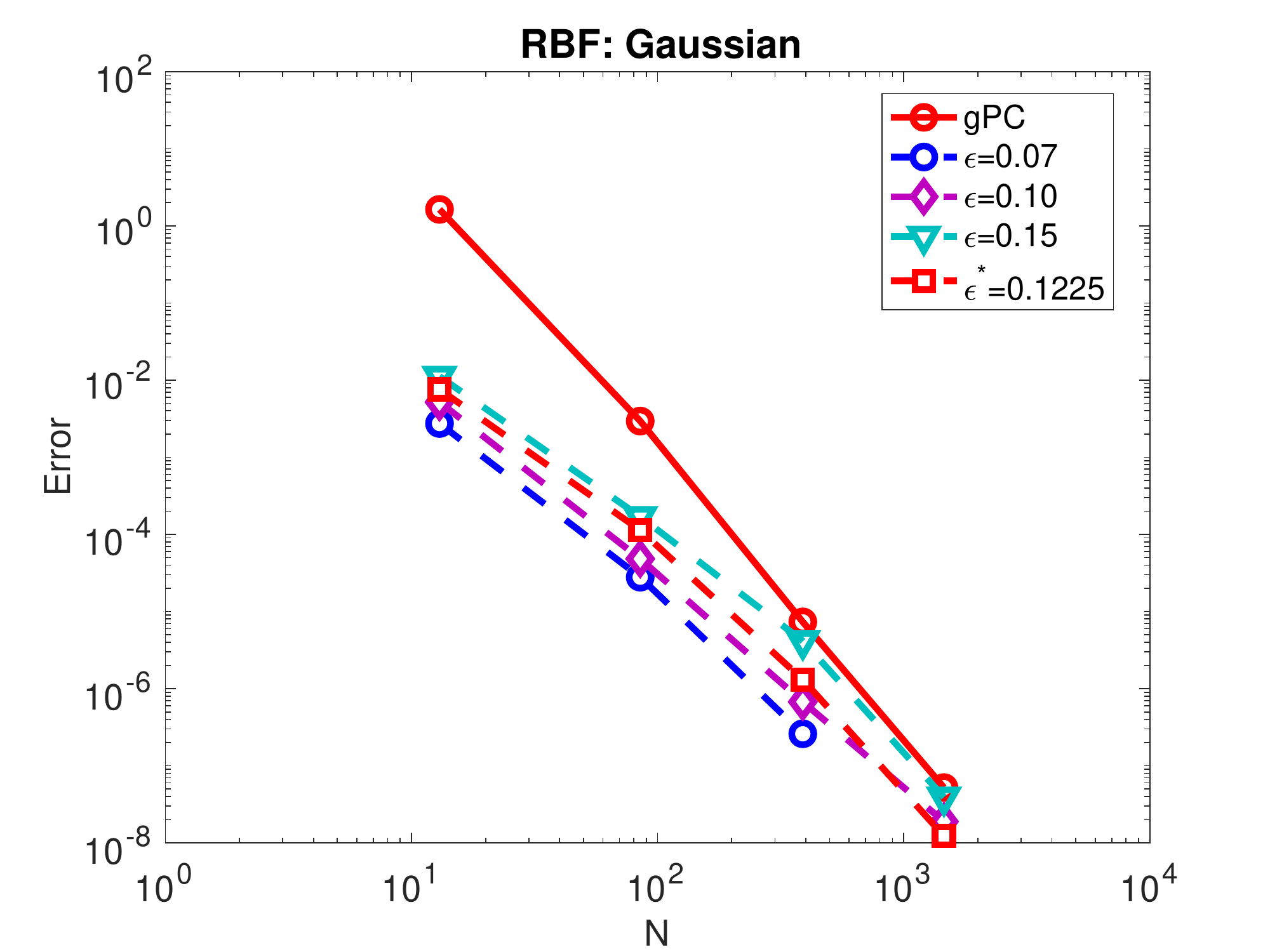}
   \includegraphics[width=4.4cm]{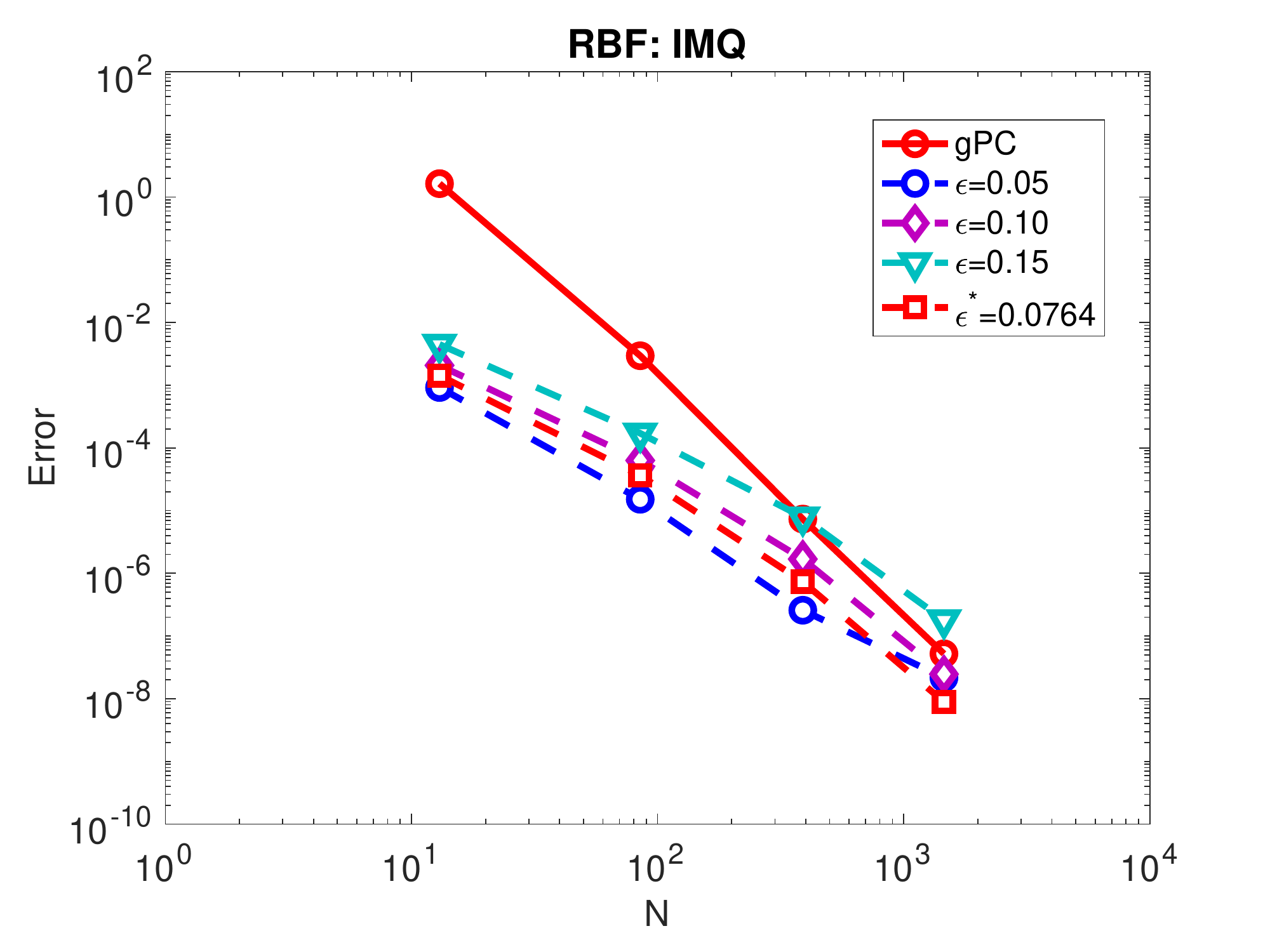}
    \includegraphics[width=4.4cm]{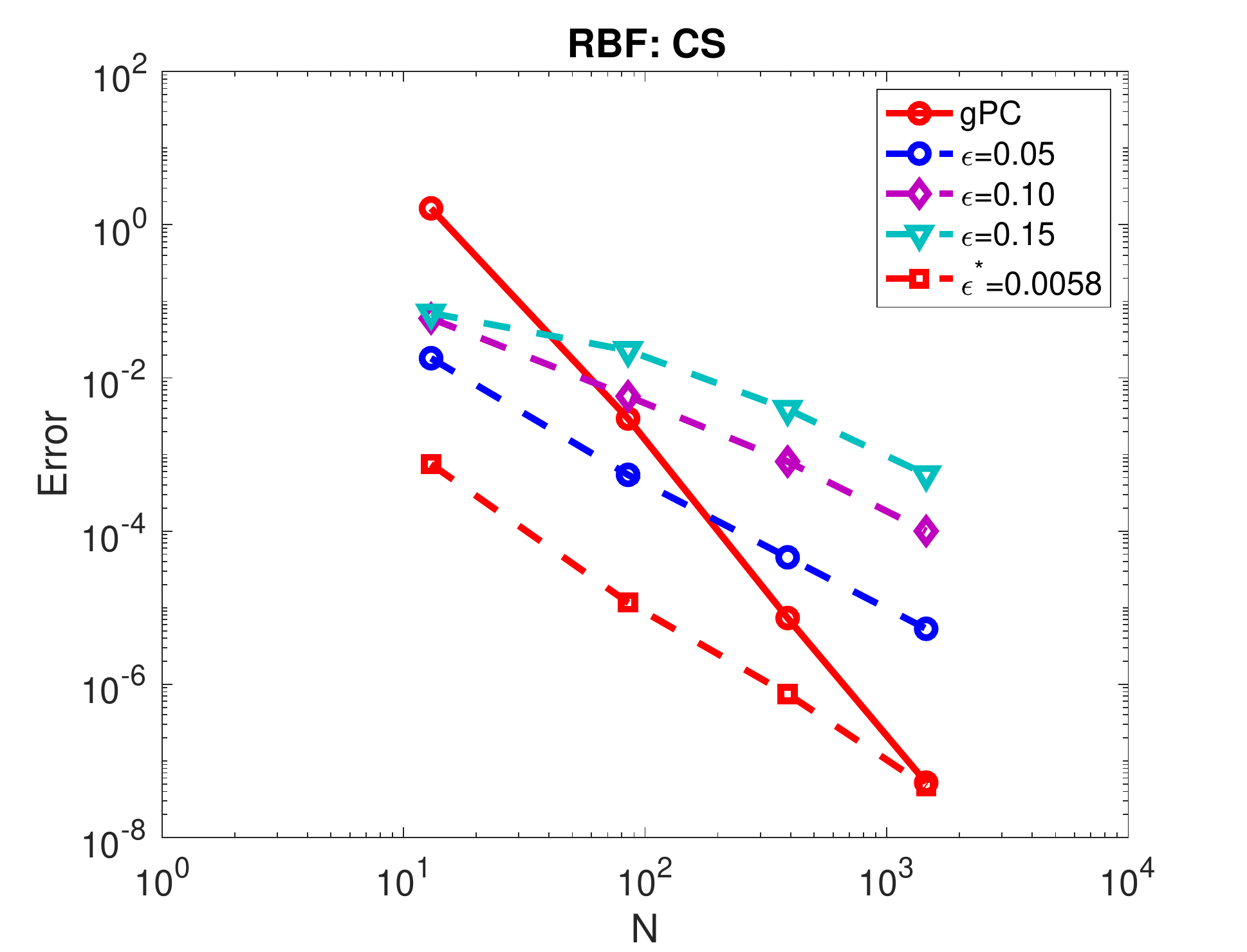}
\end{center}
  \caption{ Approximation error against the number of sample points $N$ (d=6). Left: Gaussian;  Middle: IMQ; Right: CS.
    \label{fig:d6_diff}
  }
  \end{figure}

Here we approximate the solution $u(z)=u(0.5,z)$, while the uncertain inputs $z^{(i)} \sim U[-1, 1], k=1,\cdots, d$.
For further comparison, we also employ the sparse grids stochastic collocation method using Legendre polynomial of total order $k=8 (d=3)$ and $k=4 (d=6)$ to solve the problem. The numerical results are presented in Figs. \ref{fig:d3_diff} and  \ref{fig:d6_diff}. The RBFs approximation methods,  are notably superior to the sparse grids method.


\subsection{Gradient-enhanced approach for UQ}

We now illustrate the symmetric approach to Hermite interpolation with a set of numerical experiments.

\textbf{Example 1: }  We use the corner peak
$$u(z) = (1+\sum^d_{i=1}\omega_i z^{(i)})^{-(d+1)}, \, z\in \Omega =[0, 1]^d$$ to generate the data.  In this example, we  use $d=2, \omega_i=\frac{1}{i^2}$.

\begin{figure}[htbp]
\begin{center}
    \includegraphics[width=8cm]{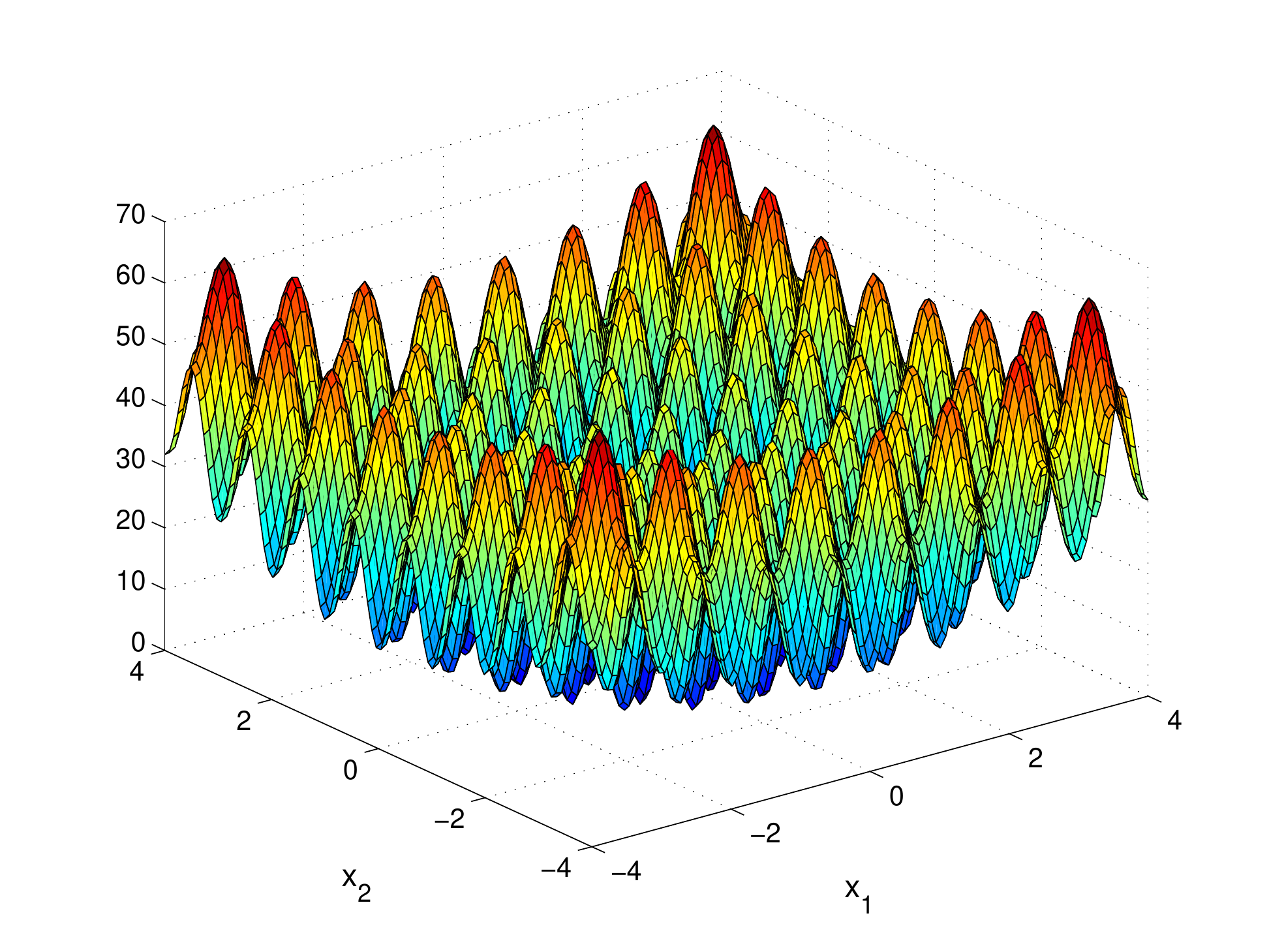}
\end{center}
  \caption{2-dimensional Rastrigin function.
  \label{fig:frasfun}
    }
\end{figure}

\textbf{Example 2:}  We consider the following 2-dimensional Rastrigin function:
$$u(z)= 20+\sum^2_{i=1}\Big((z^{(i)})^2-10 \cos(2\pi z^{(i)})\Big), \quad z^{(i)}\in [-4,4].$$
The 2-D Rastrigin function is plotted in Fig. \ref{fig:frasfun}.

\textbf{Example 3:} Consider the following 5-dimensional Friedman function
$$u(z)=10\sin(\pi z^{(1)} z^{(2)} )+20(z^{(3)}-0.5)^2+10 z^{(4)}+5 z^{(5)}.$$

\begin{figure}[htbp]
\begin{center}
    \includegraphics[width=6cm]{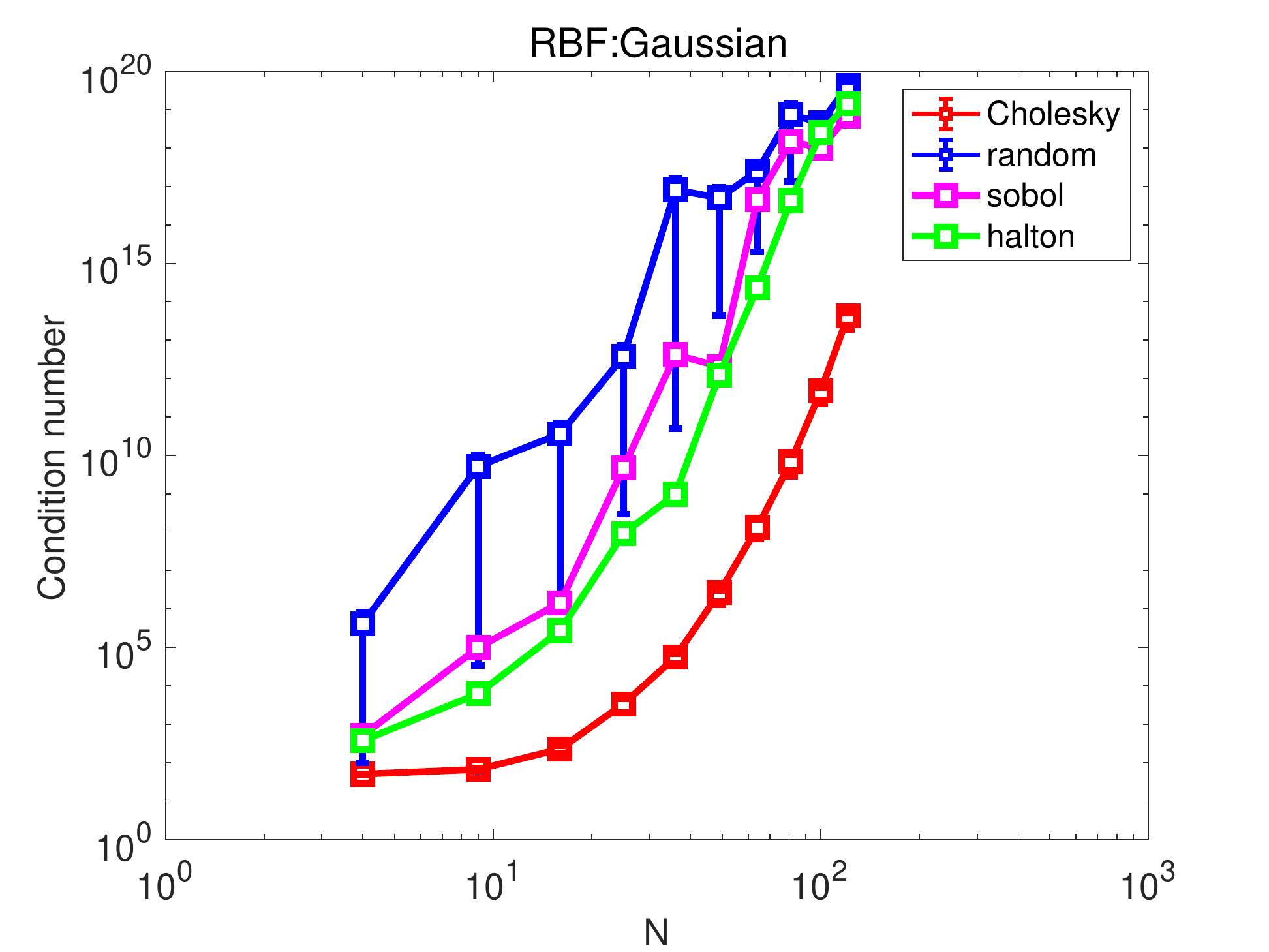}
    \includegraphics[width=6cm]{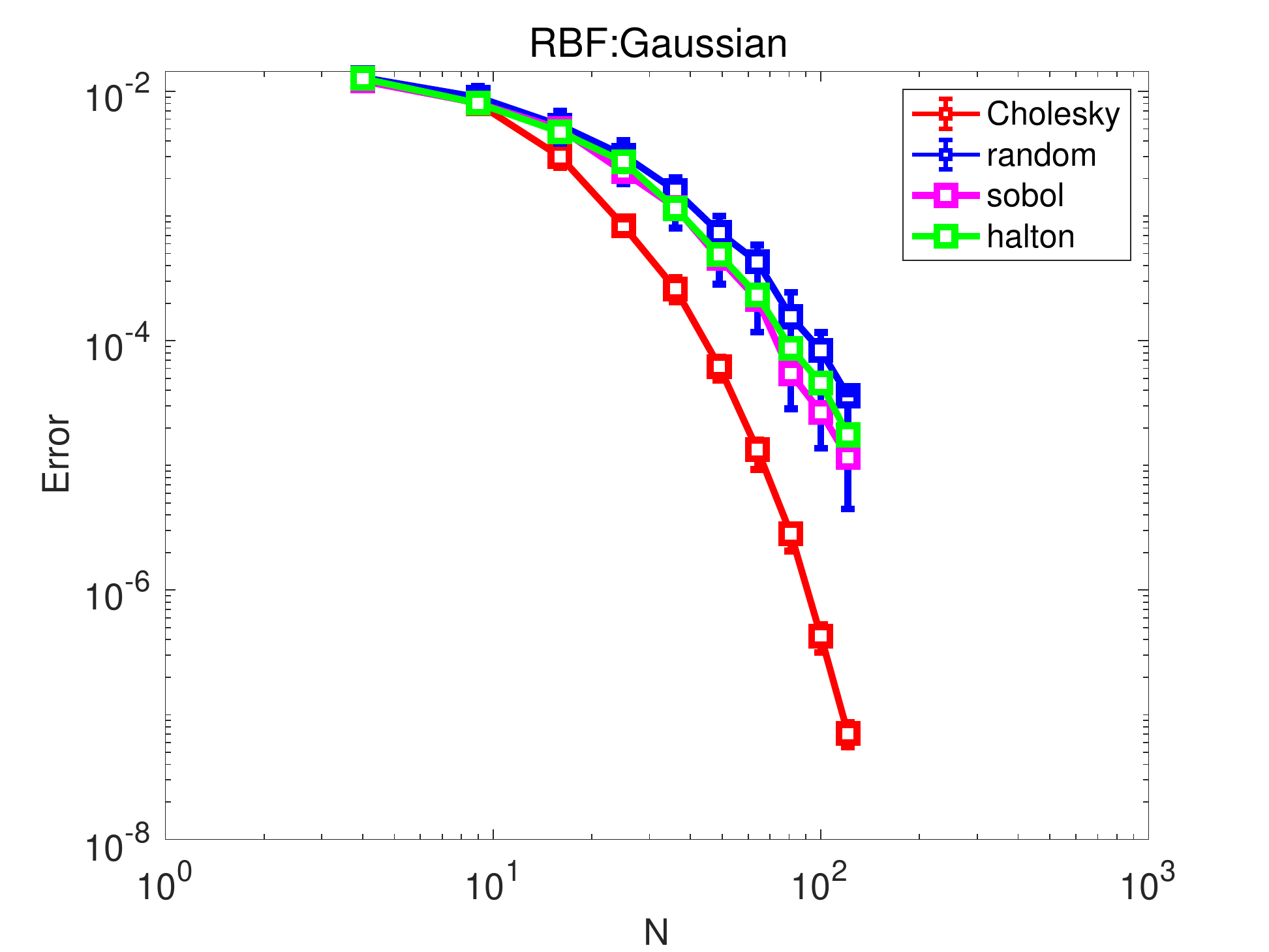}
     \includegraphics[width=6cm]{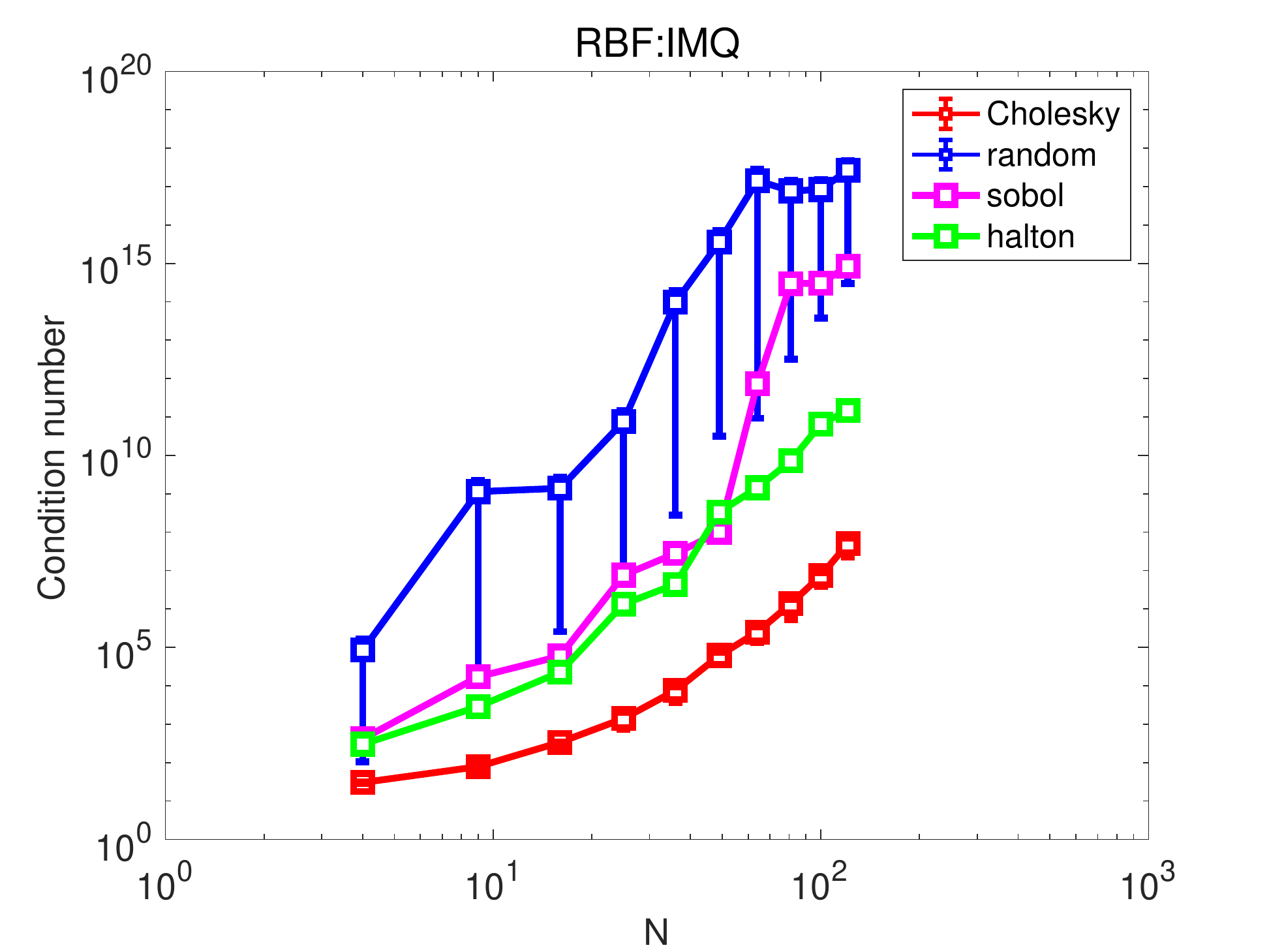}
    \includegraphics[width=6cm]{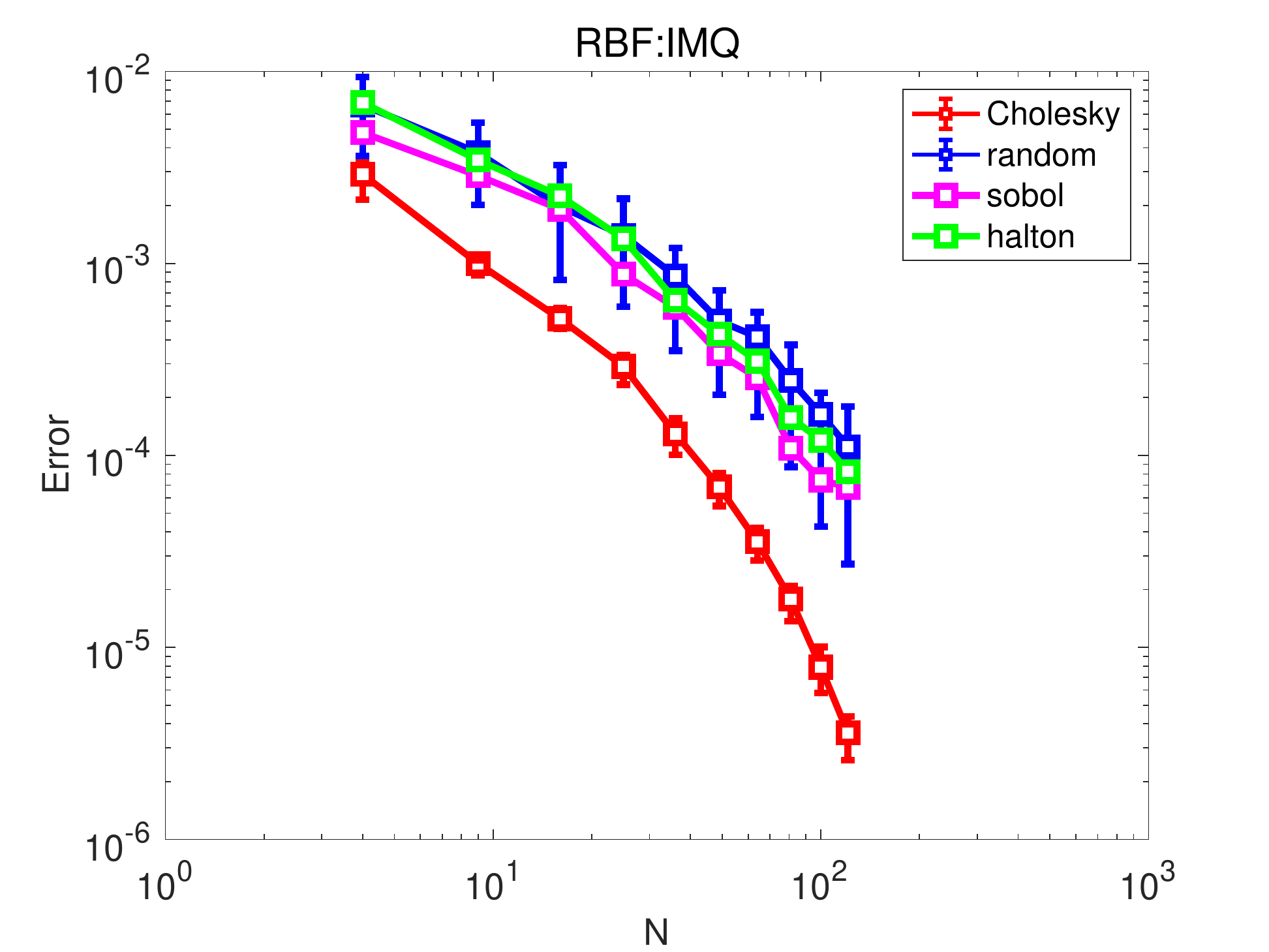}
      \includegraphics[width=6cm]{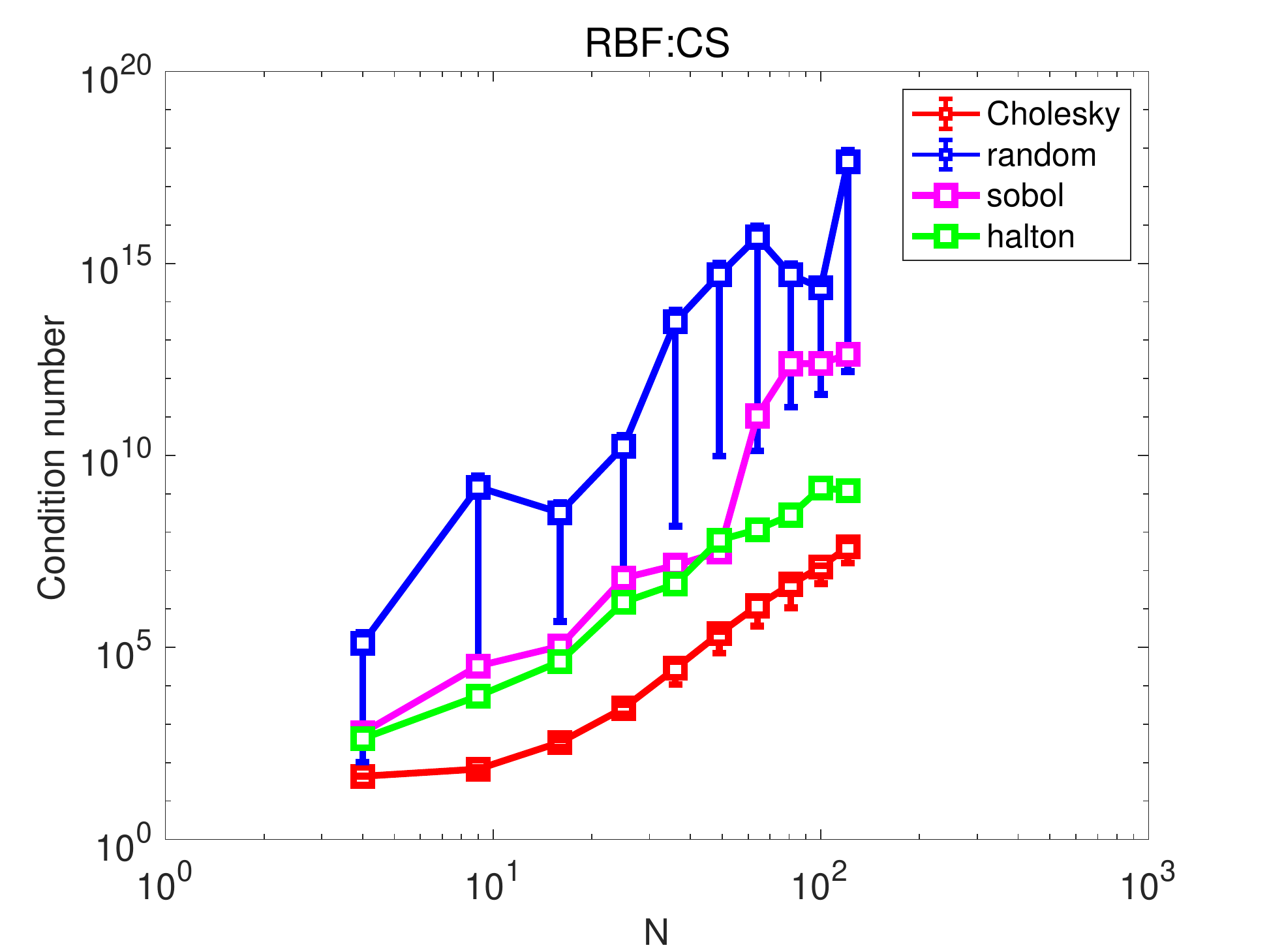}
    \includegraphics[width=6cm]{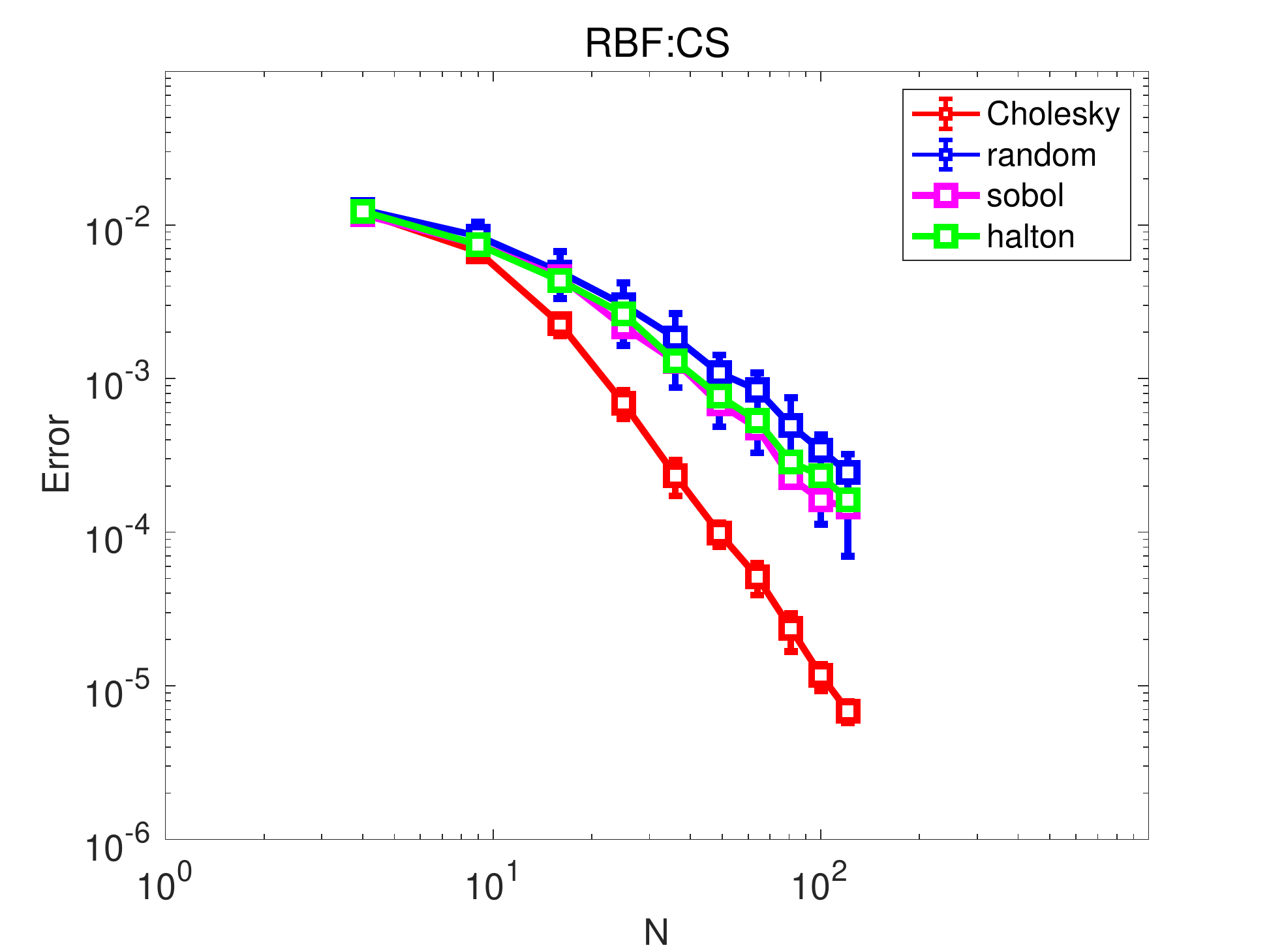}
\end{center}
  \caption{Numerical results with respect to the number of sample points $N$ for 2-dimensional corner peak.
  \label{fig:corner_2d}
    }
\end{figure}

\begin{figure}[htbp]
\begin{center}
       \includegraphics[width=6cm]{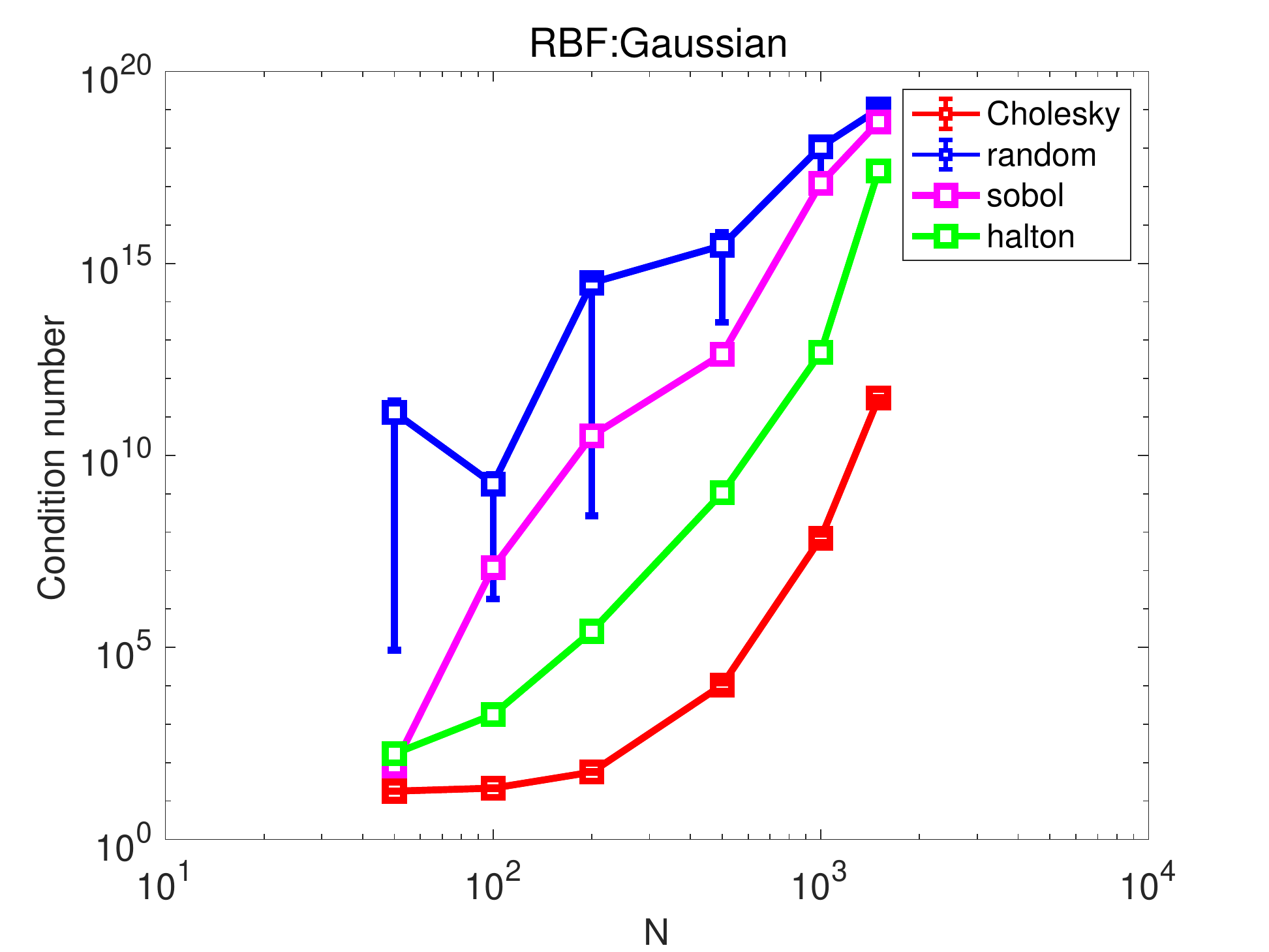}
    \includegraphics[width=6cm]{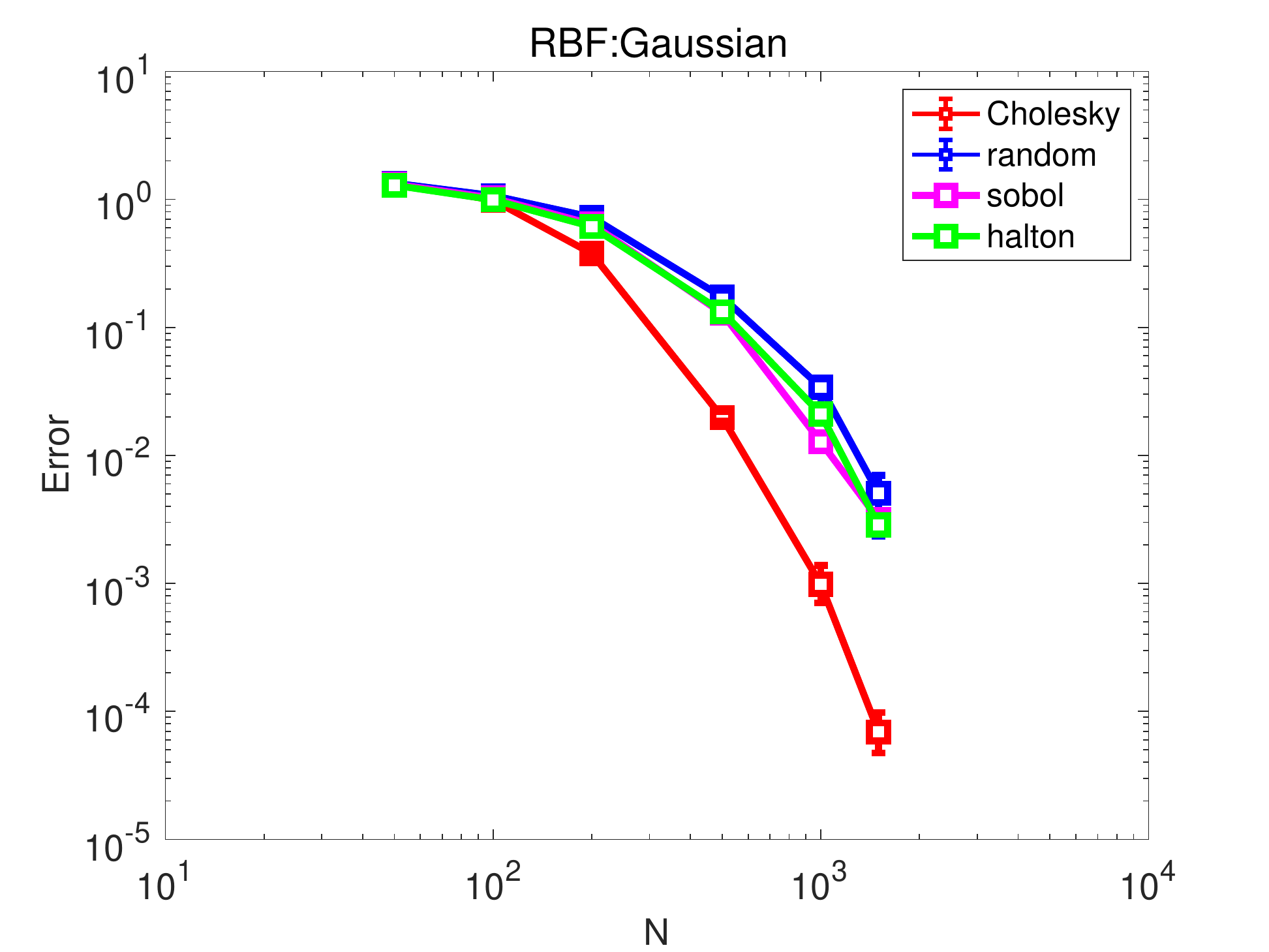}
     \includegraphics[width=6cm]{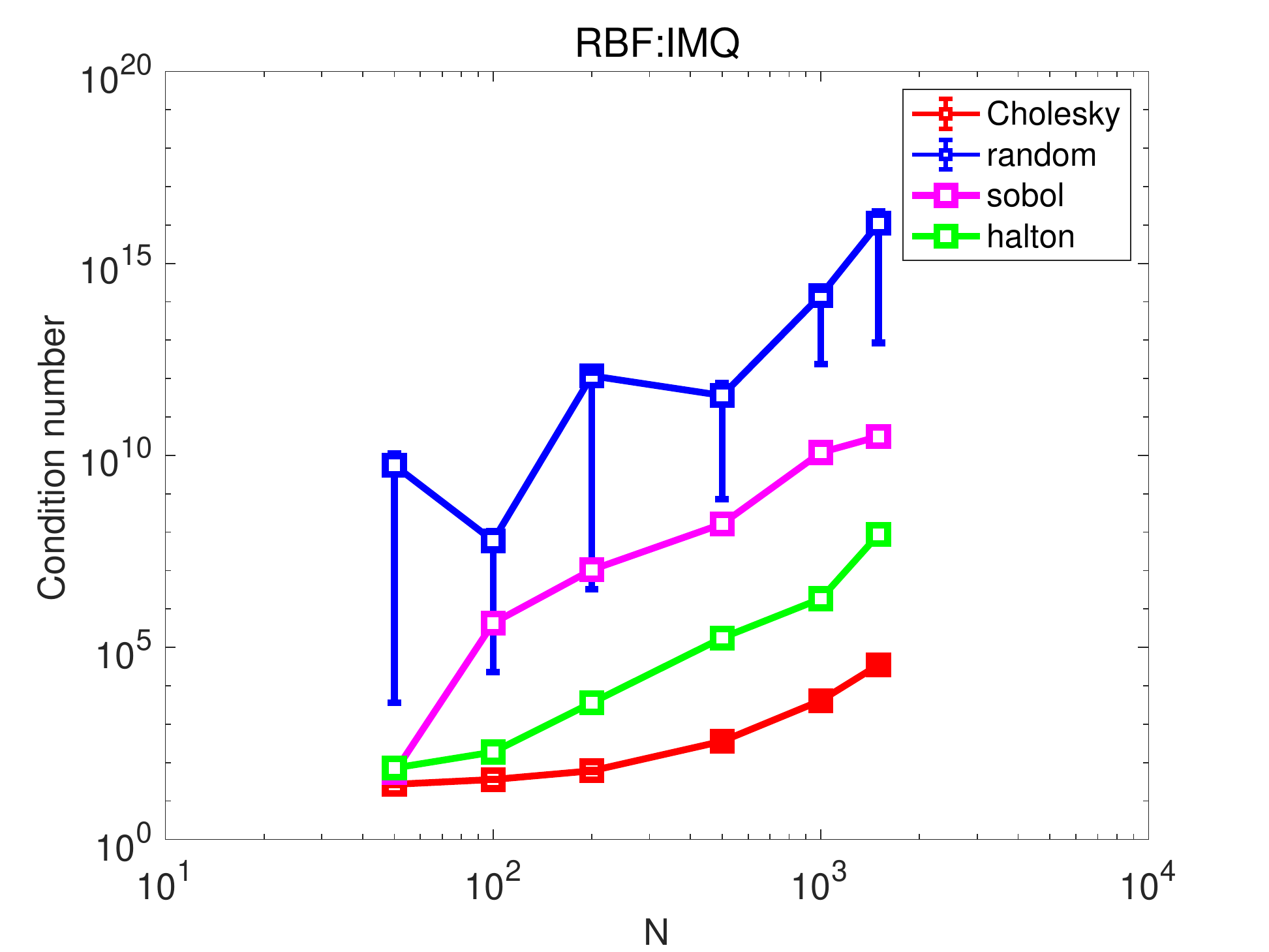}
    \includegraphics[width=6cm]{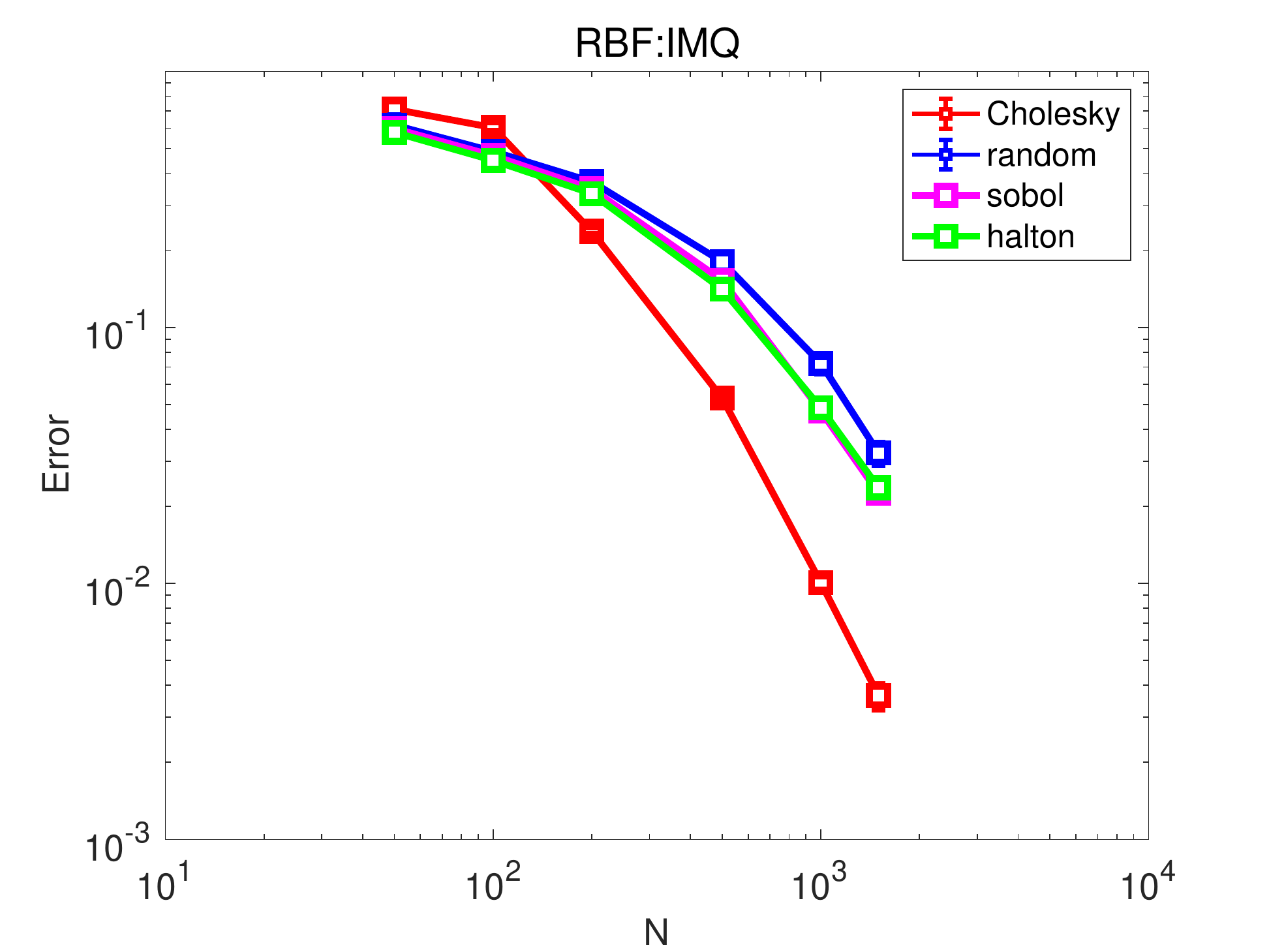}
      \includegraphics[width=6cm]{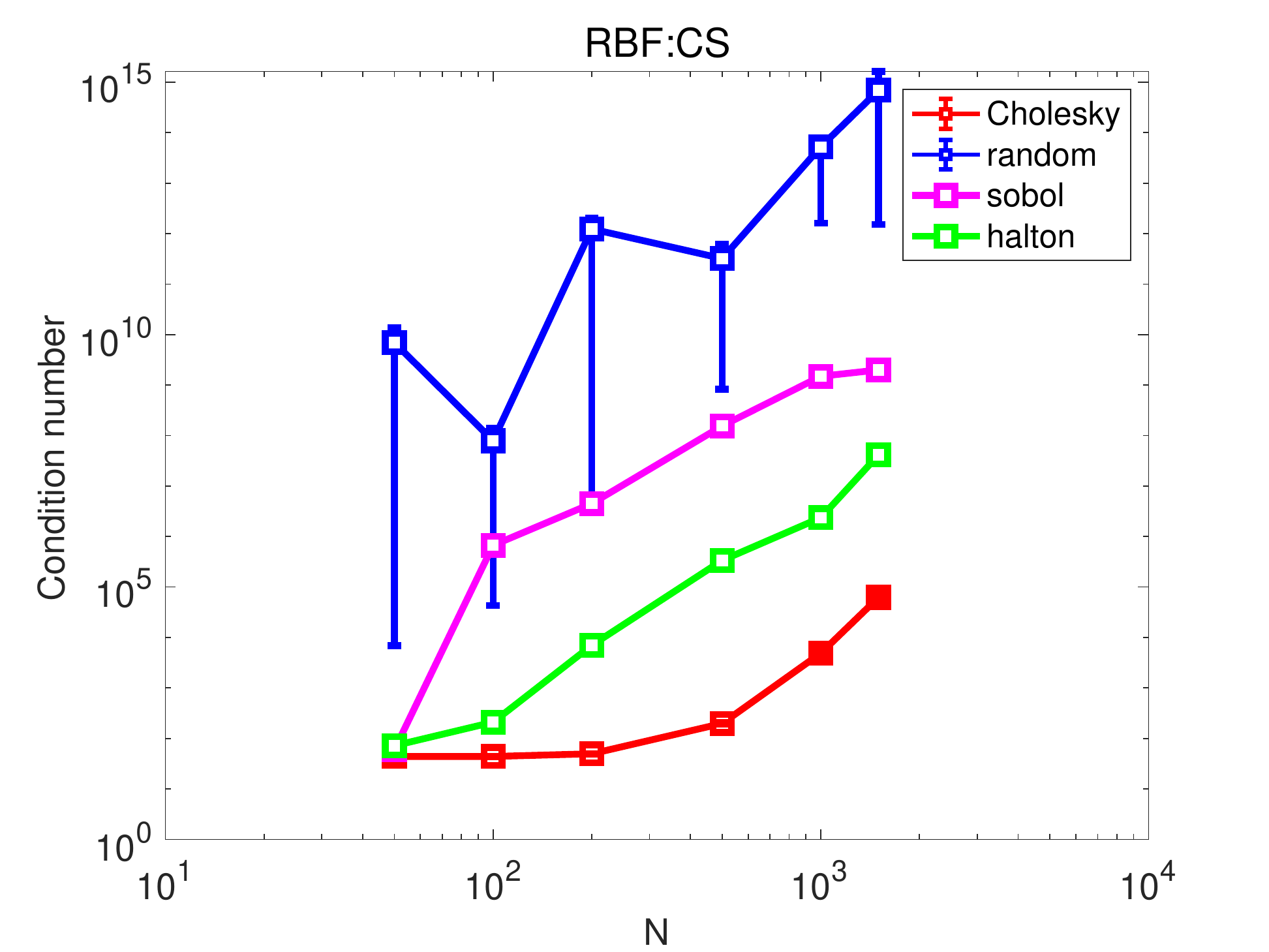}
    \includegraphics[width=6cm]{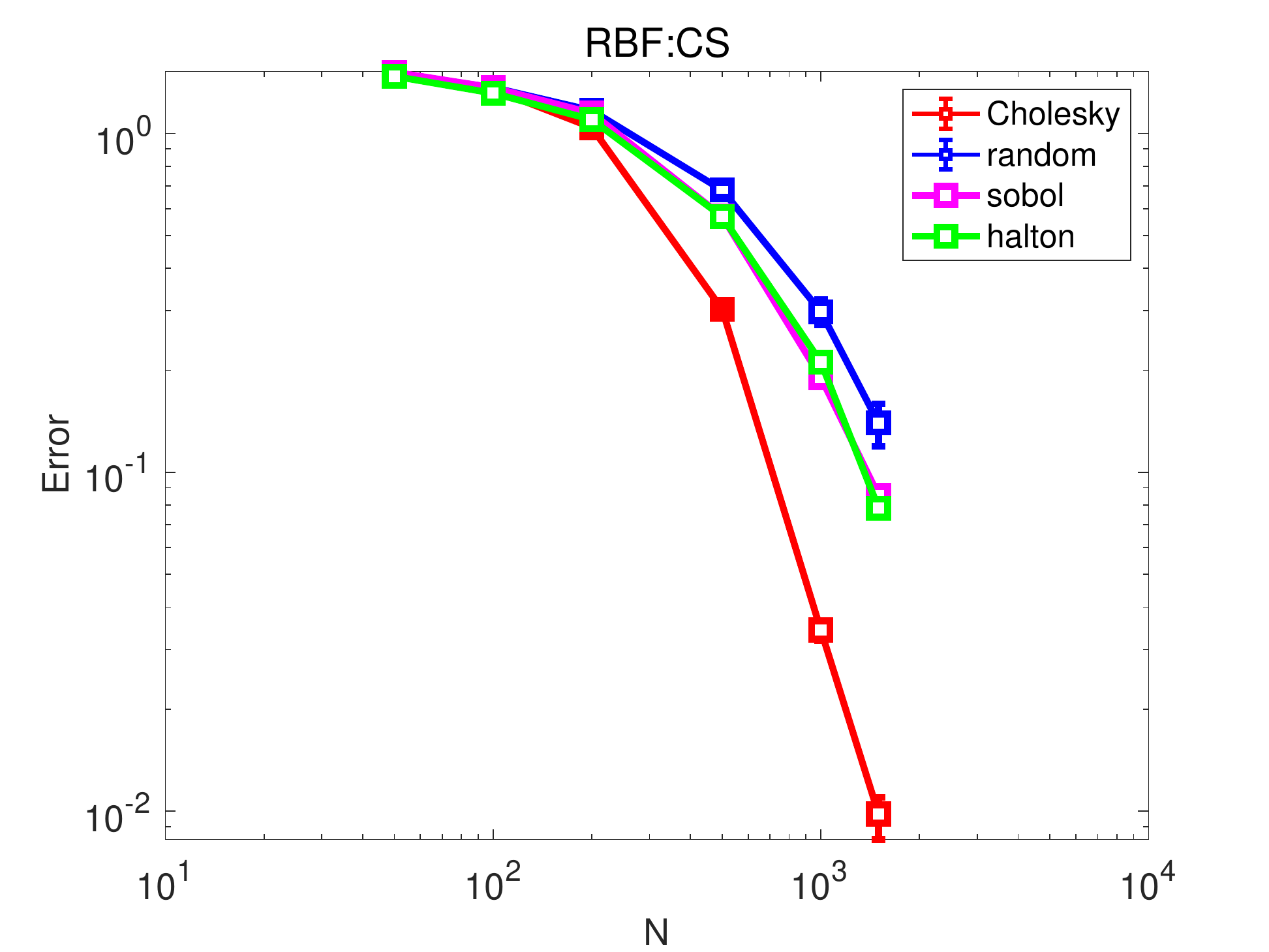}
\end{center}
  \caption{Numerical results with respect to the number of sample points $N$ for 2-dimensional Rastrigin function.
  \label{fig:fras2d}
    }
\end{figure}

\begin{figure}[htbp]
\begin{center}
    \includegraphics[width=6cm]{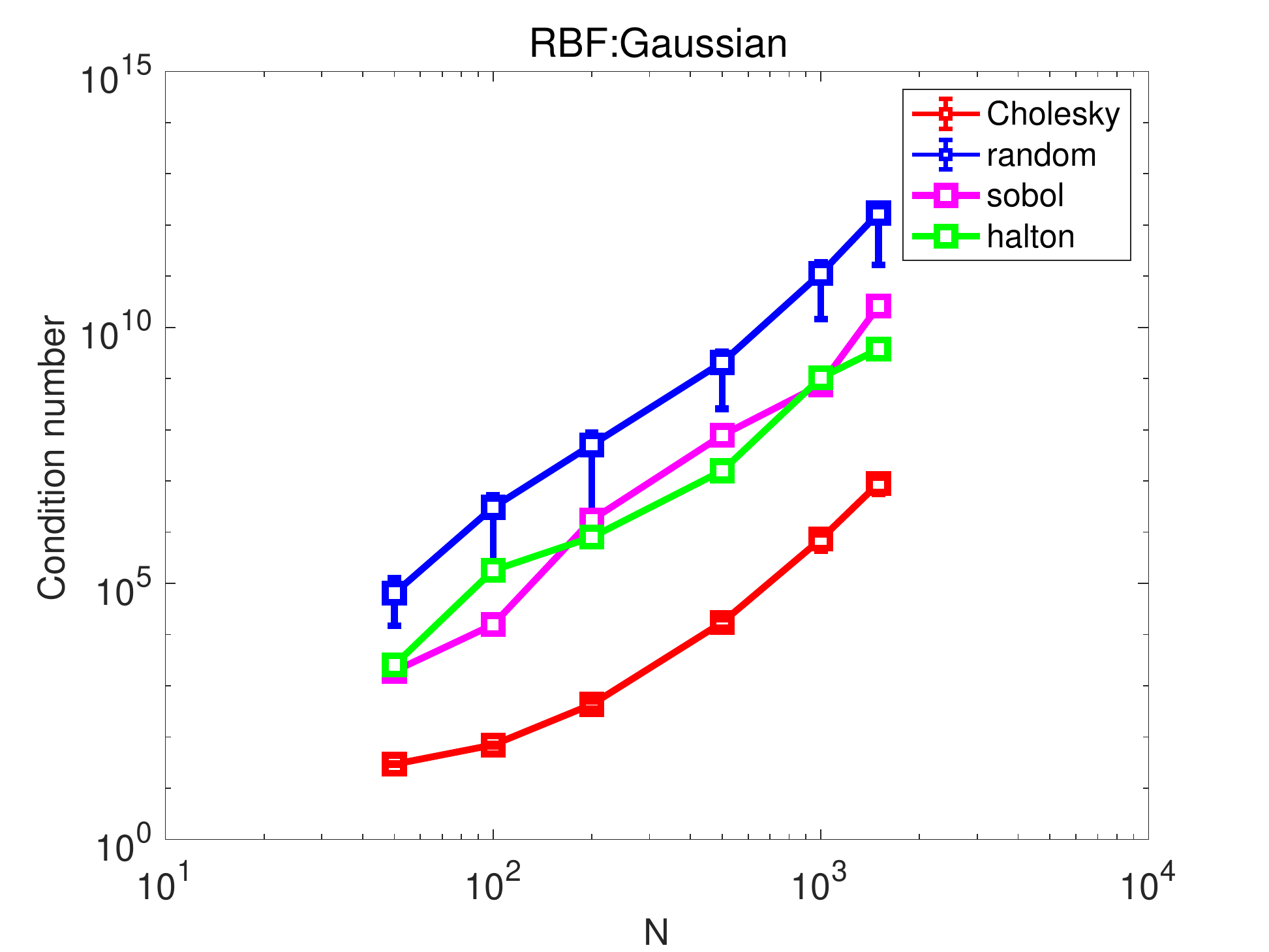}
    \includegraphics[width=6cm]{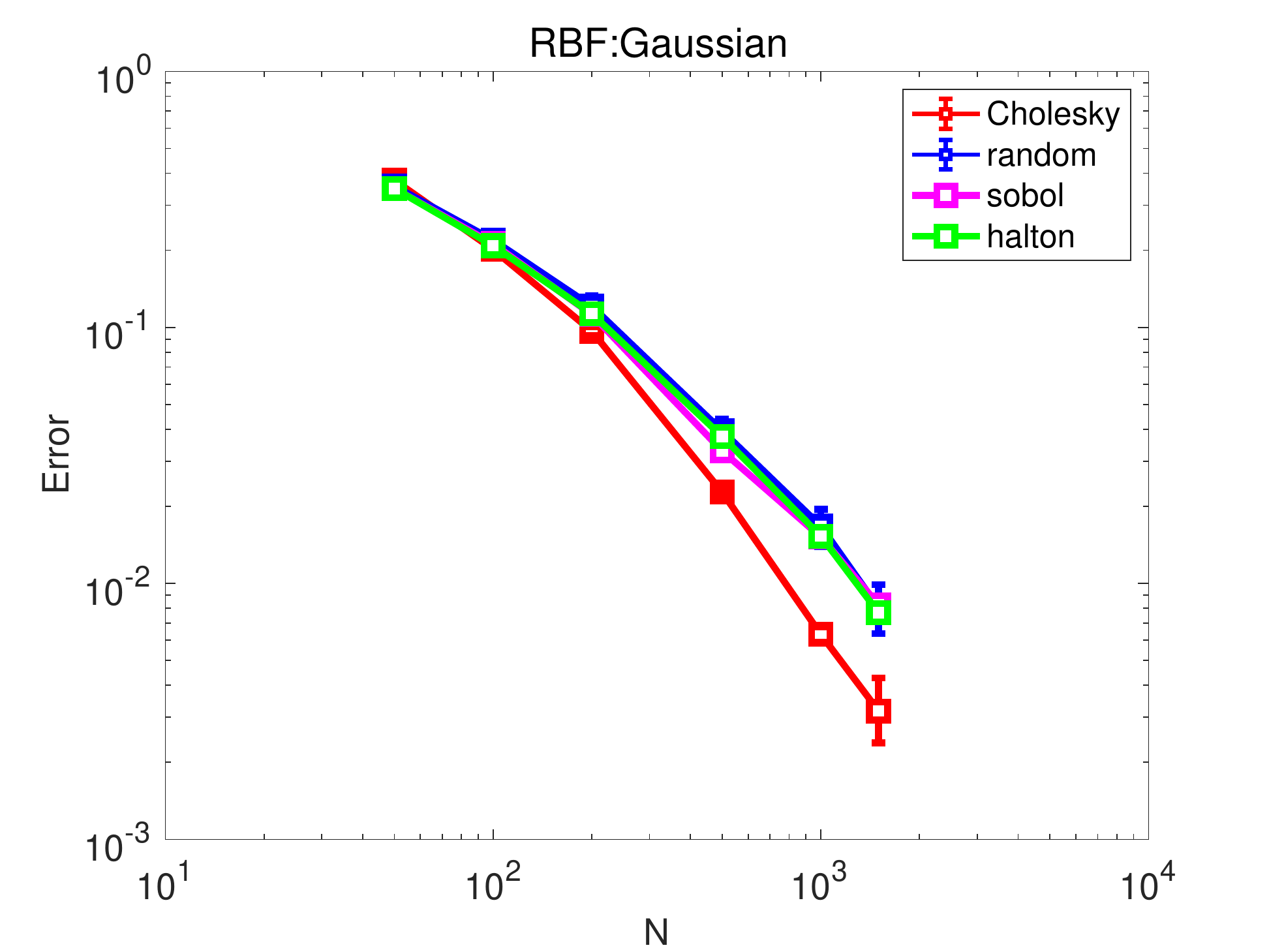}
     \includegraphics[width=6cm]{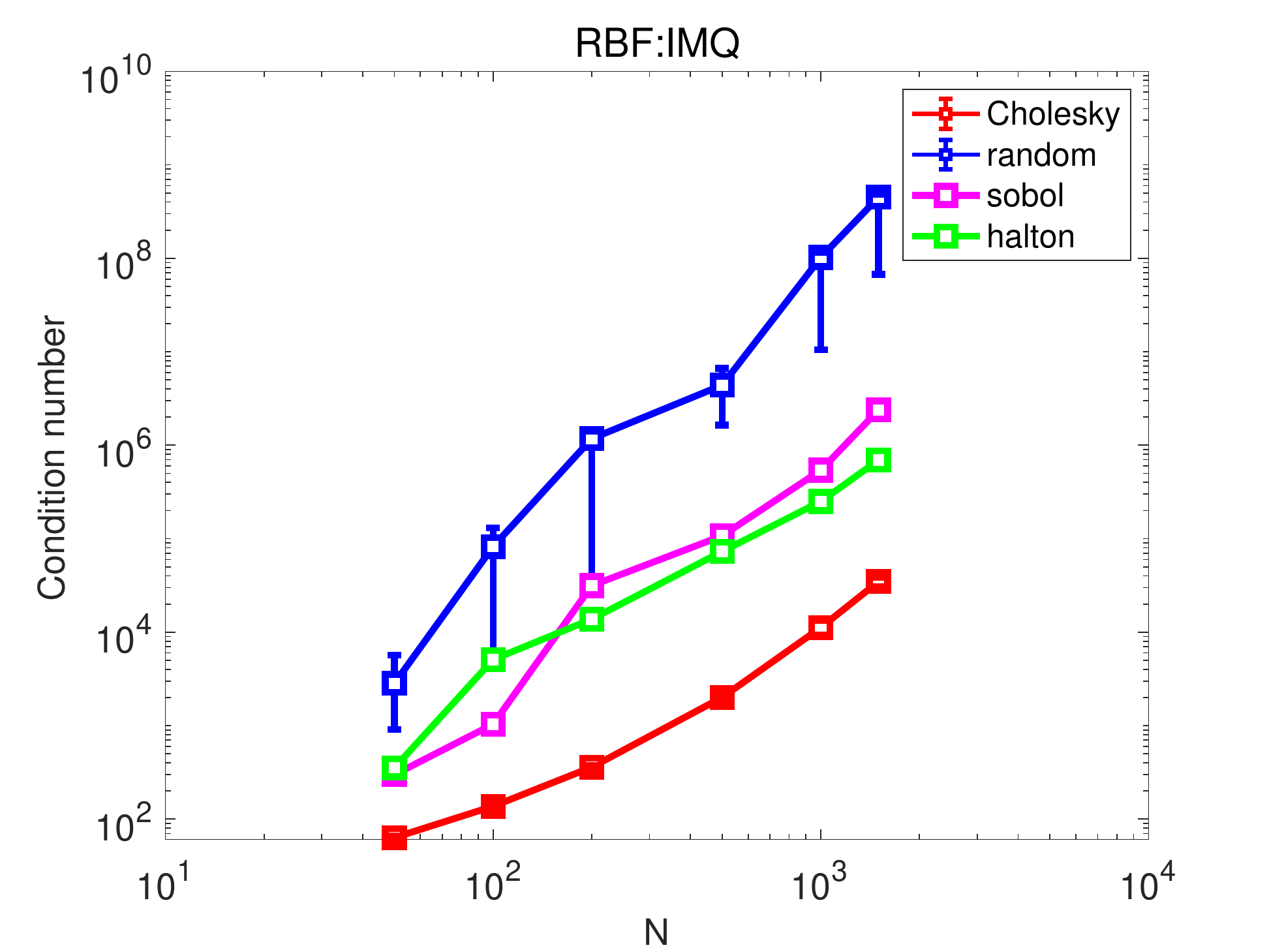}
    \includegraphics[width=6cm]{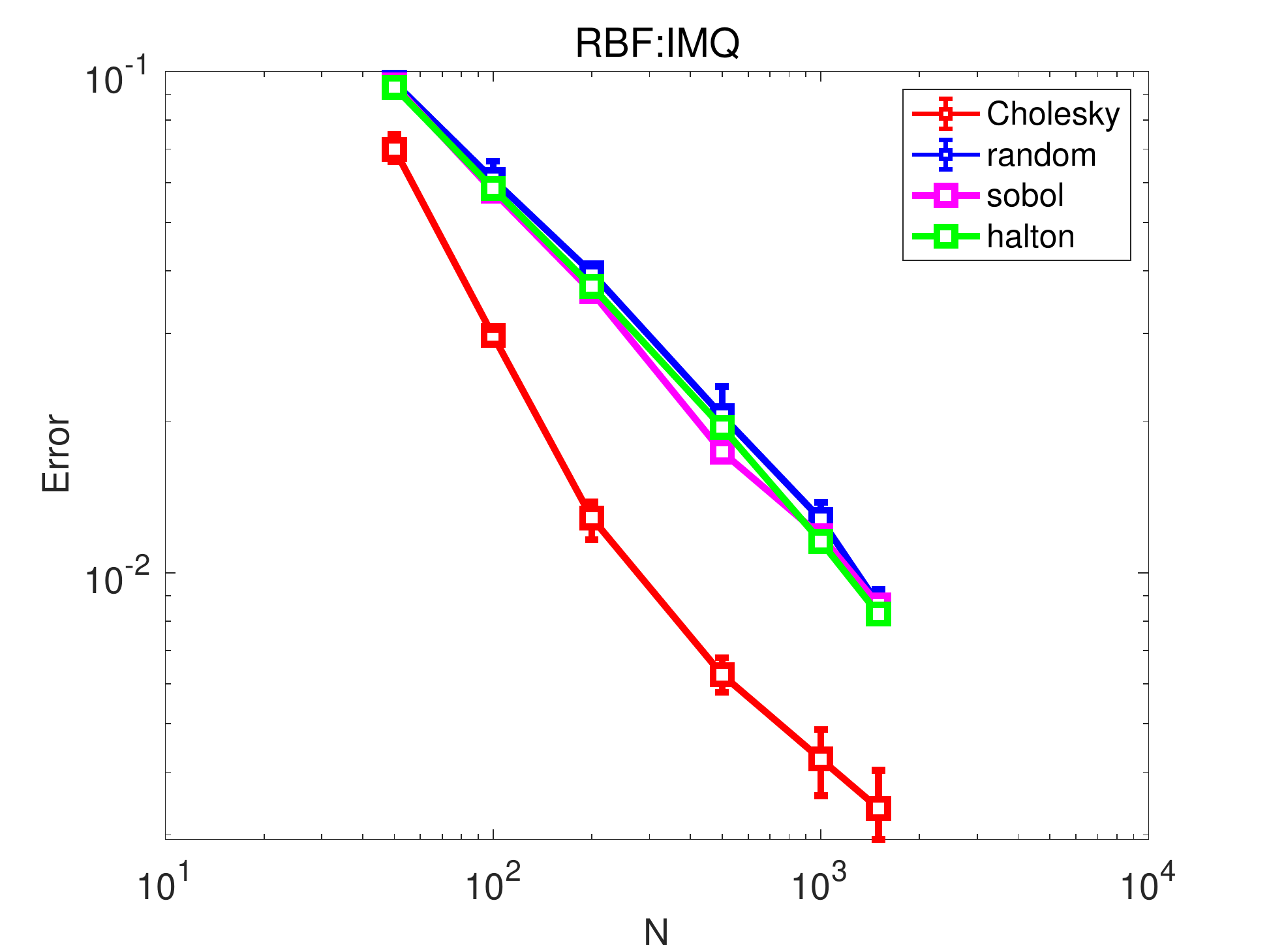}
      \includegraphics[width=6cm]{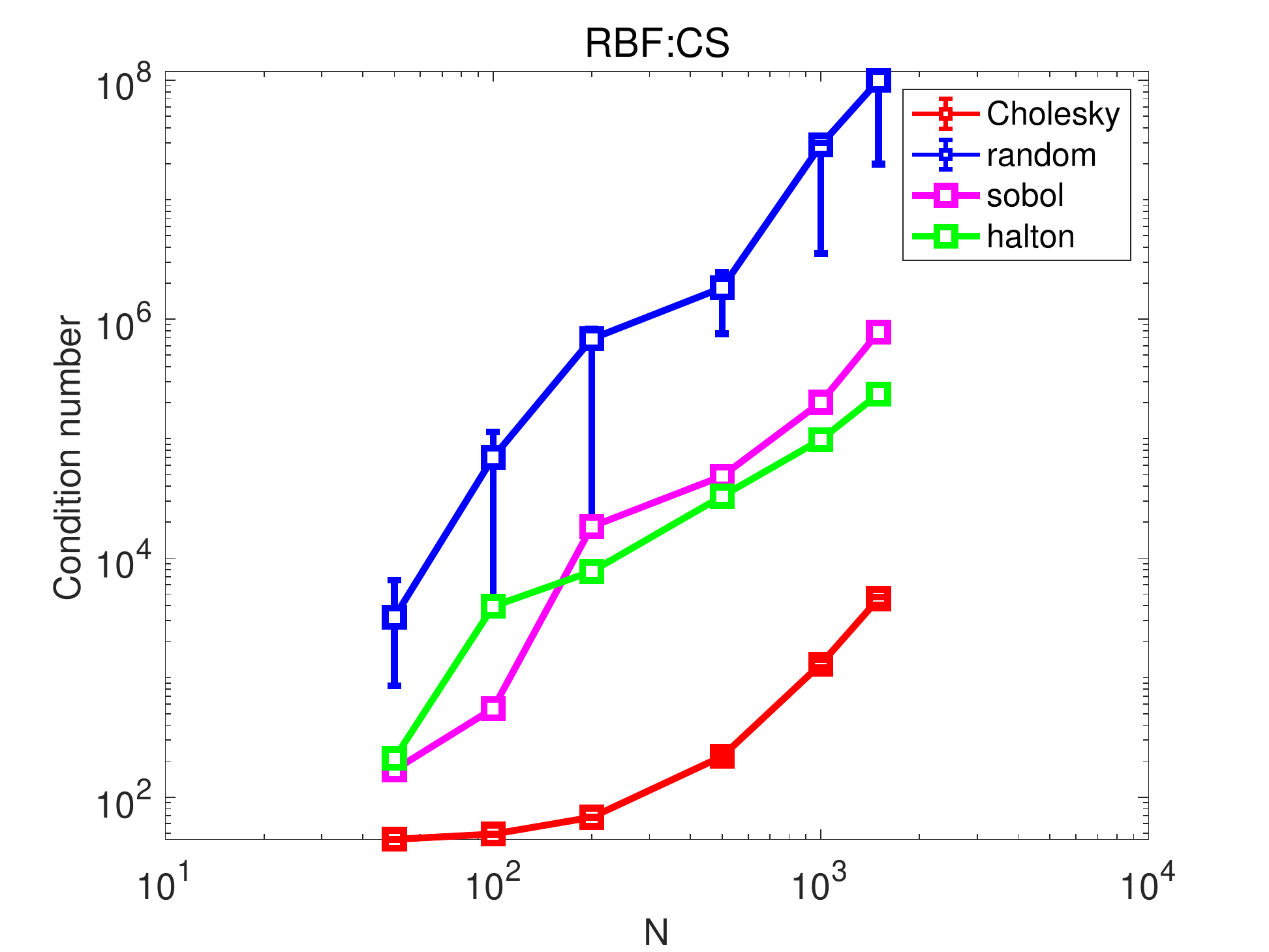}
    \includegraphics[width=6cm]{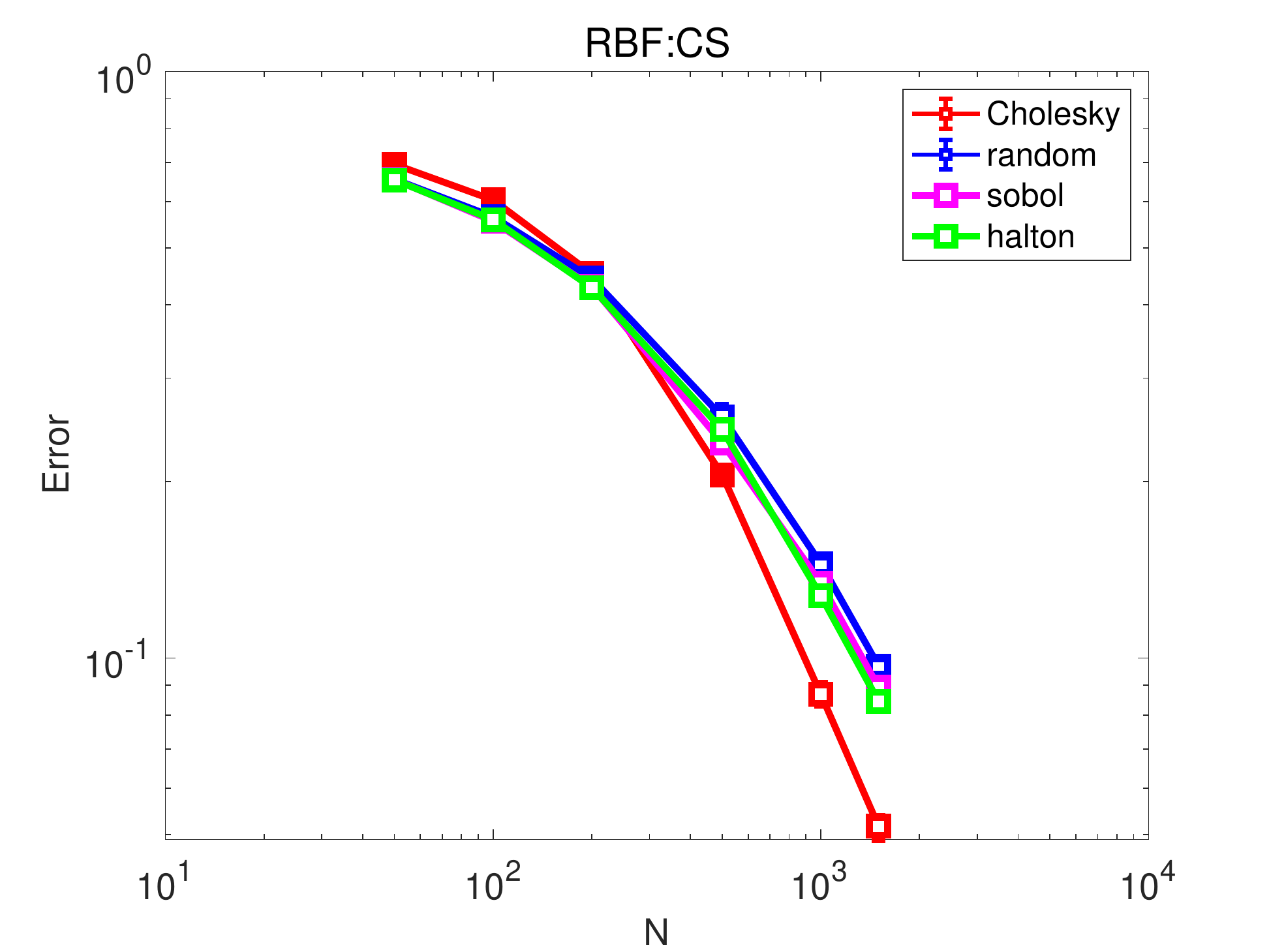}
\end{center}
  \caption{Numerical results with respect to the number of sample points $N$ for 5-dimensional Friedman function.  \label{fig:ffr5d}
    }
\end{figure}

The corresponding numerical results are shown in Figs. \ref{fig:corner_2d} - \ref{fig:ffr5d}. We can conclude that the design matrix $\mb{B}$ can be well conditioned under the proposed algorithm. Again,  our proposed algorithm  produces superior results than the other sampling methods.

\section{Summary} \label{sec:summary}

In this work, we presented a strategy for selecting a quasi-optimal sample points for kernel interpolation. We first demonstrate that the traditional sampling methods suffer from the numerical instability in the sense that the condition number will increase fast as  the total number of points becomes large.  To improve this, we propose to use the efficient Cholesky decomposition with pivoting for the Vandermonde-like interpolation matrix to choose the sample points. Then the quasi-optimal sample points are used to control the condition number of the design matrices. It is demonstrated that with the new approach the stability can be much improved. On the other hand, for problems inclusion of gradient measurements in the kernel interpolation are also considered.   Applications to parametric UQ problems are illustrated.

\appendix
\section{Gaussian process and kernel interpolation} \label{app: GP}

In order to help readers get deeper insights into the connections between kernel methods and Gaussian process regressions, we give a simple example in this section. More detains can be found in \cite{Scheuerer+Schaback2013}.

Suppose that we have a set of observed data values $\{(z_j,u_j)_{j=1}^N| z_j \in \Xi\}$. We aim to estimate the unknown value at some other locations $Z$ based on the observed data at $\Xi$. In the Gaussian process regression \cite{Rasmussen2006,Scheuerer+Schaback2013}, a.k.a the simple kriging method, $u$ is  assumed as a realization of a random field $\mb{Y}$, which is a collection $\{\mb{Y}(Z): Z\in I_Z\}$ of random variables over the probability space $(\Omega, \mathcal{F}, \mathcal{P})$. The observations $\{u(z_j)\}$ are then realizations of the random variables $\mb{Y}(z_j)$. To predict $\mb{Y}$ at some location $Z\in I_Z$, one considers all linear predictors of the form
\begin{eqnarray*}
\mb{Y}(Z) = \sum^N_{i=1}\omega_i(Z)\mb{Y}(z_i)
\end{eqnarray*}
which are themselves random variables.

To determine optimal weights $\mb{\omega(Z)}:=(\omega_1(Z),\cdots, \omega_N(Z))^T$, additional structural assumptions on $\mb{Y}$ are needed. In this work, we  only discuss the simple  kriging, due to their close connection to kernel interpolation.
Suppose $\mb{Y}$ is a Gaussian process with a mean 0 and symmetric covariance kernel $K$, i.e. $\mathbb{E}(\mb{Y}(Z))=0$,  $Cov(\mb{Y}(z_i),\mb{Y}(z_j))=K(z_i,z_j):=\phi(z_i-z_j)$. We usually assume that such covariance kernel $K$ is positive definite. The simple kriging method provides the best linear unbiased estimator $\hat{f}$ of $\mb{Y}$ in the form of
$$\hat{f}: = \omega_1(Z)u_1+\cdots+\omega_N(Z)u_N$$
where the weighting $\mb{\omega(Z)}$ is uniquely determined by the linear system (\cite{Scheuerer+Schaback2013})
\begin{eqnarray*}
\begin{bmatrix}
K(z_1,z_1) & \cdots & K(z_1,z_N) \\
\vdots &  \ddots & \vdots\\
K(z_N,z_1) & \cdots & K(z_N,z_N)
\end{bmatrix}
\begin{bmatrix} \omega_1(Z)\\ \vdots \\ \omega_N(Z) \end{bmatrix}
=\begin{bmatrix} K(Z,z_1)\\ \vdots \\ K(Z,z_N) \end{bmatrix}
\end{eqnarray*}
or in matrix form $K(\Xi,\Xi) \mb{\omega(Z)}= K(Z,\Xi)^T$. Thus, the estimator $\hat{f}(Z)$ can be written as
\begin{eqnarray}\label{gpeq}
\hat{f}(Z): = \mb{\omega(Z) u}
\end{eqnarray}

Note that if the RBF kernel and the covariance kernel coincide, then (\ref{interplant}) and (\ref{gpeq}) suggest that
$u_N(Z)=K(Z,\Xi)^TK(\Xi,\Xi)^{-1} \mb{u}=(K(\Xi,\Xi) ^{-1}K(Z,\Xi))^T\mb{u}=\hat{f}(Z)$. It means that the RBF estimator and the unbiased estimator are indeed identical.

\section{Gradient-enhanced Gaussian process}\label{app: ge-GP}

If the gradient information $\{z_i, u'_m(z_i)\}_{i=1}^N \, (m=1,2,\ldots, d)$ are included in the ordinary GP model, the covariance matrix becomes  a block matrix that includes the covariance between derivative observations in addition to the covariance between function observations. The block matrix can be represented  as \cite{Lockwood2012}:

\begin{eqnarray*}
\mb{K}=
\begin{bmatrix}
Cov(\mb{Y},\mb{Y}) & Cov(\mb{Y}, \nabla\mb{Y}) \\
Cov(\nabla\mb{Y},\mb{Y}) & Cov(\nabla\mb{Y}, \nabla\mb{Y})
\end{bmatrix},
\end{eqnarray*}
where $Cov(\mb{Y},\mb{Y}) $ represents an $N\times N$ covariance matrix of the function values at the sample points, $Cov(\mb{Y}, \nabla\mb{Y})$ is an $Nd\times N$ cross covariance matrix between gradient components and the function values at the sample points, and $Cov(\nabla\mb{Y}, \nabla\mb{Y})$ is an $Nd\times Nd$ covariance matrix of the gradients. If $\mb{Y}$ is a Gaussian process with a mean 0 and symmetric covariance kernel $K$ (or $\phi$), then the joint covariance matrix of the outputs is given by

\begin{eqnarray*}
\mb{K}=
\begin{bmatrix}
\mb{P}_{0,0}& \mb{P}_{0,1}  & \cdots & \mb{P}_{0,d} \\
\mb{P}_{1,0} & \mb{P}_{1,1} & \cdots & \mb{P}_{1,d} \\
\vdots & \vdots & \ddots & \vdots\\
\mb{P}_{d,0}  & \mb{P}_{d,1} & \cdots & \mb{P}_{d,d} \\
\end{bmatrix}
\end{eqnarray*}
where $\mb{P}_{k,l}$ is the $N\times N$ matrix with the $(i,j)$th element
\begin{eqnarray}
\begin{cases}
  \phi(z_i-z_j),   &  k=l=0, \\
  -\phi'_l(z_i-z_j),  & k=0, l\neq 0,  \\
 \phi'_k(z_i-z_j),  & k\neq0, l=0,  \\
   -\phi''_{k,l}(z_i-z_j),  & k\neq0, l\neq 0,
\end{cases}
\end{eqnarray}
and $\mb{K}$ is a symmetric $(d+1)N \times (d+1)N$ matrix.

Note that if the RBF kernel $\Phi$ and the covariance kernel $\phi$ coincide, then  the joint covariance matrix $\mb{K}$ and the design matrix $\mb{B}$  are indeed identical.

\end{document}